\documentstyle[12pt,twoside]{article}
\input{amssym.def}
\input{amssym}
\newtheorem{theorem}{Theorem}
\newtheorem{proposition}{Proposition}
\newtheorem{lemma}{Lemma}
\newtheorem{corollary}{Corollary}
\newtheorem{definition}{Definition}
\newtheorem{remark}{Remark}
\renewcommand{\thetheorem}{\thesubsection.\arabic{theorem}}
\renewcommand{\theproposition}{\thesubsection.\arabic{proposition}}
\renewcommand{\thelemma}{\thesubsection.\arabic{lemma}}
\renewcommand{\thedefinition}{\thesubsection.\arabic{definition}}
\renewcommand{\thecorollary}{\thesubsection.\arabic{corollary}}
\renewcommand{\theremark}{\thesubsection.\arabic{remark}}
\renewcommand{\theequation}{\thesubsection.\arabic{equation}}
\renewcommand{\thesubsection}{\thesection.\arabic{subsection}}
\begin{document}
\def\R {{\Bbb R}}
\def\N {{\Bbb N}}
\def\C {{\Bbb C}}
\def\Z {{\Bbb Z}}
\def\T {{\Bbb T}}
\def\D {{\cal D}}
\def\E {{\cal E}}
\def\ipar {\par \qquad}
\def\apbp {|\alpha'+\beta'|}
\def\atbt {\tilde{\alpha}+\tilde{\beta}}
\def\xx{(x,\xi)}
\def\x0x0{(x_0,\xi_0)}
\def\yy{(y,\eta)}
\def\xxx{(x_0,\xi_0)}
\def\xd{(x,D_x)}
\def\yd{(y,D_y)}
\def\gs{$G^s$}
\def\gsp{$G^s$(pseudo -- )}
\def\ct {T^*{\Bbb R}^n\setminus\{0\}}
\def\lc {<\kern-2.4ex\lower1.3ex\hbox{$\sim$}}
\def\cc{\subset\!\subset}
\def\dt{d\kern -.30em\lower -.52em\hbox{\_}}
\def\ls{\lesssim}
\def\f{\varphi}
\def\ps{p_1^s\xx }
\def\tp{\hbox{{\rm Tr}}^+\,}
\def\sp{\hbox{{\rm sp}}\,}
\def\im{\hbox{{\rm Im}} \,}
\def\Re{\hbox{{\rm Re}} \,}
\def\tphi#1#2{(T^{#1})_{#2}}
\def\dzr#1{\|\|\Lambda^\deltaZ^{#1}Rv\|\|}
\def\zr#1{\|\|Z^{#1}Rv\|\|}
\def\normL2#1{\|{#1}\|_{L^2}}
\def\normomL2#1{\|{#1}\|_{L^2(\Omega)}}
\def\normdelta#1{\|{#1}\|^2_\delta}
\def\normtwodelta#1{\|{#1}\|^2_{2\delta}}
\def\norms#1#2{\|{#1}\|^2_{#2}}
\def\ip#1#2{({#1},{#2})_{L^2}}
\def\Lbar{\overline{L}}
\def\tphidef#1#2#3{\sum_{|\alpha+\beta|\leq{#1}}{{(-1)^{|\alpha|}}
\over{\alpha!\beta!}}X''^\alpha
X'^\beta
T^{{#1}-|\alpha+\beta|}\circ
ad^\alpha_{X'}ad^\beta_{X''}(#2#3)}
\def\tphiC2def#1#2#3{\sum_{a+b\leq{#1}}{{(-1)^a}
\over{a!b!}}Y^a X^b
T^{{#1}-(a+b)}\circ
	ad^a_{X}ad^b_{Y}(#2#3)}
\def\be{\begin{equation}}
\def\ee{\end{equation}}
\def\subJ{_{|J|\leqk_0}}
\def\ident {\equiv}
\def\czi{C_0^\infty}
\def\op{{\rm Op} }
\def\caution{ \marginpar{\att S} }
\title{Propagation of Gevrey Regularity for a Class of
Hypoelliptic Equations}
\author{Antonio Bove\\
Dipartimento di Matematica\\
Universit\`a di Bologna\\
40127 Bologna, Italy\\ \&
\\David S. Tartakoff\\ Department
of Mathematics\\University of Illinois at
Chicago\\ 851 S. Morgan, m/c 249\\
 Chicago Illinois 60607-7045, U.S.A.}
\date{
}
\maketitle
\begin{abstract}
We prove results on the propagation of Gevrey and analytic wave
front sets
for a class of $C^\infty$ hypoelliptic equations with double
characteristics.
\end{abstract}
\section{Introduction}
\setcounter{equation}{0}
\setcounter{theorem}{0}
\setcounter{proposition}{0}
\setcounter{lemma}{0}
\setcounter{corollary}{0}
\setcounter{definition}{0}
\setcounter{remark}{0}

It is well known that a (pseudo -- differential) operator with semi --
definite principal part and at most double characteristics is not, in
general, $C^\infty$ hypoelliptic; however if the
lower order terms satisfy some supplementary assumptions then there is
$C^\infty$ hypoellipticity (see e.g. \cite{hormanderbook} vol. III,
\cite{hormander1}).

As far as analytic hypoellipticity is concerned, the situation is more
involved. There are in fact examples of operators being $C^\infty$
hypoelliptic (i.e. whose lower order terms satisfy the ``Levi"
conditions) but {\em not \/} analytic hypoelliptic and, at the same time,
it has been proved that if the operator's principal part vanishes exactly
of order 2 on a manifold $\Sigma$ in the cotangent bundle, if its lower
order terms satisfy the $C^\infty$ hypoellipticity conditions and
if $\Sigma$  is symplectic (i.e. the symplectic form $\sigma = d\xi \wedge
dx $ has  maximal rank on $T\Sigma$) then there is analytic hypoellipticity
(see  e.g. \cite{metivier 1981}, \cite{tartakof80}, \cite{treves},
\cite{tartakofNA}).

The situation becomes more involved if the symplectic form has not maximal rank
or if it can degenerate on a submanifold (subset) of the double characteristic
manifold. It has been proved by M\'etivier (\cite{metivier3-1982},
\cite{metivier4-1980}) that actually there is propagation of the analytic
singularities on the leaves of the characteristic manifold if the operator
satisfies the conditions for $C^\infty$ hypoellipticity with loss of one
derivative. Essentially M\'etivier constructs null
microlocal null solutions for
certain microlocal models having a non empty analytic wave front set. Moreover
analizing the same micolocal model he proved a theorem of propagation of the
analytic regularity.

A much deeper analysis has been carried out by Sj\"ostrand, \cite{sjostrand},
using F.B.I.S. transform with Lipschitz Lagrangiean manifold, considering also
the case in which the rank of the fundamental matrix of the pricipal symbol may
degenerate on a submanifold of the double characteristic manifold. In
particular he gave another proof of  M\'etivier's theorem on the propagation of
the regularity.

The purpose of the present paper is to give another proof of those microlocal
hypoellipticity results for operators with double characteristics
satisfying the conditions for $C^\infty$ hypoellipticity with loss of one
derivative. A first result essentially states that the Gevrey wave
front set $WF_s$  (or rather its complementary, i.e. the set of points of
$G^s$ - regularity), $s \geq 1$, propagates  along the leaves of the
characteristic manifold $\Sigma$ (in particular  for symplectic manifolds we
get analytic hypoellipticity). The second  result says that an operator in
the above mentioned class is actually
$G^s$ hypoelliptic if $s \geq 2$. This second result, using different
techniques, has been obtained by Kajitani -- Wakabayashi in \cite{kw}.

Our technique is to deduce {\em a priori\/} $G^s$ bounds moving from the
starting point of {\em a priori\/} hypoellipticity estimates in the
$C^\infty$ case. To do this we need a careful microlocalization procedure
in the directions tangent to $\Sigma$. This is accomplished by a
technique due to the second author, \cite{tartakof80}, \cite{tartakofNA},
and already used to give an alternative proof of M\'etivier's theorem
\cite{metivier 1981},  although the full details are still unpublished. We feel
that this technique
can be useful in more generality and in degenerate situations.

The first five sections of the paper are devoted to establishing the
notation, introducing the microlocalization and proving the first
theorem. The sixth section is concerned with the proof of the second
thorem. An appendix collects some general - purpose material used
throughout the paper.

Finally the first author would like to take this opportunity to thank the
Department of Mathematics of the University of Illinois at Chicago, where
he stayed for three weeks during the preparation of the final version of
this paper: this allowed him to enjoy lots of mathematics and to short
cut the clumsiness of e-mail!

\section{Preparations and Statement of Results}
\label{s0}
\setcounter{equation}{0}
\setcounter{theorem}{0}
\setcounter{proposition}{0}
\setcounter{lemma}{0}
\setcounter{corollary}{0}
\setcounter{definition}{0}
\setcounter{remark}{0}
\renewcommand{\thetheorem}{\thesection.\arabic{theorem}}
\renewcommand{\theproposition}{\thesection.\arabic{proposition}}
\renewcommand{\thelemma}{\thesection.\arabic{lemma}}
\renewcommand{\thedefinition}{\thesection.\arabic{definition}}
\renewcommand{\thecorollary}{\thesection.\arabic{corollary}}
\renewcommand{\theequation}{\thesection.\arabic{equation}}
\renewcommand{\theremark}{\thesection.\arabic{remark}}

Let $P\xd = P_m\xd + P_{m-1}\xd + \cdots$ be a classical $G^s$
(pseudo) -- differential operator of order $m$, $s \geq 1$, and denote by
$p_{m-j}$ the symbols of the $P_{m-j}$, which are
(positively)  homogeneous of degree $m - j$ with respect to $\xi$. We shall
make the  following assumptions:
\begin{description}
	\item[(H${}_1$)]
\begin{description}
	\item[(a)]  $p_m\xx \geq 0$, for every $\xx \in \ct$.
	
	\item[(b)]  Let $\Sigma = \{p_m\xx = 0\}$. Then $\Sigma$ is a real \gs
	manifold in $\ct$.
	
	\item[(c)]  $p_m$ vanishes on $\Sigma$ exactly of order 2, i.e. $p_m\xx
	\geq $ Const $d_{\Sigma}^2\xx$, where $d_{\Sigma}\xx$ denotes the
	distance of the point $\xx \in \ct$ from $\Sigma$.
	
	\item[(d)]  Let $F_{p_m}\x0x0$ denote the Hamilton map of $p_m$ at
	$\x0x0 \in \Sigma$, defined by
	\[
	\langle F_{p_m}\x0x0 t, d\varphi\rangle = \frac{1}{2}\langle t,
	d(H_{p_m}\varphi)\x0x0 \rangle,
	 \]
	 where $t \in T_{\x0x0}(\ct)$ and $\varphi$ is a smooth function. Then
	 \[
	 \dim \ker F_{p_m}\x0x0 = \dim  T_{\x0x0}\Sigma
	  \]
	  and
	  \[
	  \sigma_{|\Sigma} \quad {\rm has \ constant\  rank.}
	   \]
\end{description}
\end{description}
Without any loss of generality we shall suppose henceforth that $m = 2$,
the general case being recovered multiplying by an elliptic pseudo --
differential factor.

Denote by $\ps$ the subprincipal symbol of $P$, defined by
\[
\ps = p_1\xx + \frac{i}{2} \sum_{j=1}^{n}\frac{\partial^2 p_2\xx}{\partial
x_j
\partial \xi_j};
 \]
 it is invariantly defined at points $\xx$ belonging to $\Sigma$.
 Furthermore we shall denote by
\[
 \tp F_{p_m}\xx = \sum_{{\scriptstyle i\mu \in \hbox{{\scriptsize \rm
sp}}\, F_{p_m}\xx} \atop {\scriptstyle \mu > 0}}\mu .
\]
Since we are interested in micolocal results we shall always work in a
microlocal neighborhood, $U$, of $\x0x0 \in \Sigma$.

We make the following assumption on the lower order terms:
\begin{description}
	\item[(H${}_2$)]  If $\x0x0 \in \Sigma$ then $p_1^s\x0x0 + \tp
	F_{p_2}\x0x0 \notin \overline{\R^-}$.
\end{description}

Because of (H${}_1$) and (H${}_2$) we can find a canonical \gs\
transformation, $\Phi$, defined in $U$, such that
\begin{equation}
	\Phi \x0x0 = (0, e_n);
	\label{0.1}
\end{equation}
\begin{equation}
	{\rm if \quad} \yy = \Phi\xx, \quad y = (y', y'', y''') \in
	\R^{k+\ell+n-k-\ell},
	\label{0.2}
\end{equation}
$2 k$ being the rank of $\sigma_{|\Sigma}$, $\ell = \dim \ \im F_{p_2}\xx -
2k$.

Moreover in the coordinates $\yy$, $P$ can be written as
\begin{equation}
	\sum_{i,j=1}^{k}a_{ij}\yy X_i X_j^* + \sum_{j=1}^{k}\sum_{s=1}^{\ell}
	\left(b_{js}\yy X_j Y_s + b_{js}^*\yy X_j^* Y_s\right)
	\label{0.3}
\end{equation}
\begin{displaymath}
	+ \sum_{r,s=1}^{\ell}c_{rs}\yy Y_r Y_s + \tilde{p}_1\yy
	+ \sum_{j=1}^{k}\alpha_j'\yy + \sum_{j=1}^{k} \alpha_j''\yy X_j^*
\end{displaymath}
\begin{displaymath}
	+ \sum_{s=1}^{\ell}\beta_s\yy Y_s + \tilde{p}_0\yy ,
\end{displaymath}
where, denoting by $A\yy = \left[ a_{ij}\yy \right]_{i,j= 1, \ldots ,k}$,
$C\yy = \left[ c_{rs}\yy \right]_{r,s= 1, \ldots ,\ell}$, and
$B\yy = \left[ b_{js}\yy \right]_{{j= 1, \ldots ,k}\atop {s= 1, \ldots
,\ell}}$, the matrix
\[
\left[
\begin{array}{cc}
	A & B  \\
	B^* & C
\end{array}
\right]
 \]
is a self--adjoint, positive definite matrix of pseudo -- differential
operators of order $0$; $\tilde{p}_1\yy$ is a first order pseudo
differential operator such that $\tilde{p}_1\yy_{|\Sigma} = \left(\ps +
\tp F_{p_2}\right)_{|\Sigma}$; $\alpha_j'$, $\alpha_j''$, $j = 1, \cdots
,k$, $\beta_s$, $s = 1, \cdots ,\ell$, and $\tilde{p}_0$ are pseudo
differential operators of order $0$; moreover
\[
\sqrt 2 X = D_{y'} - i y' |D_{y_n}|, \quad Y = D_{y''}.
\]

Using M\'etivier's technique of addition of variables and making another
analytic canonical transformation we can write the operator in
(\ref{0.3}) as
\begin{equation}
	\langle {\cal A}\xx \left[
	\begin{array}{c}
		X  \\
		Y
	\end{array}
	\right], \left[
	\begin{array}{c}
		X^*  \\
		Y
	\end{array}
	\right]\rangle +
	\langle L\xx , \left[
	\begin{array}{c}
		X^*  \\
		Y
	\end{array}
	\right]\rangle + \tilde{p}_1\xx + \tilde{p}_0\xx ,
	\label{0.4}
\end{equation}
where
\begin{equation}
	\left\{
	\begin{array}{lc}
		X_j = \frac{1}{i}\frac{\textstyle\partial}{\textstyle \partial x_j}-
x_{k+j} |D_{x_n}|, & 1
		\leq j \leq k,  \\
		X_{k+j} = \frac{\textstyle\partial}{\textstyle\partial x_{k+j}} & 1 \leq j
\leq k,  \\
		Y_s = \frac{\textstyle\partial}{\textstyle\partial x_{2k+s}} & 1 \leq s
\leq \ell,
	\end{array}
	\right.
	\label{0.5}
\end{equation}
$X = \left(X_1,\ldots , X_k, X_{k+1}, \ldots , X_{2k}\right)$, $Y =
\left(Y_1, \ldots , Y_\ell\right)$, ${\cal A}$ is a self adjoint positive
definite matrix, of size $2 k + \ell$, of pseudo -- differential operators
of order $0$, and $L$ is a (complex) $2 k + \ell$ dimensional vector
of pseudo -- differential operators of order $0$.

If $\x0x0 \in \Sigma$, we denote by $P_{\x0x0}\xd$ the pseudo -- differential
operator obtained by freezing the coefficients of (\ref{0.4}) at $\x0x0$:
\begin{equation}
	P_{\x0x0}\xd = \langle {\cal A}\x0x0 \left[
	\begin{array}{c}
		X\xd  \\
		Y\xd
	\end{array}
	\right], \left[
	\begin{array}{c}
		X^*\xd  \\
		Y\xd
	\end{array}
	\right]\rangle
	\label{0.6}
\end{equation}
\begin{displaymath}
	+ \langle L\x0x0 , \left[
	\begin{array}{c}
		X\xd  \\
		Y\xd
	\end{array}
	\right]\rangle + \tilde{p}_1\x0x0 T + \tilde{p}_0\x0x0 ,
\end{displaymath}
where $T = \frac{\textstyle\partial}{\textstyle \partial x_n}$.

Due to (H${}_2$) it is easy to get an a priori estimate for $P_{\x0x0}\xd$:
\begin{equation}
	\sum_{|\alpha| \leq 2} \|X^\alpha u\|^2  + \sum_{\beta \leq
	2}\|Y^\beta\|^2 + \|u\|_1^2
	\label{0.7}
\end{equation}
\begin{displaymath}
	\leq C \left(\| P_{\x0x0}\xd u \|^2 + \|u\|^2\right), \qquad u \in
	C_0^\infty,
\end{displaymath}
where $C > 0$ is independent of $u$, $\alpha \in \Z_+^{2k}$, $\beta \in
\Z_+^\ell$ are multiindices and $\| \cdot \|_s$ means the microlocal
$H^s$ norm near $\x0x0$.

By Assumption (H${}_1$), (b)--(d), we know that $\Sigma$ is canonically
foliated with leaves of dimension $\ell$. If $\x0x0 \in \Sigma$ let us
denote by $\Gamma_{\x0x0}$ the leaf through $\x0x0$ and by $T_{\x0x0}
\Gamma_{\x0x0}$ its tangent space at $\x0x0$. Since we work in a
neighborhood of $\x0x0$ in $\ct$, in the sequel we will identify $\Gamma$
and its tangent space.

We are ready to state our results:
\begin{theorem}
\label{th0.1}
	Let $P$ be as above, verifying (H${}_1$) and (H${}_2$). Let $\x0x0 \in
	\Sigma$ and $W$ be a neighborhood of $\x0x0$. Suppose $1 \leq s < 2$ and
	that $\x0x0 \notin WF_s(Pu)$; then if $\Gamma_{\x0x0} \cap
	\left(W \setminus \{\x0x0\}\right) \cap WF_s(u) = \varnothing$ we have
	$\x0x0 \notin WF_s(u)$.
\end{theorem}
\begin{theorem}
	Under the same assumptions as in Theorem \ref{0.1}. Let $s \geq 2$ and
$\x0x0
	\notin WF_s(Pu)$. Then $\x0x0 \notin WF_s(u)$, i.e. $P$ is
	\gs--micro\-hy\-po\-el\-lip\-tic.
	\label{th0.2}
\end{theorem}
\begin{remark}
	When $s = +\infty$ Theorem \ref{0.2} is the well known result of
	H\"or\-man\-der \cite{hormander1}, Boutet -- Grigis -- Helffer \cite{bgh}.
	When $s = 1$ Theorem  \ref{0.1} is due to M\'etivier \cite{metivier 1981}.
\end{remark}

\section{Technical Considerations}
\label{s1}
\setcounter{equation}{0}
\setcounter{theorem}{0}
\setcounter{proposition}{0}
\setcounter{lemma}{0}
\setcounter{corollary}{0}
\setcounter{definition}{0}
\setcounter{remark}{0}
\renewcommand{\thetheorem}{\thesubsection.\arabic{theorem}}
\renewcommand{\theproposition}{\thesubsection.\arabic{proposition}}
\renewcommand{\thelemma}{\thesubsection.\arabic{lemma}}
\renewcommand{\thedefinition}{\thesubsection.\arabic{definition}}
\renewcommand{\thecorollary}{\thesubsection.\arabic{corollary}}
\renewcommand{\theremark}{\thesubsection.\arabic{remark}}
\renewcommand{\theequation}{\thesubsection.\arabic{equation}}
\renewcommand{\thesubsection}{\thesection.\arabic{subsection}}

\begin{definition}
	We say that $\x0x0 \notin WF_s(u)$, $u \in \D'$, $s \geq 1$, if there
	exists an open conic neighborhood of $\x0x0$, $V_0 \times \Gamma_0$, and
	a constant $C > 0$, such that for every $N \in \N$ there exists $v_N \in
	\E'$, $v_N = u$ in $V_0$, with
	\begin{equation}
		|\hat{v}_N(\xi)| \leq C^N \left(1 + \frac{|\xi|^{1/s}}{N}\right)^{-N},
		\quad \xi \in \Gamma_0.
		\label{1.0.1}
	\end{equation}
\end{definition}
In order to prove Theorems \ref{th0.1} and \ref{th0.2} we shall use the
estimate (\ref{0.6}) and some microlocalizations of high derivatives with
respect to characteristic directions.

\subsection{The Localizing Functions}
\label{s1.1}
\setcounter{equation}{0}
\setcounter{theorem}{0}
\setcounter{proposition}{0}
\setcounter{lemma}{0}
\setcounter{corollary}{0}
\setcounter{definition}{0}
\setcounter{remark}{0}

The definition of $WF_s$ given in (\ref{1.0.1}) refers to a fixed conic
neighborhood $V_0 \times \Gamma_0$ of $\x0x0$, valid for all $N$, and we
shall prove (\ref{1.0.1}) for $u$ and $V_0 \times \Gamma_0$, given
(\ref{1.0.1}) for $Pu$ and $\tilde{V}_0 \times \tilde{\Gamma}_0$, with
$V_0 \subset \tilde{V}$ and $\Gamma_0 \subset \bar{\tilde{\Gamma}}$. In
doing this we shall need to nest many such neighborhoods, using carefully
chosen cut--off functions, $\f(x)$, $\psi(x)$. In fact, once the
definition of $WF_s$ has been rewritten in terms of $L^2$ norms,
(\ref{0.7}) will allow us to estimate derivatives of order $\leq N$ in
$V_0 \times \Gamma_0$ in terms of those of order $\leq N/2$ without
changing cut--off functions; at this point we switch to a new pair for
the reduction to order $N/4$ etc., requiring $\log_2 N$ pairs in all.
Actually it is helpful to have two pairs for each step, with two
additional pairs initially.

Thus we shall use $2\log_2 N + 2$ pairs of functions $\{\f , \psi\}$:
$\{\f_j , \psi_j\}$, $\{\f_j' , \psi_j'\}$, $j = -1, 0, 1, \ldots ,
\log_2 N$, satisfying the following properties:
\begin{equation}
	\varphi_j\quad (\psi_j) \equiv 1 \qquad \mbox{{\rm near}}\ \mbox{{\rm
	supp}} \varphi_j \quad (\mbox{{\rm supp}}\ \psi_j),\quad j= 0,1,\ldots ,
	\log_2 N.
	\label{1.1.1}
\end{equation}
\begin{equation}
	\varphi_{-1} \equiv 1 \qquad \mbox{{\rm near}}\quad V_0, \quad \varphi_j \
(\varphi_j')
	\in C_0^\infty(\tilde{V}), \ \forall j.
	\label{1.1.2}
\end{equation}
\begin{eqnarray}
	\psi_{-1} \equiv 1 & \mbox{{\rm near}} & \Gamma_0 \cap \{|\xi| \geq 2
	N\} \quad \mbox{{\rm and}}
	\label{1.1.3} \\
	\psi_{-1} \equiv 0 & \mbox{{\rm for}} & |\xi| \leq N. \nonumber
\end{eqnarray}
\begin{eqnarray}
	\psi_j(\xi) \ \left(\psi_j'(\xi)\right) \equiv 0 & \mbox{{\rm for}} &
	|\xi| \leq N/2^j,
	\label{1.1.4} \\
	\psi_j \in C_0^\infty\left(\tilde{\Gamma}\right) & \mbox{{\rm for \
	every}} & j\geq 1.
\nonumber
\end{eqnarray}
\begin{equation}
	|D^\alpha\varphi_j^{(\prime)}(x)| \leq \left(K_j
N/2^j\right)^{|\alpha|},\quad
	|\alpha | \leq 3 N/2^j,\quad j = -1, 0, \ldots , \log_2 N.
	\label{1.1.5}
\end{equation}
\begin{equation}
	|D^\alpha\psi_j^{(\prime)}(\xi)| \leq \left(K_j N/\left(2^j |\xi| \right)
	\right)^{|\alpha|},\quad
	|\alpha | \leq 3 N/2^j,\quad j = -1, 0, \ldots , \log_2 N.
	\label{1.1.6}
\end{equation}
\begin{equation}
	\prod_{j=-1}^{\log_2 N} K_j^{N/2^j} \leq C^N.
	\label{1.1.7}
\end{equation}
\begin{equation}
	\mbox{{\hskip -\parindent \rm All the constants}}\quad K_j, \ C \  \
\mbox{{\rm are
	independent of $N$ and depend only on}}\ V_0.
	\label{1.1.8}
\end{equation}
We shall also need a form for the $\varphi_j^{(')}$, more adapted to the
geometry involved in our problem. If $x \in \R^n$, let $x = \left(x',
x'', x''', x_n\right)$ denote a partition of the variables according to
(\ref{0.5}), i.e. $x' \in \R^{2 k}$, $x'' \in \R^\ell$, $x''' \in \R^{n -
2 k - \ell -1}$. Then for every $j$
\begin{equation}
	\varphi_j^{(\prime)}(x) = \varphi_j^{(\prime)\#}(x'')
\tilde{\varphi}_j^{(\prime)}(x',
	x''', x_n).
	\label{1.1.9}
\end{equation}
The construction of the $\varphi_j^{(\prime)}(x)$ and
$\psi_j^{(\prime)}(\xi)$ is  easy; moreover these pairs of functions have
been used by the second  author before (\cite{tartakof80}), although the
idea of using  cut--off functions which behave in an analytic fashion up to
a given  order is due to L.~Ehrenpreis and has been exploited by
L.~H\"ormander,  K.G.~Andersson and others (\cite{hormander71},
\cite{andersson}). We give here a sketch of the construction.
Let $\Psi_i(t) \equiv 0$ for $|t| \leq N_i$, $\Psi_i(t) \equiv 1$ for
$|t| \geq 2 N_i$ and satisfy
\begin{equation}
	|D^\alpha_t \Psi_i (t)| \leq C^{|\alpha|+1}, \quad \mbox{{\rm for}}\
	\alpha \leq 3 N_i,
	\label{1.1.10}
\end{equation}
where $N_i = N/2^i$, $\Psi_i \in C^\infty(\R^n)$ and  the constant $C$ is
independent of $N$ and $i$. Such a function is obtained by convolving the
characteristic function of $\left\{|t| \leq 3 N_i/2\right\}$ with $3 N_i$
identical non--negative functions of integral one and support in $\left\{
|t| \leq 1/6\right\}$.

Let $U_0\cc V_{-1} \cc V_{-1}' \cc V_0 \cc \ldots \cc V_{\log_2 N}' \cc
\tilde{V}$, with
\[
\max \left\{\mbox{{\rm dist}}\ \left(V_i,{V'}_i^{{\rm \scriptstyle
comp}}\right), \mbox{{\rm dist}}\ \left ( V_i', V_{i+1}^{{\rm
\scriptstyle comp}} \right)\right\} = d_i,
\]
such that
\begin{equation}
	d_i = d_0/ 2^i.
	\label{1.1.11}
\end{equation}
Without loss of generality we may take each of the above mentioned nested
open sets in product form: $V_j^{(\prime)\#} \times
\tilde{V}_j^{(\prime)}$, where $V_j^{(\prime)\#} \subset \R^\ell$,
$\tilde{V}_j^{(\prime)} \subset \R^{n-\ell}$. Then we construct
$\varphi_i^{(\prime)\#}(x'')$,
$\tilde{\varphi}_i^{(\prime)}(x',x''',x_n)$ (and hence
$\varphi_i^{(\prime)}(x)$) just as we did $\Psi_i(t)$, but scaling by a
factor $d_i N_i$ (any positive order derivative of $\Psi_i$ had
support on a set of size $N_i$; now this distance is $d_0/2^i$).

The $\psi_i^{(\prime)}(\xi)$ are similarly defined by nesting open cones
$\Gamma_0 \cc \tilde{\Gamma}_{-1} \cc \tilde{\Gamma}_{-1}' \cc \ldots \cc
\tilde{\Gamma}_{\log_2 N}' \cc \tilde{\Gamma}$, with separations on the
unit sphere of $e_i = e_0/2^i$. Thes repeating on the unit sphere the
construction of the $\varphi_i^{(\prime)}$ --- but disregarding the
product form --- we construct the $\tilde{\psi}_i^{(\prime)}$; we then
extend the $\tilde{\psi}_i^{(\prime)}$ to be homogeneous of degree zero
and then take the product with $\Psi_i(|\xi|)$.

Note that the functions thus obtained, but {\em not\/} the constants,
depend on $N$.

\subsection{Constants}
\label{s1.2}
\setcounter{equation}{0}
\setcounter{theorem}{0}
\setcounter{proposition}{0}
\setcounter{lemma}{0}
\setcounter{corollary}{0}
\setcounter{definition}{0}
\setcounter{remark}{0}

A note on the use of constants. Any use of the letter $C$ denotes a
constant different from line to line, depending only on the dimension of
the space and the operator $P$, but independent of $N$ and $u$.

The constants $K_i$ or $\tilde{K}_i$ will be reserved for constants
satisfying (\ref{1.1.7}) and, like $C$, may change from line to line.

\subsection{Underlining}
\label{s1.3}
\setcounter{equation}{0}
\setcounter{theorem}{0}
\setcounter{proposition}{0}
\setcounter{lemma}{0}
\setcounter{corollary}{0}
\setcounter{definition}{0}
\setcounter{remark}{0}

Often we will be interested in the number of terms of a given form which
appear in an expansion. This is denoted by underlining a coefficient, and
may denote an upper bound rather than the exact count :
\begin{eqnarray}
	\left(\frac{\partial}{\partial x}\right)^a \left(x^b f(x)\right) & = &
	\sum_{c\leq \max\{a, b\}} {a \choose c} {b \choose c} c! x^{b-c}
	f^{(a-c)}(x)
	\nonumber \\
	 & = & \sum_{c\leq \max\{a, b\}} \underline{2^a b^c}\ x^{b-c} f^{(a-c)} (x)
	\nonumber
\end{eqnarray}
Another example, proved in the Appendix and often used below is
\[
x^{\alpha_1}\left({\frac {\partial} {\partial x}}\right )^{\beta_1} \ldots
x^{\alpha_r}  \left ({\frac {\partial} {\partial x}}\right )^{\beta_r} =
\sum_{\delta
\leq
\alpha ,
\beta} C^{|\beta|} |\alpha|^{|\delta|} x^{\alpha - \delta} \left ({\frac
{\partial} {\partial x}}\right )^{\beta - \delta},
\]
where $\alpha_i$, $\beta_i$ are multiindices, $|\alpha| = \alpha = \sum
\alpha_i$, $|\beta| = \beta = \sum \beta_i$.

\subsection{Pseudo -- differential Operators}
\label{s1.4}
\setcounter{equation}{0}
\setcounter{theorem}{0}
\setcounter{proposition}{0}
\setcounter{lemma}{0}
\setcounter{corollary}{0}
\setcounter{definition}{0}
\setcounter{remark}{0}

We shall use the Gevrey $s$ pseudo -- differential operators ($G^s$ pdo's)
of  Boutet de Monvel and Kr\'ee (\cite{boutetkree}); i.e. if
$v
\in  C_0^\infty(\R^n)$,
\begin{displaymath}
	p\xd v(x) =  \int e^{i x\cdot \xi} p\xx \hat{v}(\xi) \dt \xi,
\end{displaymath}
where $\dt \xi = (2\pi)^{-n} d\xi$, $p\xx \sim \sum_{k} p_k\xx$ and $p\xx$
($p_k\xx$) is real analytic  in $\Omega \times (\R^n \setminus \{0\})$,
$\Omega \subset \R^n$ open,  and the $p_k\xx$ are positively homogeneous of
degree $r - k$ with  respect to  $\xi$ and satisfy the following:
\par \noindent
$\forall K \cc \Omega$, $K$ compact, there exist constants $C$, $A$ such
that for any integer $k$, any $\alpha$, $\beta \in \Z_+^n$ and any $x \in
K$ we have
\begin{equation}
	|D_x^\alpha D_\xi^\beta p_k\xx | \leq C A^{k+|\alpha +\beta |}
	|\xi|^{r-k-|\beta |} (k+|\alpha |)!^s |\beta |!
	\label{1.4.1}
\end{equation}
and in addition, for any integer $N$ we have
\begin{equation}
	|D_x^\alpha D_\xi^\beta \left(p\xx - \sum_{k=0}^{N-1} p_k\xx \right) |
	 \leq C A^{N+|\alpha +\beta |}
	|\xi|^{r-N-|\beta |} (N+|\alpha |)!^s |\beta |!	
	\label{1.4.2}
\end{equation}
(with $|\xi|^{r-N-|\beta |}$ replace by $( 1 + |\xi|)^{r-N-|\beta |}$ if
$r - N - |\beta| > 0$).
Here $D_x = \frac{1}{i} \partial / \partial x$ and
$D_\xi = \frac{1}{i} \partial / \partial \xi$.
	
\section{Gevrey Hypoellipticity}
\label{s2}
\setcounter{equation}{0}
\setcounter{theorem}{0}
\setcounter{proposition}{0}
\setcounter{lemma}{0}
\setcounter{corollary}{0}
\setcounter{definition}{0}
\setcounter{remark}{0}

\subsection{Preliminary remarks}
\label{s2.1}
\setcounter{equation}{0}
\setcounter{theorem}{0}
\setcounter{proposition}{0}
\setcounter{lemma}{0}
\setcounter{corollary}{0}
\setcounter{definition}{0}
\setcounter{remark}{0}

\begin{proposition}
	\label{p2.1.1}
	To prove (\ref{1.0.1}) for $u$ it suffices to show that for $p \leq
	\frac{\textstyle n}{\textstyle s} + \frac{\textstyle n + 1}{\textstyle 2}$
	\begin{displaymath}
		\Vert T^p \psi_{-1} (D) \varphi_{-1} (x) u \Vert \leq C C^N N^{p s}.
	\end{displaymath}
	In other words it suffices to show that for $p \leq N$, $|\lambda | \leq
	N$,
\begin{equation}
	\Vert T^p \psi_{-1}^{(\lambda)} (D) \varphi_{-1 (\lambda)} (x) u \Vert
	\leq C C^N N^{s (p + |\lambda|)}
	\label{2.1.1}
\end{equation}
and that
\begin{eqnarray}
R_{1,N}(u) &  \equiv & \Vert T^p \left(\psi_{-1} (D) \varphi_{-1} (x) -
	\sum_{|\lambda | \leq N} \frac{1}{\lambda !}
(D_x^\lambda \varphi_{-1}) (x)
 (\partial_\xi^\lambda
	\psi_{-1}) (D)   \right) u \Vert
\nonumber \\
 & \leq  & \left(C K_{-1}\right)^N N!^s \label{2.1.2}
\end{eqnarray}
where $K_{-1}$ satisfies (\ref{1.1.8}).
\end{proposition}

{\bf Proof.}
First we show that (\ref{1.0.1}) is implied by the first estimate
\begin{displaymath}
	\Vert T^p \psi_{-1} (D) \varphi_{-1} (x) u \Vert \leq C C^N N^{p s},
	 p \leq \frac{N}{s} + \frac{n+1}{2}.
\end{displaymath}
In fact, replacing $u_N$ by $\varphi_{-1} (x) u$ - we may assume that
$(x_0,\xi_0) = (0,e_n)$ and that $u_N = u$ in the open neighborhood of $0$,
$U_0$,  with compact support there - we have the relation that
\begin{displaymath}
	|\hat{u}_N(\xi) | \leq C^N (1 + \frac{|\xi|^{1/s}}{N})^{-N}
\end{displaymath}
is equivalent to
\begin{equation}
	|\xi_n^k \hat{u}_N (\xi) | \leq C^N N^{k s} , k \leq \frac{N}{s}
	\label{2.1.3}
\end{equation}
This is easy to show since, possibly adjusting the constant C,
$c |\xi| \leq |\xi_n| \leq c'|\xi|$
in a small cone near $(0,e_n)$. The last inequality allows us to show that
\begin{displaymath}
	| \xi_n^k \psi_{-1}(\xi) \widehat{\varphi_{-1}(x) u}|
	\leq C^N N^{k s} , k \leq \frac{N}{s} , \quad {\rm if}\quad |\xi| \geq 2 N.
\end{displaymath}
When $|\xi| \leq 2N$ we have:
\begin{displaymath}
	|\xi_n^k \int \hat{\varphi}_{-1} (\eta) \hat{u} (\xi - \eta) \dt \eta|
\leq
 C (1 + |\xi|)^{k+M} \cdot
	\sup_{\eta} \left[(1 + |\eta|)^{M+n+1} |\hat{\varphi}_{-1}(\eta)| \right]
 \end{displaymath}
\begin{displaymath}
 \leq   C^N N^{(k+M) s} C^{M+n+1} N^{M+n+1}
	\leq  C^N N^{k s} ,
 \end{displaymath}
where, by the Paley-Wiener theorem,
 $|\hat{u}(\eta) | \leq C (1 + |\eta |)^M , M > 0$ depending on $n$. The
 insertion of a converging factor $(1 + |\xi|)^{-n-1}$ together with the
 above remark that $|\xi| \sim |\xi_n|$ yields the first assertion of the
 Proposition.

 Let us now turn to the second part of the Proposition. It suffices to
 prove (\ref{2.1.2}), since then (\ref{2.1.1}) will follow by general
 arguments of the calculus of pdo' s. We have
 \begin{displaymath}
 	\psi_{-1} (D_x) \varphi_{-1} (x)
  \sim
 	Op \left(\sum_{|\alpha| \geq 0} \frac{1}{\alpha !}
 	\partial_\xi^\alpha \psi_{-1} D_x^\alpha \varphi_{-1} (x) \right)
 \end{displaymath}
 \begin{displaymath}
  \sim
 \sum_{|\alpha| \leq N} \frac{1}{\alpha !}
 	\varphi_{-1}^{(\alpha)} (x) \psi_{-1(\alpha)}(D_x) +
 	\psi_{-1} (D_x) \varphi_{-1} (x)  -  \{ \psi_{-1} \circ \varphi_{-1} \}_N
 	(x, D_x)  ,
 \end{displaymath}
 where we used the notation
 \begin{equation}
 	\{ \psi_{-1} \circ \varphi_{-1} \}_N (x,\xi) = \sum_{0 \leq |\alpha| \leq
N}
 	\frac{1}{\alpha !} 	 D_x^\alpha \varphi_{-1} (x) \partial_\xi^\alpha
\psi_{-1}
 	\label{2.1.4}
 \end{equation}
 Thus
\begin{eqnarray*}
\lefteqn{\Vert T^p \psi_{-1}(D_x) \varphi_{-1} (x) u \Vert _{L^2}}    \\
 &\leq & \sum_{0 \leq |\lambda | \leq N} \frac{\textstyle
1}{\textstyle \lambda !}
\Vert T^p \varphi_{-1 (\lambda )} (x)
 	\psi_{-1}^{(\lambda)} (D_x) u \Vert _{L^2}   +  R_{1,N} (u),
\end{eqnarray*}
 where $R_{1,N} (u)$ is given by (\ref{2.1.2}). It is then enough to show
 that  $N^{|\lambda |}  \leq C^N \lambda !$ for some positive constant
 $C$, and for $|\lambda | \leq N$.
\par\noindent
 This completes the proof of the Proposition.\hfill$\blacksquare$

\begin{proposition}
\label{p2.1.2}
Using the same notation of Proposition \ref{p2.1.1}, we have that the
estimate (\ref{2.1.2}) holds.
\end{proposition}
 {\bf Proof.}
Recalling Lemma \ref{lA.1} of the Appendix we have that
\begin{eqnarray*}
\lefteqn{\Vert R_{1,N} (u) \Vert _{L^2}}    \\
 & = & \Vert T^p \left(\psi_{-1}(D_x) \varphi_{-1}(x) - \sum_{0\leq
 |\beta |\leq N}\frac{1}{\beta !} \varphi_{-1 (\beta)}(x) \psi_{-1}
 ^{(\beta)} (D_x) \right)u \Vert _{L^2} \\
 & = & \sum_{c_1}\underline{C^{p+n}} \left\| \int e^{i x \cdot \xi}
 \left(\int r_{c_1}(x,\eta,\xi) \dt \eta \right) \hat{u}(\xi) \dt \xi
\right\|
\end{eqnarray*}
where $c_1$ takes the values specified below and
\begin{displaymath}
	r_{c_1}(x,\eta,\xi) = \sum_{|\varepsilon |= N} \frac{(1 + |\eta
	|)^{-n-1}}{\varepsilon !} e^{i x \cdot \xi}
\left (\varphi_{-1 (\varepsilon + c_1 + b_2)}\right)\hat{}
 \; (\eta)
\end{displaymath}
\begin{displaymath}
	\times \left( \int_{0}^{1} \psi_{-1}^{(\varepsilon)} (\xi + \rho \eta) (\xi
+
	\rho \eta)^{p e_n - c_1 + b_1} (1 - \rho)^{|c_1 + \varepsilon |}
	d \rho \right);
\end{displaymath}
here $|b_1 + b_2| = n + 1$, $c_1 \leq p e_n + b_1$, and we have to study
the symbol $\int r_{c_1}(x,\eta,\xi) \dt \eta$. Now $D_x^\mu r_{c_1}$ has
a similar expression where $\left (\varphi_{-1  (\varepsilon + c_1 +
b_2)}\right )\hat{} \ $ is replaced by $\eta^\mu
\left (\varphi_{-1 (\varepsilon + c_1 +  b_2)}\right )\hat{}\ $, which is
rapidly decreasing in $\eta$ by the Paley-Wiener  theorem. Furthermore the
part under the integral sign is bounded by
\begin{displaymath}
	\sup_{\sigma \in \R^n ,|\varepsilon |=N}
	| \sigma^{\varrho e_n -c_1 + b_1} \psi_{-1}^{(\varepsilon)} (\sigma)|
	\leq
	C^N \; N!
\end{displaymath}
by (\ref{1.1.3}) and (\ref{1.1.6}). Hence we consider
\begin{displaymath}
	\sum_{|\mu | \leq n+1} \int |D_x^\mu \  r_{c_1}(x,\eta,\xi)| dx \dt \eta
\end{displaymath}
\begin{displaymath}
     =   \sum_{|a| \leq n+1} \int \frac{1}{(1+|x|^2)^{(n+1)/2}}
      (1+|x|^2)^{(n+1)/2} |D_x^\mu \  r_{c_1}(x,\eta,\xi)| \dt \eta dx,
\end{displaymath}
it is easily seen that integrating by part with respect to $\eta$ we get
an absolutely convergent integral. We may then conclude that the symbol
$\int r_{C_1} \dt \eta$ is $L_2$ continuous; thus
\begin{equation}
	\Vert R_{1,N} (u) \Vert _{L^2} \leq C^N \
	 \Vert u \Vert _{L^2} \  N!
	\label{2.1.5}
\end{equation}
where $ C = C' K_1$, $K_1$ given by (\ref{1.1.5}) - (\ref{1.1.8}).
\hfill$\blacksquare$
\vskip 1.cm

To deduce the microlocal Gevrey regularity of $u$  we shall obtain an
upper bound for a slightly more complicated expression. First some
notation.
Recalling the definition (\ref{0.5}) let us denote by $X = (X',X'')$ the
symplectic vector fields : $X_j'= X_j = \frac{\partial}{\partial x_j} -
x_{j+k} \frac{\partial}{\partial x_n}$, $X_j''= X_{j+k}=
\frac{\partial}{\partial x_{j+k}}$, $j=1,..,k$, $Y_s =
\frac{\partial}{\partial x_{2k+s}}$, $1\leq s \leq \ell$,
$T=\frac{\partial}{\partial x_n}$. We shall also denote by $Z$
either an $X-$ or an $Y-$ vector field.

If $I \in \Z_+^{2k+\ell}$ is a multi-index, the expression $Z^I$ means
$X_1^{I_1}$ $X_2^{I_2} \ldots X_{2k}^{I_{2k}}$
$Y_{2k+1}^{I_{2k+1}} \ldots Y_{2k+\ell}^{I_{2k+\ell}}$, and similarly $X^I$,
$I \in \Z_+^{2k}$,  means $X_1^{I_1} \ldots X_{2k}^{I_{2k}}$. Analogously
we write
$X'^{\alpha'}$, $X''^{\alpha''} \ldots $ etc, for $\alpha'$,$\alpha'' \in
\Z_+^k$.

In the proof of the microlocal regularity the quantity in (\ref{2.1.1})
will be replaced by
\begin{equation}
	\sup_{|J|\leq 2} \Vert Z^J T^p \varphi_{-1  (\lambda)}(x)
	\psi_{-1}^{(\lambda)} (D_x) u \Vert _{L_2}
	\leq
	C^{N+1} N^{s(p+|\lambda|)}, |\lambda| \leq N, p \leq N.
	\label{2.1.6}
\end{equation}
To proceed further a more effective microlocalization will be needed.
\subsection{Effective localization of powers of $T$}
\label{s2.2}
\setcounter{equation}{0}
\setcounter{theorem}{0}
\setcounter{proposition}{0}
\setcounter{lemma}{0}
\setcounter{corollary}{0}
\setcounter{definition}{0}
\setcounter{remark}{0}

\begin{definition}
\label{d2.2.1}
For any $\varphi(x)$, $\psi (D_x)$, $p \geq 0$ set
	\begin{equation}
		(T^p)_{\varphi \psi} = \sum_{|\alpha + \beta | \leq p}
		\frac{(-1)^{|\alpha|}}{\alpha ! \beta !} X''^{\alpha} X'^{\beta}
		T^{p-|\alpha +\beta |}\ \circ \ ad_{X'}^{\alpha} \
		 \ ad_{X''}^{\beta} (\varphi \psi) ,
		\label{2.2.1}
	\end{equation}
where $X''^\alpha $ has been defined above,
\begin{displaymath}
	ad_{X_i}^{\alpha _i} (W) =
	 \underbrace{[X_i,[X_i,...,[X_i,W]...]]}_{\alpha _i \ {\rm \scriptstyle
times}}
	\end{displaymath}
and $ad_{X'}^\alpha (W) = ad_{X_{1}}^{\alpha _{1}}\
ad_{X_{k}}^{\alpha _{k}} (W)$.
\end{definition}
\begin{remark}
	Recalling the definition (\ref{1.1.9}) of $\varphi _j$, we have
	\begin{equation}
		(T^p)_{\varphi _j \psi _j} = \varphi_j^{\#}(x'')
		(T^p)_{\tilde{\varphi}_j \psi _j}
		\label{r2.2.2}
	\end{equation}
\end{remark}
\begin{definition}
\label{d2.2.3}
	We shall write $Z^p$ for $Z^I$ with $|I| = p$ and $X^p$ for $X^J$, $|J| = p$
	and so on. Analogously $ad_X^p$ means $ad_X^\alpha$, $|\alpha | = p$ etc.
	If  $F$ is an operator we shall write
	\begin{equation}
		F \equiv_p 0
		\label{2.2.3}
	\end{equation}
	if
	\begin{equation}
		F = \underline{C^p} \  \frac{1}{p!} \  X^p \  ad_X^p(\varphi \psi)
		\label{2.2.4}
	\end{equation}
where $C$ is a universal constant.
	The index $p$ will be dropped if there is no possibility of
misunderstanding.	
\end{definition}
\begin{proposition}
	\label{p2.2.4}
We have the following properties of the microlocalizations (\ref{2.2.1}):
\begin{equation}
	[X_j' , \tphi p{\varphi\psi}] \equiv{}_p \; 0, \qquad  1 \leq j \leq k;
	\label{2.2.5}
\end{equation}
\begin{equation}
	[X_j'' , \tphi p{\varphi\psi}] \equiv{}_p \ X_j'' \tphi
	{p-1}{ad_T(\varphi\psi)} , \qquad  1 \leq j \leq k;
	\label{2.2.6}
\end{equation}
\begin{equation}
	[Y_j , \tphi p{\varphi\psi}] =  \tphi
	{p}{ad_Y(\varphi) \psi} , \qquad  1 \leq j \leq \ell.
	\label{2.2.7}
\end{equation}
\end{proposition}
{\bf Proof.}
By a calculation. (\ref{2.2.7}) is almost obvious since the vector field
$Y_j$ commutes with all the $X'$, $X''$. Thus we are left with the
verification of (\ref{2.2.5}) and (\ref{2.2.6}). As for (\ref{2.2.5}) we
have
\begin{displaymath}
	[X_j' , \tphi p{\varphi\psi}] = [X_j' , \sum_{|\alpha + \beta | \leq
	p}A_{\alpha \beta}^{(p)} B_{\alpha \beta} (\varphi\psi) ] ,
\end{displaymath}
where
\begin{displaymath}
	A_{\alpha \beta}^{(p)} = \frac{(-1)^{|\alpha |}}{\alpha ! \beta !}
	X''^{\alpha} X'^{\beta} T^{p - |\alpha + \beta |},
\end{displaymath}
\begin{displaymath}
	B_{\alpha \beta} (\varphi \psi) = ad_{X'}^{\alpha} ad_{X''}^{\beta}
	(\varphi \psi).
\end{displaymath}
Now
\begin{eqnarray*}
	[X'_j , A_{\alpha\beta}^{(p)}] & = & -\frac{(-1)^{|\alpha - e_j
	|}}{(\alpha - e_j)! \beta !} X''^{\alpha - e_j} X'^{\beta} T^{p -
	|\alpha - e_j + \beta |}  \\
	 & = & - A_{\alpha - e_j \ \beta}^{(p)},
\end{eqnarray*}
\begin{displaymath}
	[X'_j , B_{\alpha \beta} (\varphi \psi) ] = ad_{X'}^{\alpha + e_j}
	ad_{X''}^{\beta} (\varphi \psi) = B_{\alpha + e_j \ \beta} (\varphi \psi),
\end{displaymath}
so that
\begin{eqnarray*}
	[X_j' , \tphi p{\varphi\psi}] & = & - \sum_{|\alpha + \beta | \leq p-1}
	A_{\alpha\beta}^{(p)}  B_{\alpha + e_j \ \beta} (\varphi \psi) \\
	 &   & + \sum_{|\alpha + \beta | \leq p} A_{\alpha\beta}^{(p)}  B_{\alpha
	 + e_j \ \beta} (\varphi \psi)   \\
	 & =  & \sum_{|\alpha + \beta | = p} A_{\alpha\beta}^{(p)}  B_{\alpha
	 + e_j \ \beta} (\varphi \psi) \equiv_p 0.
\end{eqnarray*}
Let us now turn to (\ref{2.2.6}). This time
\begin{eqnarray*}
	[X''_j , A_{\alpha\beta}^{(p)} ] & = & - \frac{(-1)^{|\alpha |}}{\alpha
	! (\beta - e_j)!} X''^{\alpha} X'^{\beta - e_j} T^{p - |\alpha + \beta -
	e_j |}  \\
	 & = & - A_{\alpha \ \beta - e_j}^{(p)}
\end{eqnarray*}
and since $ ad_{X''_j} ad_{X'}^{\alpha}(v) = ad_{X'}^{\alpha} ad_{X''_j}
(v) - \alpha_j ad_{X'}^{\alpha - e_j} ad_T (v)$,
\begin{displaymath}
	[X''_j , B_{\alpha \beta} (\varphi \psi) ] = ad_{X'}^{\alpha}
	ad_{X''}^{\beta + e_j} (\varphi \psi) - \alpha_j ad_{X'}^{\alpha - e_j}
	ad_{X''}^{\beta} ad_{T} (\varphi \psi) ,
\end{displaymath}
so that
\begin{eqnarray*}
	[X''_j , \tphi p{\varphi\psi} ] & = &  - \sum_{|\alpha + \beta | \leq
	p-1} A_{\alpha \beta}^{(p)} B_{\alpha \ \beta + e_j} (\varphi\psi) +
	\sum_{|\alpha + \beta | \leq p}  A_{\alpha \beta}^{(p)} B_{\alpha \
	\beta + e_j} (\varphi\psi)  \\
	 &  & + \sum_{|\alpha + \beta | \leq p} \frac{(-1)^{|\alpha -
	 e_j|}}{(\alpha - e_j)! \beta !} X''_j X''^{\alpha - e_j} X'^{\beta}
	 T^{p - 1 - |\alpha -e_j + \beta |}  \\
		&  & \circ \  ad_{X'}^{\alpha - e_j} ad_{X''}^{\beta} ad_T (\varphi\psi)
\\
	 & \equiv_p & X''_j \tphi{p-1}{ad_T (\varphi \psi)},
\end{eqnarray*}
which proves the assertion.
\hfill$\blacksquare$
\vskip 1.cm \noindent
Proposition \ref{p2.2.4} can be iterated giving
\begin{proposition}
	\label{p2.2.5}
	We have, for $r$, $p \in \Z_+$
	\begin{eqnarray}
		[X^r , \tphi p{\varphi\psi} ] & =  &  \sum_{\ell = 1}^{r} (-1)^{\ell -
		1} { r \choose \ell} X^r \tphi{p-\ell}{ad_T^{\ell} (\varphi \psi)}
		\label{2.2.8} \\
		 & + & \sum_{\ell=0}^{r-1} \sum_{k=0}^{r-\ell} \underline{C}^{p-\ell}
		 \frac{1}{(p-\ell)!} X^{p+r-\ell - k -1} ad_{X,T}^{p + k + 1} (\varphi
		 \psi),
		\nonumber
	\end{eqnarray}
	where $ad_{X,T}^k (v) = ad_X^{k_1} ad_T^{k_2} (v)$, $k_1 + k_2 = k$.
	When $X^r = X'^r$, the first term on the right hand side of
	(\ref{2.2.8}) is missing.
\end{proposition}
{\bf Proof.}
By induction with a direct calculation. (\ref{2.2.8}) is obviously true
when $r = 1$. Suppose (\ref{2.2.8}) holds for a certain value of $r$ and
that every $X = X''$. Then
\begin{eqnarray*}
	[X^{r+1}, \tphi p{\varphi\psi} ] & = & X [X^r , \tphi p{\varphi\psi} ] +
	[X , \tphi p{\varphi\psi} ] X^r  \\
	 & = & \sum_{\ell = 1}^{r} (-1)^{\ell - 1} {r \choose \ell} X^{r + 1}
	 \tphi{p - \ell}{ad_T^{\ell} (\varphi \psi)} + X^{r + 1} \tphi{p -
	 1}{ad_T (\varphi \psi)}  \\
	 &  & - \sum_{\ell = 1}^{r} (-1)^{\ell - 1} {r \choose \ell} X^{r + 1}
	 \tphi{p - \ell - 1}{ad_T^{\ell + 1}(\varphi \psi)} \\
	 &  & - \sum_{\ell = 0}^{r -1} \sum_{k=0}^{r - \ell} \underline{C}^{p - 1
- \ell} \frac{1}{(p -
	 1 - \ell)!} X^{p + r - \ell - k - 1} ad_{X,T}^{p + k + 1} (\varphi
\psi)  \\
	 &  &  + \sum_{\ell = 0}^{r -1} \sum_{k=0}^{r - \ell}  \underline{C}^{p -
\ell}
	 \frac{1}{(p - \ell)!} X^{p + r - \ell - k} ad_{X,T}^{p + k + 1}
(\varphi \psi)  \\
	 &  & + \underline{C}^p \frac{1}{p!} X^{p+r} ad_X^{p + 1} (\varphi\psi) -
	 \underline{C}^p \frac{1}{p!} X^p [X^r , ad_X^{p+1} (\varphi \psi)].
\end{eqnarray*}
Now taking into account that
\begin{displaymath}
	[X^r , ad_X^{p+1} (\varphi \psi)] = \sum_{k=1}^{r} (-1)^{k-1} {r \choose
	k} X^{r-k} ad_X^{p+1+k} (\varphi \psi),
\end{displaymath}
we get the conclusion.\hfill$\blacksquare$
\subsection{Using suitably microlocalized norms}
\label{s2.3}
\setcounter{equation}{0}
\setcounter{theorem}{0}
\setcounter{proposition}{0}
\setcounter{lemma}{0}
\setcounter{corollary}{0}
\setcounter{definition}{0}
\setcounter{remark}{0}

The next step will be to replace $T^p \varphi_{-1 (\lambda)}
\psi_{-1}^{(\lambda)}$ in (\ref{2.1.1}) by the localization $\tphi
p{\varphi_0
\psi_0}$.

Recall that $\varphi_0 \equiv 1$ on a neighborhood of $\mbox{{\rm supp}}\
\varphi_{-1}$ and that $\psi_0 \equiv 1 $ near $\mbox{{\rm supp}}\
\psi_{-1}$; then we may write:
\begin{eqnarray}
	Z^J T^p \varphi_{-1 (\lambda)} \psi_{-1}^{(\lambda)}(D_x) & = & \sum_{p'
	\leq p} {p \choose p'} Z^J \varphi_{-1 (\lambda + (p - p')e_n)}(x)
	\psi_{-1}^{(\lambda)}(D_x) T^{p'}
	\label{2.3.1} \\
	 & = & \sum_{p'\leq p}\sum_{J'\leq J}{p\choose p'}{J\choose J'}
	 \varphi_{-1 (\lambda + p - p' + J - J')}(x) Z^{J'} \psi_{-1}^{(\lambda)}
T^{p'},
	\nonumber
\end{eqnarray}
where $\varphi_{-1 (\lambda + p - p' + J - J')} (x) = Z^{J - J'} T^{p -
p'} \varphi_{-1 (\lambda)}(x)$. From now on we shall often denote like
a derivative of $\varphi_{-1}$ such a blend of derivatives of
$\varphi_{-1}$.

Using the above notation we have
\begin{proposition}
	\label{p2.3.1}
	In order to prove the estimate (\ref{2.1.6}) it suffices to show that
	\begin{eqnarray}
	\lefteqn{\sup_{{p'\leq p}\atop {J'\leq J}} \| \varphi_{-1 (\lambda + p
	- p' + J - J')}(x) \psi_{-1}^{(\lambda)}(D_x) Z^{J'} \tphi{p'}{\varphi_0
	\psi_0} u \|_{L^2}}
	\label{2.3.2} \\
		 & \leq &  C^{N+1} N^{s(p + |\lambda | +
|J|)},\makebox[2in]{}
		\nonumber
	\end{eqnarray}
	where $|J| \leq 2$, $|\lambda | \leq N$, $p \leq N$.
\end{proposition}
{\bf Proof.}
>From (\ref{2.3.1}) we obtain
\begin{eqnarray}
	\lefteqn{Z^J T^p \varphi_{-1 (\lambda)}(x) \psi_{-1}^{(\lambda)}(D_x)}
	\label{2.3.3} \\
	 & = & \sum_{p'\leq p} \sum_{J'\leq J} {p\choose p'} {J\choose J'}
	 \varphi_{-1 (\lambda + p - p' + J - J')} (x) \left(
\psi_{-1}^{(\lambda)}(D_x) Z^{J'} T^{p'} \right.
	\nonumber \\
&  &   \left . + [Z^{J'} , \psi_{-1}^{(\lambda)}(D_x)]
	 T^{p'} \right)
\nonumber \\
	 & = & \sum_{p'\leq p} \sum_{J'\leq J} {p\choose p'} {J\choose J'}
	 \varphi_{-1 (\lambda + p - p' + J - J')} (x)  \psi_{-1}^{(\lambda)}(D_x)
Z^{J'} \tphi{p'}{\varphi_0 \psi_0}
	\nonumber \\
	 &  &  + \sum_{{p'\leq p}\atop {J'\leq J}} {p\choose p'} {J\choose J'}
	 \varphi_{-1 (\lambda + p - p' + J - J')} (x) \left( \psi_{-1}^{
	 (\lambda)} (D_x) Z^{J'}\left(T^{p'} - \tphi{p'}{\varphi_0
	 \psi_0}\right)\right .
	\nonumber \\
	&	 & + \left .  [Z^{J'} , \psi_{-1}^{(\lambda)}(D_x)] T^{p'}\right)
\nonumber \\
	 & = & A_1 + A_2 + A_3
	\nonumber
\end{eqnarray}
Clearly (\ref{2.3.2}) means that $A_1 \leq C^{N+1} N^{s (p+|\lambda | +
|J|)}$ which, in turn, implies (\ref{2.1.6}) provided we show that
\begin{equation}
	A_i \leq C^{N+1} N^{s (p + |\lambda | + |J|)}, \quad i = 2, 3.
	\label{2.3.4}
\end{equation}
Let us start considering $A_3$. Modulo constants bounded by $C^N$, the
generic term appearing in the sum contained in $A_3$ has the form
\begin{equation}
	\varphi_{-1 (\lambda_1)} (x) [ Z^{J'} , \psi_{-1}^{ (\lambda)}] T^{p'},
	\label{2.3.5}
\end{equation}
where $\lambda_1$ is a suitable multi--index. The term in (\ref{2.3.5})
can be rewritten as (see Equation (\ref{A.2}) for the definition of $\{
\cdot \}_M$)
\begin{displaymath}
	-\left[ \varphi_{-1 (\lambda_1)} (x) \psi_{-1}^{ (\lambda)} (D_x) \circ
	Z^{J'} - \left\{ \varphi_{-1(\lambda_1)} (x) \psi_{-1}^{ (\lambda)} (D_x)
\circ Z^{J'} \right\}_{1} \right] T^{p'}.
\end{displaymath}
Applying this operator to a smooth function $u$, by lemma (\ref{lA.1}) we
obtain that the result can be expressed as
\begin{displaymath}
	\underline{C^{p'+n}} \int e^{ix\xi} \int r_{c_1}(x,\eta ,\xi)
	\dt \eta \hat{u}(\xi) \dt \xi ,
\end{displaymath}
where
\begin{displaymath}
 r_{c_1}(x,\eta ,\xi) = \sum_{|\varepsilon |+1} (1 + |\eta |)^{-n-1}
 \widehat{Z^{J'(\varepsilon + b_2 + c_1)}} (\eta , \xi) e^{ix\eta}
 \end{displaymath}
 \begin{displaymath}
 	\cdot \ \int_{0}^{1} \varphi _{-1 (\lambda _1)} (x)
 	\psi _{-1}^{ (\lambda + \varepsilon)} (\xi + \varrho \eta)
 	(\xi + \varrho \eta)^{b_1 - c_1} (1-\varrho)^{|c_1 + \varepsilon |}
 	d\varrho.	
 \end{displaymath}
 Some explanation is in order here: $|b_1 + b_2| = p' + n +1, c_1 \leq
 b_1$ and by $ \widehat{Z^{J'(\varepsilon + b_2 + c_1)}} (\eta , \xi) $
we
 denoted the Fourier transform with respect to $X$ of the
 $(\varepsilon + b_2 +c_1) $ - derivative of the symbol $Z^{J'}(x,\xi)$
 of the $J'-$th power of $X$ or $Y$. It is worth noting that this Fourier
 transform could give a derivative of the $\delta -$distribution in the
 $\eta -$variable, since the $X'$s have polynomial coefficients. On the
 other hand, going back to (\ref{2.1.6}), since $Z^J T^p$ are local
 operators and $\varphi _{-1}^{\lambda _1} (x) $ has compact support
 contained in $V_{-1}$, we may think of smearing out each $Z^J$ by
 premultiplying it by a smooth function with compact support being
 identically $1$ in a neihborhood of supp $\varphi_{-1}$, e.g. by
 $\varphi _{-1}'$; in such a way we obtain symbols $Z^{J'}(x,\xi)$ which
 are polynomials with respect to $\xi$ and compactly supported in $x$.
 Thus $\widehat{Z^{J'(\varepsilon + b_2 + c_1)}} (\eta , \xi)$ is a
 plynomial in $\xi$ of order $\leq |J'|$ with rapidly decreasing
 coefficients in the $\eta$ variables. Hence $r_{C_1}$ is given by a sum
 of terms of the form:
 \begin{equation}
 	e^{i x\cdot \eta} \  \chi(\eta) \ \xi^\alpha
 	\int_{0}^{1}\varphi _{-1 (\lambda _1)} (x)
 	\psi _{-1 }^{(\lambda + \varepsilon)} (\xi + \varrho \eta)
 	(\xi + \varrho \eta)^{b_1 - c_1} (1-\varrho)^{|c_1 + \varepsilon |}
 	d\varrho,
 	\label{2.3.6}
 \end{equation}
 where $|\alpha| \leq |J'|$.  Writing $\xi = \xi + \varrho \eta - \varrho
\eta$ we may write (\ref{2.3.6}) as a sum, modulo constants of the strength
  $C^N$, of terms of the form
\begin{eqnarray}
\lefteqn{e^{i x \cdot \eta} \  \chi(\eta)  \eta^\beta \varphi _{-1
(\lambda_1)} (x) }
\label{2.3.7} \\
 & \cdot &  	\int_{0}^{1} 	\psi _{-1}^{ (\lambda + \varepsilon)} (\xi +
\varrho \eta) (\xi + \varrho \eta)^{b_1 - c_1 + \alpha - \beta}
 	(1-\varrho)^{|c_1 + \varepsilon |} \varrho^{|\beta|} d \varrho
\nonumber 	
\end{eqnarray}
 where $|\beta| \leq |\alpha| \leq |J'| \leq |J|$. The term under the
  integral sign in (\ref {2.3.7}) can be estimated by
  $N^{s(p+|\lambda|+|J|)}$ modulo the $N-$th power of a constant and since
  $\varphi _{-1}$ has compact support the symbol
  $\int r_{c_1}(x,\eta,\xi) \dt \eta $ is $L^2$ continuous, thus yielding
  the desired conclusion. Next we show the estimate in (\ref {2.3.4})
  when $i = 2$. Again applying the generic term of $A_2$ in (\ref
  {2.3.3}) to a smooth function $u$ we have to estimate the $L^2$ norm of
  a sum of terms of the type
  \begin{equation}
  	\varphi_{-1 (\lambda + p - p' + J - J')} (x) \psi_{-1}^{(\lambda)}
  	(D_x) \  Z^{J'} (T^{p'} - (T^{p'})_{\varphi_0 \psi_0}) (u)
  	\label{2.3.8}
  \end{equation}
 Recalling the definition of $(T^{p'})_{\varphi_0 \psi_0}$, modulo $N-$th
 power of a suitable constant, the quantity in (\ref{2.3.8}) is a sum of
 terms of the form
 \begin{equation}
 	\frac{1}{r!} \  \varphi_{-1 (\lambda + p - p' + J - J')} (x)
 	\psi_{-1}^{(\lambda)} (D_x) T^{p'-r} Z^{r+J'}\ ad_X^r (1 - \varphi_0
 	\psi_0) u,
 	\label{2.3.9}
 \end{equation}
 where $r\leq p'\leq p\leq N$. Denote by $h_r(x,\xi)$ the symbol of
 the operator $Z^{r + J'} ad_X^r (1 - \varphi_0 \psi_0)$.
Then supp  $h_r(x,\xi) \subset \complement (V_0 \times \Gamma _0)$,
 while
\begin{displaymath}
\hbox{\rm supp} \  \varphi_{-1 (\lambda + p - p' + J - J')} (x) \psi
_{-1}^{(\lambda)} (\xi)
\subsetneqq Int[V_0 \times \Gamma _0] .	
\end{displaymath}
hence the quantity in (\ref{2.3.9}) can be rewritten as
\begin{displaymath}
 \frac{1}{r!} \  \varphi_{-1 (\lambda + p - p' + J - J')} (x)
 \left[ \psi_{-1}^{(\lambda)} (D_x) , h_r(x,D_x) \right] u ,	
\end{displaymath}
since
\begin{displaymath}
\varphi_{-1 (\lambda + p - p' + J - J')} (x) h(x,D_x)
\psi_{-1}^{(\lambda)} (D_x) = 0 .
\end{displaymath}
Now by the same argument on the support of $h$ and $\varphi_{-1}
\psi_{-1}$ we have that
\begin{eqnarray*}
 \lefteqn{\varphi_{-1 (\lambda + p - p' + J - J')} \
 \sigma([\psi_{-1}^{(\lambda)}(D_x),h_r(x,D_x)])} \\	
&	= &  \varphi_{-1 (\lambda + p - p' + J - J')} \
	\sigma\left( \psi_{-1}^{(\lambda)}(D_x) \ h_r(x,D_x)
\phantom{\frac{1}{\beta!}}\right .\\
 &  &  \left .  - \sum_{|\beta | < N}
\frac{1}{\beta!} \psi_{-1}^{(\lambda +
\beta)} (D_x) \ h_{(\beta)} (x,D_x)\right)
\end{eqnarray*}
since
\begin{displaymath}
	 \varphi_{-1 (\lambda + p - p' + J - J')} \sum_{1 < |\beta | < N}
 \frac{1}{\beta!} \psi_{-1}^{(\lambda + \beta)} (\xi) \  h_{(\beta)} (x,\xi)
= 0.
\end{displaymath}
The conclusion then follows using Lemma \ref{lA.1} in the Apendix and
arguing as above.
\hfill$\blacksquare$
\begin{corollary}
  In order to prove (\ref{2.3.2}) it suffices to show that
  \begin{equation}
  	\sup_{|J|\leq 2} \Vert Z^J (T^p)_{\varphi_0 \psi_0} u \Vert _{L^2}
  	\leq
  	C^{N+1} \  N^{sp}.
  	\label{2.3.10}
  \end{equation}
	\label{c2.3.2}
\end{corollary}
{\bf Proof.}
>From (\ref{1.1.5}) and (\ref{1.1.6}) we have that
\begin{displaymath}
	\Vert \varphi_{-1 (\lambda_1)} (x) w \Vert _{L^2} \leq
	C^{|\lambda_1|} \ K_{-1}^{|\lambda_1|} \ N^{|\lambda_1|}
	 \Vert w \Vert _{L^2} \ \
\end{displaymath}
and
\begin{displaymath}
	\Vert \psi_{-1}^{(\lambda_2)} (x) w \Vert _{L^2} \leq
	C^{|\lambda_2|} \ K_{-1}^{|\lambda_2|} \
	 \Vert w \Vert _{L^2}
\end{displaymath}
The assertion then follows combining these two estimates.
\hfill$\blacksquare$
\vskip 1.cm

Actually we shall need to estimate more general derivatives of the
microlocalizing functions, due to the interactions of the vector fields
$Z$ with the cut-off functions. First some notation:
\begin{definition}
    Let $Z$ be a vector field on $T^* \R^n$. We denote by $Z_\sigma$ a
    map from $OpS^m(\R^m) \longrightarrow OpS^m(\R^n)$ defined by
    \begin{equation}
    	Z_\sigma h(x,D_x) = \op (Z(h(x,\xi))) (x,D_x)
    	\label{2.3.10bis}
    \end{equation}
    when $h(x,\xi) \in S^m(\R^n)$
	\label{d2.3.3}
\end{definition}

\begin{remark}
   \label{r2.3.4}
We have
\begin{displaymath}
	ad_{X_j'}(h) = \tilde{X}_{j\sigma}' \  h ,\ \  1 \leq j \leq k,
\end{displaymath}
	where
\[
\tilde{X}_j' = \frac{\partial}{\partial x_j} -
	x_{j+k} \  \frac{\partial}{\partial x_n} +
	\xi_n \  \frac{\partial}{\partial \xi_{j+k}},
\]
\begin{displaymath}
	ad_{X_j''}(h) = \tilde{X}_{j\sigma}'' \  h ,\ \  1 \leq j \leq k,
\end{displaymath}
	where
\[
\tilde{X}_j'' = X_j'' =  \frac{\partial}{\partial x_{j+k}},
\]
and
\begin{displaymath}
	ad_{X_s}(h) = \tilde{Y}_{s\sigma} \  h ,\ \  1 \leq s \leq k,
\end{displaymath}
where
\[
\tilde{Y}_s = \frac{\partial}{\partial x_{2k+s}}.
\]
\end{remark}
\begin{remark}
\label{r2.3.5}
By the preceding remark we have
\begin{equation}
	ad_Y^\alpha \ ad_{X'}^\beta \ ad_{X''}^\gamma \ (h) =
	\tilde{Y}_\sigma^\alpha \ \tilde{X}_\sigma' \ \tilde{X}_\sigma'' \ (h)
	\label{2.3.11}
\end{equation}
and in particular
\begin{equation}
	B_{\alpha \beta}(\varphi \psi) =
	ad_{X'}^\alpha \ ad_{X''}^\beta (\varphi \psi) =
	\tilde{X}_\sigma'^\alpha \ \tilde{X}_\sigma''^\beta (\varphi \psi).
	\label{2.3.12}
\end{equation}
Equation (\ref{2.3.11}) can be rewritten as
\begin{equation}
	ad_Z^\alpha (h)\ = \ \tilde{Z}_\sigma^\alpha \ h.
	\label{2.3.13}
\end{equation}
\end{remark}
To define our generalized derivatives we shall need the following vector
fields:
\begin{definition}
	\label{d2.3.6}
	Write
	\begin{displaymath}
		W_j' \ = \ \frac{\partial}{\partial x_j},\ \ 1 \leq j \leq k \ , 2k < j
		\leq n,
	\end{displaymath}
	\begin{displaymath}
		W_{j+k}' \ = \ \frac{\partial}{\partial x_{j+k}} -
		 x_j \frac{\partial}{\partial x_n}, \ \  1 \leq j \leq k,
	\end{displaymath}
	\begin{displaymath}
		W_j'' \ = \ \xi_n \frac{\partial}
{\partial \xi_{j+k}},\ \ 1 \leq j \leq k \ , 2k < j
		\leq n,
	\end{displaymath}
	\begin{displaymath}
		W \ = \ (W',W'')
	\end{displaymath}
	\begin{displaymath}
		\Xi_j \ = \ \frac{\partial}{\partial \xi_j},\ \ 1 \leq j \leq n-1
	\end{displaymath}
	\begin{displaymath}
		\Xi_n \ = \ \frac{\partial}{\partial \xi_n} - i
		\sum_{j=1}^{n} x_j \frac{\partial}{\partial \xi_{k+j}}.
	\end{displaymath}
\end{definition}
Note that ${\rm span} \ [ W, \Xi ] \ = {\rm span} \
[ \tilde{Z}, \frac{\partial}{\partial  x} ]$  and ${\rm span} \  [ \Xi] \ = \
{\rm span} \ [ \frac{\partial} {\partial \xi} ]$. Moreover
\begin{equation}
	[W_j,\tilde{Z}_k] = 0,
	\label{2.3.14}
\end{equation}
so that
\begin{eqnarray}
	W_{j,\sigma}B_{\alpha \beta}(\varphi \psi) & = &
	B_{\alpha\beta}(W_{j,\sigma}(\varphi\psi))
	\label{2.3.15} \\
	\Xi_{j,\sigma}B_{\alpha \beta}(\varphi \psi) & = &
	B_{\alpha\beta}(\Xi_{j,\sigma}(\varphi\psi))
	\nonumber \\
	W'_{j}\circ B_{\alpha \beta}(\varphi \psi) & = &
	B_{\alpha\beta}(W'_{j}\circ (\varphi\psi))
	\nonumber
\end{eqnarray}
where the notation $W_j \circ h\xd $ means the usual composition of
(pseudo) differential operators.
\begin{definition}
	\label{d2.3.7}
	Using the fields $W_j$ and $\Xi_j$, denote by
	$\left(\varphi^{(\prime)}_j \psi^{(\prime)}_j\right)^{(s)}$ any sum of
	$C^s$ operators of the form
	\begin{eqnarray}
		\lefteqn{\left(\varphi_j \psi_j\right)^{(s)}\xd}
		\label{2.3.16} \\
		 & = & \underline{\left(2^{-j}N\right)^{r_1 - r_2}} \Xi_\sigma^{r_1}
		 W_\sigma^{r_3} \left(W'^{r_2} \circ \varphi_j \psi_j\right) \xd ,
		\nonumber
	\end{eqnarray}
	where $s = | r_1 + r_3 |$, $|r_2| \leq |r_1|$, $C$ is a universal constant
and the operations $\Xi_\sigma$, $W_\sigma$ and $W' \circ \ $ may occur in
any order; in this sense Equation (\ref{2.3.16}) is a formal equation.
\end{definition}
\begin{proposition}
	\label{p2.3.8}
	For any multi--indices $a$, $b$ and for any $s$, such that $|a+b| + s
	\leq 2^{-j+1} N$, we have
	\begin{equation}
		|\partial_x^a \partial_\xi^b (\varphi_j \psi_j)^{(s)}\xx | \leq C (C
		K_j)^{|a+b|+s}(2^{-j} N)^{|a| +s}
		\left(\frac{2^{-j}N}{|\xi|}\right)^{|b|},
		\label{2.3.17}
	\end{equation}
	$j= 0, 1, \ldots$
\end{proposition}
{\bf Proof.}
For sake of brevity put $N_j = 2^{-j} N$. The quantity on the left hand
side of (\ref{2.3.17}) can be written as a sum of terms of the form
\begin{displaymath}
	\partial_x^a \partial_\xi^b N_j^{r_1 - r_2} \Xi_\sigma^{r_1}
	W_\sigma^{r_3} (W'^{r_2} \circ) (\varphi_j \psi_j) \xd .
\end{displaymath}
Let us take a look at the symbol of this operator: since $W'$ has a
symbol containing either $\xi_j$ or $x_j \xi_n$ we see that $W'^{r_2}
\circ h$ has a symbol of the form $\sum_{r_2' + r_2'' = r_2} \sum_{\ell
\leq r_2''} a_{r_2', r_2'', \ell} x^{r_2'' - \ell} \xi^{r_2''}
\partial_x^{r_2'} h$, so that if we apply $\Xi^{r_1}$ to this symbol and
remark that $\Xi_\sigma$ means either a $\xi$ -- derivative or a $\xi$ --
derivative multiplied by an $x_j$, we obtain that
\begin{eqnarray*}
	\lefteqn{\partial_\xi^b \Xi_\sigma^{r_1} (W'^{r_2} \circ ) h} \\
	 & = & \sum_{{r_2' + r_2'' = r_2}\atop {\ell \leq r_2''}} \sum_{\rho_1 ,
	 s_2 } a_{r_2', r_2'', \ell, \rho_1, s_2} \xi^{r_2'' - s_2}
	 \partial_\xi^{r_1 - s_2} \partial_x^{r_2'} h,
\end{eqnarray*}
provided $s_2 \leq r_1 + b$, $s_2 \leq r_2''$. Analogously, applying
$\partial_x^a W_\sigma^{r_3}$ to the above symbol, we obtain
\begin{eqnarray*}
	\lefteqn{\partial_x^a W_\sigma^{r_3} \partial_\xi^b \Xi_\sigma^{r_1}
	(W'^{r_2} \circ ) h} \\
	 & = & C \sum_{{\scriptstyle r_2' + r_2'' = r_2}\atop {\scriptstyle \ell
\leq r_2''}}
\sum_{\rho_1 ,
	 s_2 } \sum_{\rho_3 , s_1}
	 x^{r_1 + r_3 +r_2'' - \ell - s_1 - \rho_1 - \rho_3} \xi^{r_2'' - s_2}
	 \partial_\xi^{r_1 + b - s_2} \partial_x^{r_2'+ a + r_3 - s_1} h,
\end{eqnarray*}
where $s_1 \leq r_2' + r_3 + a$ and $s_2$ has the same bounds as above.
Roughly speaking we may say that the quantity in the left hand side of
(\ref{2.3.17}) can be written as a sum of terms of the form
\begin{displaymath}
	N_j^{r_1 - r_2} x^{(r_1 + r_3 + r_2'' - s_1)^\leq} \xi^{r_2'' - s_2}
	\partial_\xi^{r_1 + b - s_2} \partial_x^{r_2' + a + r_3 - s_1}
	(\varphi_j \psi_j) ,
\end{displaymath}
where $s_2 \leq r_1 + b$, $s_2 \leq r_2''$, $s_1 \leq r_2' + r_3 + a$,
$r_2' + r_2'' = r_2$, $|r_2| \leq |r_1|$ and $(r)^\leq  =$ integer that
can be bounded by $r$.

Now taking into account (\ref{1.1.3}) -- (\ref{1.1.6}), since $r_2 \leq
r_1$, $|\xi| \geq N_j$ on the support of $\psi_j$, and $|\xi| \leq 2 N_j$
on the support of any derivative of $\psi_j$, we obtain that the latter
term can be estimated by
\begin{displaymath}
	K_j^{|r_1 + b - s_2 + r_2' + a + r_3 - s_1|}
	\left(\frac{N_j}{|\xi|}\right)^{r_1+b-r_2''} N_j^{|r_1+r_3+a|}
\end{displaymath}
and this can be estimated by
\begin{displaymath}
	K_j^{|a+b|+s} N_j^{|a|+s} \left(\frac{N_j}{|\xi|}\right)^{|b|},
\end{displaymath}
thus proving the assertion.\hfill$\blacksquare$
\begin{corollary}
	\label{c2.3.9}
	Let $s\leq 2 N$, $w \in L^2(\R^n)$; then
	\begin{equation}
		\|(\varphi_j \psi_j)^{(s)}\xd w \|_{L^2} \leq C (C K_j)^{s+n+1}
		N_j^{s+n+1} \|w\|_{L^2}.
		\label{2.3.18}
	\end{equation}
\end{corollary}
{\bf Proof.}
By Theorem 18.1.11${}^\prime$ in \cite{hormanderbook}, if $H\xd $ is a
pseudo  differential operator of order $0$ whose symbol has $x$ -- support
contained in a fixed compact set $\bar{\Omega}$ of $\R^n$, we have
$$
\| H\xd \|_{L^2 \to L^2} \leq C \ {\rm vol}(\Omega)
\sup_{{x,\xi}\atop{|\rho|
\leq n+1}} |\partial_x^\rho h\xx |.
$$
This fact, together with the
preceding proposition give the assertion.
\hfill$\blacksquare$

\section{The a priori estimate}
\label{s3}
\setcounter{equation}{0}
\setcounter{theorem}{0}
\setcounter{proposition}{0}
\setcounter{lemma}{0}
\setcounter{corollary}{0}
\setcounter{definition}{0}
\setcounter{remark}{0}

\subsection{Preparations}
\label{s3.1}
\setcounter{equation}{0}
\setcounter{theorem}{0}
\setcounter{proposition}{0}
\setcounter{lemma}{0}
\setcounter{corollary}{0}
\setcounter{definition}{0}
\setcounter{remark}{0}

Let us write, using (\ref{0.4}),
\begin{equation}
	P(x,\xi) = \sum_{|\alpha| \leq 2} a_\alpha (x,\xi) Z^\alpha +
	           b(x,\xi) T
	\label{3.1.1}
\end{equation}
and
\begin{equation}
	P_{(x_0,\xi_0)}(x,\xi) = \sum_{|\alpha| \leq 2}
a_\alpha (x_0,\xi_0) Z^\alpha +
	                        b(x_0,\xi_0) T.
	\label{3.1.2}
\end{equation}
Then by (\ref{0.7}) we have the a priori estimate with frozen coefficients:
\begin{equation}
	\sum_{|J| \leq 2}\!{}^\prime \Vert Z^J v \Vert ^2 + \Vert v \Vert _1 ^2 \leq
	C \left(\Vert P_{(x_0,\xi_0)} (x,D) v \Vert ^2 + \Vert v \Vert ^2 \right)
	\label{3.1.3}
\end{equation}
Where $\sum_{|J| \leq 2} \!{}^\prime\Vert Z^J v \Vert ^2 =
       \sum_{|J| \leq 2 } \left( \Vert X^J v \Vert ^2 + \Vert Y^J v \Vert
       ^2 \right)$.
We may also allow $x$ to vary in a suitable small open set, e.g. supposing
$ v \in \czi(\R^n)$ with small support:
 \begin{equation}
 	\sum_{|J| \leq 2} \!{}^\prime \Vert Z^J v \Vert ^2 + \Vert v \Vert _1 ^2
\leq
	C \left(\Vert P_{(x_0,\xi_0)} (x,D) v \Vert ^2 + \Vert v \Vert ^2 \right),
 	\label{3.1.4}
 \end{equation}
 where
 \begin{equation}
 	P_{(x,\xi_0)}(x,\xi) = \sum_{|\alpha| \leq 2} a_\alpha (x,\xi_0)Z^\alpha+
	                        b(x,\xi_0) T.
 	\label{3.1.5}
 \end{equation}
 We have
 \begin{equation}
 	\Vert P_{(x,\xi_0)} (x,D) v \Vert   \leq
 	\Vert P(x,D) v \Vert + \Vert (P - P_{x,\xi_O}(x,D) \psi_i'(x,D)) v \Vert
 \label{3.1.6}
 	\end{equation}
 \begin{displaymath}
 	\leq  \varepsilon \ \sum_{|J| \leq 2}\!{}^\prime\Vert Z^J v \Vert  +
 	 C \sum_{|J| \leq 2}\!{}^\prime\Vert Z^J (1 - \psi_i' (D)) v \Vert  +
 	 \Vert P(x,D) \Vert,
  \end{displaymath}
 where we used the fact that $ a_\alpha (x,\xi_0) -  a_\alpha (x,D)$ has
 small $L^2$ norm when applied to $\psi_i' (D) v$, provided
 cone supp $\psi_i' (\xi)$ is small enough.\\
 Our purpose will be to deal with a function $r$ of the form
 \begin{equation}
 	v = Z^{I'} T^q\  \tphi p{(\varphi_i \psi_i)^{(r)}}\  u.
 	\label{3.1.7}
 \end{equation}
\begin{lemma}
	\label{l3.1.1}
	Let $r \leq 2 N_i$, $N_i = N\ 2^{-i}$, $|I'| + q + p \leq N_i$ and
	$w \in  L^2(\Omega)$. Then
	\begin{equation}
		\sum_{|J|\leq 2} \Vert Z^J (1 - \psi_i'(D)) Z^{I'}T^q\
		(T^p)_{(\varphi_i \psi_i)^{(r)}} w \Vert _{L^2(\Omega)}
		\label{3.1.8}
	\end{equation}
	\begin{displaymath}
		\leq C\ K_i^{N_i}\ N_i^{|I'|+p+q+r} \Vert w \Vert _{L^2(\Omega)}.
	\end{displaymath}
\end{lemma}
{\bf Proof.}  By Proposition \ref{A.4} we may write :
\begin{displaymath}
	Z^J (1 - \psi_i'(D)) Z^{I'}T^q\ (T^p)_{(\varphi_i \psi_i)^{(r)}} w
\end{displaymath}
\begin{displaymath}
	= \underline{C^{|J|+|I'|+p}\ N_i^\ell}(1 - \psi_i'(D))^{(J')}\
	x^{|I'+J''|+p-\ell}
\end{displaymath}
\begin{displaymath}
	\cdot D^{|I'+J''|+q+p-\ell}\ \frac{1}{p'!}\ (\varphi_i \psi_i)^{(s+p')}\ w,
\end{displaymath}
where $J' + J'' = J$, $p' \leq p$, $\ell \leq |I'+J''| + |p|$ (note
that  the increase of $s$ has been obtained from the definition of
 $( \varphi_i \psi_i)^{(r)}$, the definition of $(T^p)_{( \varphi_i
\psi_i)^{(r)}}$ and Equation ( \ref{2.3.11} ). The last quantity equals
 \begin{displaymath}
	 \underline{C^{|J|+|I'|+p}\ N_i^\ell}(1 - \psi_i'(D))^{(J')}\
	x^{|I'+J''|+p-\ell}\ \cdot
 D^{|I'+J''|+q+p-\ell}\ \frac{1}{p'!}\ (\varphi_i \psi_i)^{(s+p')}\ w
 \end{displaymath}
 \begin{displaymath}
	= C^{|J + I'| + P}\ N_i^\ell \ x^{|I'+J''|+p-\ell -\ell'}
	D^{|I'+J''|+q+p-\ell }\ \frac{1}{p'!}\
\end{displaymath}
\begin{displaymath}
\cdot
 [(1 - \psi_i')(D)^{(J'+\ell')}\ ,\ (\varphi_i \psi_i)^{(s+p')}]\ w,
\end{displaymath}
where $\ell' \leq |I'+J''|+p-\ell$. Thus in order to prove the Lemma it
suffices to prove that
\begin{displaymath}
	\Vert D_x^{\tau_1} [\psi_{i}^{(\tau_2)}, (\varphi_i \psi_i)^{(\tau_3)}]
\Vert
	\leq K_i^{N_i}\ N_i^{\tau_1+\tau_3} \Vert u \Vert,
\end{displaymath}
when $|\tau_1 + \tau_2 + \tau_3| \leq 2\ N_i$. Since $\psi_i'$ is a
function of $\xi$ only and due to the definition of the conical support of
$\psi_i, \psi_i'$ we must actually estimate
\begin{equation}
	\left\|D_x^{\tau_1} \left[ {\psi'}_{i}^{(\tau_2)} (\varphi_i
	\psi_i)^{(\tau_3)} - \sum_{|\beta| < N_i} \frac{1}{\beta !} (\varphi_i
	\psi_i)^{(\tau_3 + \beta)} {\psi'}_{i}^{ (\tau_2 + \beta)}\right] u
\right\|.
	\label{3.1.9}
\end{equation}
This is the same quantity of Lemma \ref{lA.1}, where now cone supp $f\xx$
has compact cosphere sections ($f = \psi'_i$). Using the notation of
Lemma \ref{lA.1}, we remark that
\begin{displaymath}
	\int |r_{c_1} (x, \eta, \xi)| \dt \eta \leq \frac{C}{M!} \sup
	\left|f_{(\varepsilon)}^{(a_2)} (x, \sigma) \sigma^{a_1 + b_1 -
	c_1}\right| |\left( h^{(\varepsilon + b_2 + c_1)}
	(\eta, \xi)\right)^{\hat{}} |,
\end{displaymath}
the supremum being taken over all $\varepsilon$, $|\varepsilon| = M$,
$a_1 + a_2 = a$, $| b_1 + b_2 | = |b| + n + 1$, $c_1 \leq a_1 + b_1$.

If $f$ has compact $x$ -- support or if $f$ is a function of $\xi$ only,
we obtain that, using again the notation of
Lemma \ref{lA.1}
\begin{eqnarray}
	\lefteqn{\left|D_x^a \left( F\circ H - \left\{ F\circ H\right\}_M\right)
	D_x^b w \right|}
	\label{3.1.10} \\
	 & \leq &  \frac{C^{|a+b|+n}}{M!} \sup_{{{{\scriptstyle |\xi|=M}\atop
	 {\scriptstyle |b_1+b_2|=|b|+n+1}}\atop {\scriptstyle c_1 \leq a_1 + b_1}}
\atop {\scriptstyle a_1 + a_2 = a}}
	 \left|f_{(\varepsilon)}^{(a_2)} (x, \sigma) \sigma^{a_1 + b_1 -
	c_1}\right| |\left( h^{(\varepsilon + b_2 + c_1)}
	(\eta, \xi)\right)^{\hat{}} | \| w\|.
	\nonumber
\end{eqnarray}
Applying (\ref{3.1.10}) to (\ref{3.1.9}) and taking into account the
properties (\ref{1.1.1}) -- (\ref{1.1.8}), we get that (\ref{3.1.9}) can
be estimated by
\begin{math}
	K_i^{N_i} N_i^{(\tau_1 + \tau_3)} \| u \| ,
\end{math}
and this ends the proof of the Lemma.
\hfil$\blacksquare$
\vskip 1.cm

Applying Lemma \ref{3.1.1} and the a priori inequality (\ref{3.1.3})
allows us to deduce:
\begin{eqnarray}
	\lefteqn{\sum_{|J| \leq 2}{}^\prime \| Z^J Z^{I'} T^q \tphi p{(\varphi_i
	\psi_i)^{(r)}} u \|_{L^2(\Omega)}}
	\label{3.1.11} \\
	 & \leq & C \left( \| P_i Z^{I'} T^q \tphi p{(\varphi_i
	\psi_i)^{(r)}} u \|_{L^2(\Omega)} + \|Z^{I'} T^q \tphi p{(\varphi_i
	\psi_i)^{(r)}} u \|_{L^2(\Omega)} \right .
	\nonumber \\
& & \left. + K_i^{N_i} N_i^{|I'| + q + p + r}
	\|u\|_{L^2(\Omega)} \right),
\nonumber
\end{eqnarray}
where $P_i$ is defined by
\begin{eqnarray}
	\lefteqn{P_i\xd }
	\label{3.1.12} \\
	 & = & \sum_{|\alpha|\leq2} \Phi_i(x) a_\alpha\xd \psi'_i(D) Z^\alpha +
	 \Phi_i(x) b\xd \psi'_i (D)T
	\nonumber \\
	 & = & \sum_{|\alpha|\leq 2} \tilde{a}_{\alpha i}\xd Z^\alpha +
	 \tilde{b}_i \xd T,
	\nonumber
\end{eqnarray}
where $\Phi_i(x) \equiv 1$ in a neighborhood of $\bar{\Omega}$, $\Phi_i
\in \czi (\Omega')$, $\Omega \cc \Omega'$ and $|D^\alpha \Phi_i (x)| \leq
(C N_i )^{|\alpha|}$ if $|\alpha| \leq 3 N_i$ (see e.g. (\ref{1.1.5})).
We point out that the introduction of such a function $\Phi_i$ is always
possible since the $L^2$ norms in (\ref{3.1.11}) are actually
$L^2(\Omega)$ -- norms.

\subsection{The use of the a priori estimate}
\label{s3.2}
\setcounter{equation}{0}
\setcounter{theorem}{0}
\setcounter{proposition}{0}
\setcounter{lemma}{0}
\setcounter{corollary}{0}
\setcounter{definition}{0}
\setcounter{remark}{0}

Let now $u \in \czi(\Omega)$. Assume that $i=0$ and apply (\ref{3.1.11})
for $i = 0$ with $q = 0 = r$, $|I'| = 0$:
\begin{eqnarray}
\lefteqn{	\sum_{|J| \leq 2} \| Z^J \tphi p{\varphi_0\psi_0} u
\|_{L^2(\Omega)}}
	\label{3.2.1} \\
	 & \leq  & C \left ( \| P_0\xd \tphi p{\varphi_0\psi_0} u
	 \|_{L^2(\Omega)} + \| \tphi p{\varphi_0\psi_0} u \|_{L^2(\Omega)}
	 \right .
	\nonumber \\
	 &  & + \left . K_0^{N_0} N_0^p \| u \|_{L^2(\Omega)} \right ).
	\nonumber
\end{eqnarray}
Our purpose is to iterate (\ref{3.2.1}) in the following sense: we must
estimate the term $\|P_0 \tphi p{\varphi_0 \psi_0} u \|$; in order to do
this we write $P_0 \tphi p{\varphi_0 \psi_0} u = \tphi p{\varphi_0 \psi_0}
P_0 u$  $+ [ P_0 , \tphi p{\varphi_0 \psi_0} ] u $. The first term is
known and hence it will be a good term in our estimate; on the other hand
the commutator generates a certain number of terms according to
Proposition \ref{p2.2.4}. In particular new $X'$s and $Y'$s appear
causing the number of $T'$s and the index $r$ to increase as $p$
decreases. Now in the final step of the first iteration there is a term
with $p = 0$ ; since basically a commutator of two $X'$s generates a $T$
vector field, at this point $q$ may be as large as the original $p$
divided by $2$. We are then allowed to reboot another iteration procedure
by introducing the next pair of cut off functions $\varphi_1 \psi_1$ and
use (\ref{3.1.11}) over and over.\\
Hence our main task will be to commute $P_i$ with
 $Z^{I'}\ T^q (T^p)_{(\varphi _i \psi _i)^{(r)}}$. \\
 Hence
 \begin{displaymath}
 	[P_i,Z^{I'}T^q (T^p)_{(\varphi _i \psi _i)^{(r)}}]
 \end{displaymath}
 \begin{displaymath}
   =   [P_i,Z^{I'}T^q] (T^p) _{(\varphi _i \psi _i)^{(r)}} \ + \
  Z^{I'} T^q (T^p) _{(\varphi _i \psi _i)^{(r)}} .
 \end{displaymath}
  Recalling that, from (\ref{3.1.12}), $P_i(x,D)=
   \sum_{|\alpha|\leq 2} \tilde{a} _{\alpha i} (x,D) Z^\alpha +
   \tilde{b} _i(x,D) T $, we may write
 \begin{eqnarray}
	   	\lefteqn{ [P_i ,Z^{I'}T^q (T^p)_{(\varphi _i \psi _i)^{(r)}} }
   	\label{3.2.2} \\
   	& = & \sum_{|\alpha|\leq 2}
\left [ \
   	[ \tilde{a} _{\alpha i} , Z^{I'}T^q ]
   	Z^\alpha (T^p) _{(\varphi _i \psi _i)^{(r)}}
\right .
   	\nonumber \\
   	&   & +\ \tilde{a} _{\alpha i} [Z^\alpha , Z^{I'}] T^q
   	(T^p)_{(\varphi _i \psi _i)^{(r)}}
   	\nonumber \\
    &  &  +\ Z^{I'} T^q [\tilde{a} _{\alpha i},
    (T^p)_{(\varphi _i \psi_i)^{(r)}}] Z^\alpha
   	\nonumber \\
    &  & +\
\left .  Z^{I'} T^q \tilde{a} _{\alpha i}
    [Z^\alpha ,(T^p)_{(\varphi _i \psi _i)^{(r)}}] \
\right ]
   	\nonumber \\
   	&  & +\ [ \tilde{b} _i,Z^{I'}T^q] T (T^p)_{(\varphi _i \psi _i)^{(r)}} +
   	\tilde{b} _i [T,Z^{I'}T^q] (T^p)_{(\varphi _i \psi _i)^{(r)}}
   	\nonumber \\
    &  & +\ Z^{I'}T^q [\tilde{b} _i,(T^p)_{(\varphi _i \psi _i)^{(r)}}] T
   	\nonumber \\
   	&  & +\ Z^{I'} T^q \tilde{b} _i [ T,(T^p) _{(\varphi _i \psi _i)^{(r)}} ]
   	\nonumber \\
   	& = & \sum_{|\alpha|\leq 2} \ \sum_{j=1}^{4} A_{\alpha j}^{(i)}  +
   	\sum_{j=1}^{4} B_j^{(i)},
   	\nonumber
 \end{eqnarray}
   where $\tilde{a} _\alpha , \tilde{b} $ are defined in (\ref{3.1.12}).
We deal first with the $A-$terms; the $B'$s will then be easy. First of
all we point out that if $X$ is a fixed vector field we have the identity
(easily proved by induction)
\begin{equation}
	[\tilde{a} _{\alpha i}, X^k] = \sum_{1\leq k'\leq k} {k \choose k'}
	ad_X^{k'} \tilde{a} _{\alpha i} X^{k-k'}
	\label{3.2.3}
\end{equation}
Iterating (\ref{3.2.3}) and using the multi-index notation we have
\begin{eqnarray}
	A_{\alpha 1}^{(i)} & =  & - \sum_{{{{\scriptstyle |I''| \leq |I'|} \atop {
\scriptstyle q' \leq q}}
	\atop {\scriptstyle |I''| + q' \geq 1}}} \underline{{|I'| \choose |I''|} {q
\choose q'}}
\tilde{a}_{\alpha i}^{(|I''| + q')}
	\label{3.2.4} \\
	 &  & \cdot  Z^{I' - I'' +\alpha}
	 T^{q-q'} \tphi p{(\varphi_i \psi_i)^{(r)}} ,
	\nonumber
\end{eqnarray}
where we denote by
\begin{equation}
	\tilde{a}_{\alpha i}^{(s)}\xd = ad_{Z_{i_1}} \ldots ad_{Z_{i_s}}
	(\tilde{a}_{\alpha i}\xd)
	\label{3.2.5}
\end{equation}
$Z_{i_1}$, $\ldots$ $Z_{i_s}$ being vector fields belonging to a ''fixed"
finite set of analytic vector fields (e.g. all the vector fields used
until now, i.e. the $Z$'s, the $W$'s and the $\Xi$'s may build such a
family).\\
Next
\begin{equation}
	A_{\alpha 2}^{(i)} = \underline{|\alpha| |I'|} \tilde{a}_{\alpha i}
	Z^{|\alpha|+|I'|-2} T^{q+1} (T^p)_{(\varphi _i \psi _i)^{(r)}}
	\label{3.2.6}
\end{equation}
where this term is missing if $|I'|=0$ and by $	Z^{|\alpha|+|I'|-2}$ we
denote an expression of the form $Z^\beta$, with
$|\beta|=|\alpha|+|I'|-2$; note that this term is also missing if no
commutator between $X'$-type and $X''$-type field is involved.\\
Furthermore using (\ref{2.2.8}) and making again the same conventions as
above we obtain:
\begin{eqnarray}
	A_{\alpha 4}^{(i)} & = & \sum_{\ell =1}^{|\alpha'|} (-1)^{\ell -1}
	{|\alpha'| \choose \ell} Z^{I'} T^q  \tilde{a}_{\alpha i}
	(T^{p-\ell})_{  (\tilde{\varphi} _i \psi _i)^{(r+\ell)}
	\varphi^{\sharp (\alpha-\alpha')}  }
	\label{3.2.7} \\
	& = & \sum_{\ell =0}^{|\alpha'| -1} \sum_{k =1}^{|\alpha'|- \ell}
	\underline{C}^{p- \ell} \frac{1}{(p-\ell)!} Z^{I'} T^q  \tilde{a}_{\alpha i}
	X^{p+|\alpha'|-\ell-k-1} (\tilde{\varphi} _i \psi _i)^{(r+p+k+1)}
	 \varphi^{\sharp (\alpha-\alpha')}
	\nonumber
 \end{eqnarray}
We point out explicitly that in the above equation $Z^\alpha$ has been
decomposed as $ Y^{\alpha-\alpha'} X^{\alpha'}$ ($X$ equals either $X'$
or $X''$) and the $Y$ vector fields act on $(T^p)_{(\f_i \psi_i)}$
according to (\ref{2.2.7}) - see also Remark \ref{r2.2.2} - whereas we
applied (\ref{2.2.8}) only to the $X$ vector fields (see also
(\ref{1.1.9})).
\begin{eqnarray}
	\lefteqn{ \sum_{|\alpha|\leq 2}
	( A_{\alpha 1}^{(i)} +  A_{\alpha 2}^{(i)} +  A_{\alpha 4}^{(i)} ) }
	\label{3.2.8} \\
	 & = &  - \sum_{|\alpha|\leq 2} \ \sum_{\ell =0}^{|\alpha|} \
	\sum_{{{{ \scriptstyle |I''| \leq |I'|} \atop { \scriptstyle q' \leq q}}
	\atop { \scriptstyle |I''| + q' \geq 1}}}
	\underline{ (-1)^\ell{|I'| \choose |I''|} {q \choose q'}{|\alpha|
	\choose \ell}  }
	\nonumber \\
	 &   &  \cdot \tilde{a}_{\alpha i}^{(|I''|+q')} Z^{I'-I''+\alpha}
	 T^{q-q'} (T^{p-\ell})_{(\f_i \psi_i)^{(r+\ell)}}
	 \nonumber \\
	 &  & + \sum_{|\alpha|\leq 2} \underline{|\alpha| |I'|} \tilde{a}_{\alpha i}
	 Z^{|\alpha|+|I'|-2} T^{q+1} (T^p)_{(\f_i \psi_i)^{(r)}}
	\nonumber \\
	 &  & + \sum_{ { {\scriptstyle |\alpha| \leq 2} \atop
{\scriptstyle \alpha'\leq \alpha} } }
	 \sum_{\ell =1}^{|\alpha| -1}
	 \sum_{k =0}^{|\alpha|- \ell}
	 \sum_{{{{\scriptstyle |I''| \leq |I'|} \atop {\scriptstyle q' \leq q}}
\atop {\scriptstyle |I''| + q'
\geq 1}}}
	  	\underline{ C^{p-\ell}{|I'| \choose |I''|} {q \choose q'}  }
	  	\frac{1}{(p-\ell)!}
	\nonumber \\
	 &  & \cdot \tilde{a}_{\alpha i}^{(|I''|+q')} Z^{|I'-I''|}
	 X^{|\alpha'|+p-\ell-k-1} T^{q-q'}
	 (\tilde{\f}_i \psi_
 i)^{(r+p+k+1)}
\f^{\sharp (\alpha-\alpha')}.
	 \nonumber
\end{eqnarray}
\begin{proposition}
	\label{p3.2.1}
	Assume $\tilde{a}_{\alpha i}$ is a $G^s$ pdo of order zero and assume
	that $\tilde{a}_{\alpha i}^{(\ell)}$ entails only $x-$derivatives of the
	symbol $\tilde{a}_{\alpha i}^{(\ell)}(x,\xi)$. Then if
	 $u \in L^2(\Omega)\cap {\cal E}'(\Omega)$,
	 \begin{equation}
	 	\| \tilde{a}_{\alpha i}^{(\ell)} \xd u \|_{L^2(\Omega)} \leq
	 	C^{\ell+1} \ell !^s \| u \|_{L^2(\Omega)}.
	 	\label{3.2.9}
	 \end{equation}
	 Here $s \geq 1$.
\end{proposition}
{\bf Proof.}
Since
\begin{eqnarray*}
	\lefteqn{\| \tilde{a}_{\alpha i}^{(\ell)} \xd u \|_{L^2(\Omega)}} \\
	 & \leq & C \sup_{{\scriptstyle |\rho| \leq n+1}\atop {\scriptstyle \xi \in
\R^n}} \| D_x^\rho
	 \sigma (\tilde{a}_{\alpha i}^{(\ell)}) \xx \|_{L^1(\Omega)} \| u
	 \|_{L^2(\Omega)} ,
\end{eqnarray*}
the conclusion follows at once from (\ref{1.4.1}) and the definition of
$\tilde{a}_{\alpha i}$ in (\ref{3.1.12}).
\hfil$\blacksquare$
\vskip 1.cm

Assembling the estimates for the terms in (\ref{3.2.8}) and
(\ref{3.1.11}) we obtain:
\begin{eqnarray}
	\lefteqn{\sup_{|J|\leq 2} \| Z^J Z^{I'} T^q \tphi p{(\varphi_i
	\psi_i)^{(r)}} u \|_{L^2(\Omega)}}
	\label{3.2.10} \\
	 & \leq & C \left \{  \| Z^{I'} T^q \tphi p{(\varphi_i
	\psi_i)^{(r)}} \Phi_i P u \|_{L^2(\Omega)}
  \phantom{\sup_{|a| \leq 2} A_{\alpha 3}^{(i)}{}_{L} }  \right .
	\nonumber \\
  &  & +
	\| Z^{I'} T^q \tphi p{(\varphi_i \psi_i)^{(r)}} (\Phi_i P - P_i) u
	\|_{L^2(\Omega)}
  \nonumber \\
	 &  & + \|Z^{I'} T^q \tphi p{(\varphi_i \psi_i)^{(r)}} u\|_{L^2(\Omega)}
	 + K_i^{N_i} N_i^{|I'| +p+q+r} \| u \|_{L^2(\Omega)}
	\nonumber \\
	 &  & + \sup_{ {{{\scriptstyle |\alpha|\leq 2, \ \ell \leq 2}\atop {
\scriptstyle |I''|
\leq |I'|}}
	 \atop {\scriptstyle q' \leq q}} \atop { \scriptstyle |I''| + q' \geq 1} }
C^{|I''| + q'}
	 N_i^{(|I''| + q') s} \| Z^{I' - I'' +\alpha} T^{q-q'} \tphi
	 {p-\ell}{(\varphi_i \psi_i)^{(r+\ell)}} u \|_{L^2(\Omega)}
	\nonumber \\
	 &  & + \sup_{|\alpha| \leq 2} |I'| \| Z^{\alpha + I' - 2} T^{q+1} \tphi
	 p{(\varphi_i \psi_i)^{(r)}} u \|_{L^2(\Omega)}
	\nonumber \\
	 &  & + \sup_{ {{{{ \scriptstyle |\alpha| \leq \ell, \ \alpha'\leq \alpha}
\atop {\scriptstyle 1\leq
	 \ell \leq |\alpha| - 1}} \atop { \scriptstyle 0 \leq k \leq |\alpha| -
\ell}} \atop
	 { \scriptstyle q'\leq q, \ |I''|\leq |I'|}} \atop { \scriptstyle |I''| + q'
\geq 1} } C^{|I''| +
	 q'} N_i^{(|I''| + q') s} \frac{1}{(p-\ell)!}
	\nonumber \\
		&  & \cdot  \|Z^{|I' - I''|}
	 X^{|\alpha'| +p-\ell-k-1 } T^{q-q'} (\varphi_i
	 \psi_i)^{(r+p+k+1)}\varphi^{\# (\alpha - \alpha')} u \|_{L^2(\Omega)}
   \nonumber \\
	 &  & \left . + \sup_{|\alpha| \leq 2} \| A_{\alpha 3}^{(i)}
u \|_{L^2(\Omega)}
	 \right \}
	\nonumber
\end{eqnarray}
for $|I'| + p + q + 2 \leq N_i$.
\begin{lemma}
	\label{l3.2.2}
	Let $|I'| + p + q + 2 \leq N_i$, $r \leq 2 N_i$; then
	\begin{eqnarray*}
		\lefteqn{\| Z^{I'} T^q \tphi p{(\varphi_i \psi_i)^{(r)}} (\Phi_i P - P_i) u
	\|_{L^2(\Omega)}}\\
		 & \leq &  K_i^{r+N_i+2} N_i^{r+|I'|+p+q} \| u \|_{L^2(\Omega)} ,
	\end{eqnarray*}
	$p$, $q$, $I'$ and $r$ being defined as above.
\end{lemma}
{\bf Proof.}
Recalling formula (\ref{A.4}) and remarking that $\psi'_i (\xi) \equiv 1$
on the support of $\psi_i$, we may write that
\begin{eqnarray}
	\lefteqn{Z^{I'} T^q \tphi p{(\varphi_i \psi_i)^{(r)}} (\Phi_i P - P_i)}
	\label{3.2.11} \\
	 & = &  \sum_{ {\scriptstyle r'+\tilde{r} \leq |I'|+p+q} \atop {
\scriptstyle r''\leq 2,\ |\alpha|
	 \leq 2} } \underline{C^{|I'|+p+q} N_i^{r'+r''}} O\left(|x|^{|I'|
	 +p-r'}\right) \frac{1}{p!} (\varphi_i \psi_i)^{(r+p+\tilde{r})}
	\nonumber \\
	 &  & D_x^{|I'|+p+q-r'- \tilde{r}} \left[ \Phi_i a_\alpha \ , (1 -
	 \psi_i')(D) \right] D_x^{|\alpha|- r''} O(|x|^{|\alpha|-r''}) ,
	\nonumber
\end{eqnarray}
where $D_x^k$ means a derivative with respect to $x$ of order $k$. The
first order term in $T$ in the expression of $P$ receives a completely
analogous treatment (even simpler!). Denoting by $\psi''(\xi)$ a cut-off
symbol of order zero such that supp $\psi_i''$ $\subset \{ \psi_i' \equiv
1\}$, $\psi_i' \equiv 1 $ on the support of $\psi_i$, $|\psi_i''| \leq
1$, we easily see that, due to Corollary \ref{c2.3.9} it suffices to show
that
\begin{eqnarray}
	\lefteqn{ \| \psi_i''(D) D_x^{\tau_1} \left[ \Phi_i a , \ (1 - \psi_i'(D))
	\right] D_x^{\tau_2} w \|_{L^2(\Omega)} }
	\label{3.2.12} \\
	 & \leq & K_i^{N_i} N_i^{|\tau_1 + \tau_2 |} \|w\|_{L^2(\Omega)} ,
	\nonumber
\end{eqnarray}
whenever $|\tau_1 + \tau_2 | \leq 2 N_i$, $\tau_1$, $\tau_2$ suitable
multi-indices. Note that since ${\rm supp}\ \psi_i'' \subset \{ \psi_i'
\equiv 1\}$, we have
\begin{eqnarray}
	\lefteqn{ \psi_i'' D_x^{\tau_1} \left[ \Phi_i a , \ (1 - \psi_i'(D))
	\right] }
	\label{3.2.13} \\
	 & = & \psi_i'' D_x^{\tau_1} \Phi_i a (1 - \psi_i'(D))
	\nonumber\\
	 & = & D_x^{\tau_1} \left[\psi_i''\ , \Phi_i a \right] (1 - \psi_i'(D))
	\nonumber \\
	 & = & D_x^{\tau_1} \left[ \psi_i'' \Phi_i a - \sum_{|\beta|\leq N_i}
	 \frac{1}{\beta!} \op\left( \psi_{i (\beta)}'' (\Phi_i
	 a)^{(\beta)}\right)\right] (1 - \psi_i'(D)).
	\nonumber
\end{eqnarray}
As a consequence we must estimate
\begin{eqnarray}
	\lefteqn{\left \| D_x^{\tau_1} \left[ \psi_i'' \Phi_i a - \sum_{|\beta|\leq
N_i}
	 \frac{1}{\beta!} \op\left( \psi_{i (\beta)}'' (\Phi_i
	 a)^{(\beta)}\right)\right] (1 - \psi_i'(D)) D_x^{\tau_2} w
	 \right \|_{L^2(\Omega)} }
	\label{3.2.14} \\
	 & \leq &  K_i^{N_i} N_i^{|\tau_1 + \tau_2 |} \|w\|_{L^2(\Omega)} ,
\makebox[3.4in]{}
	\nonumber
\end{eqnarray}
and the estimate in (\ref{3.2.14}) follows from Lemma \ref{lA.1} arguing
along the same lines of Lemma \ref{l3.1.1}.
\hfil$\blacksquare$

\subsection{Estimate of the term containing $A_{\alpha 3}^{(i)}$}
\label{s3.3}
\setcounter{equation}{0}
\setcounter{theorem}{0}
\setcounter{proposition}{0}
\setcounter{lemma}{0}
\setcounter{corollary}{0}
\setcounter{definition}{0}
\setcounter{remark}{0}

We will estimate in this section terms of the form $[G(x,D),(T^p)_{(\f_i
\psi_i)(r)}]$, modeling those which build the $A_{\alpha 3}^{(i)}$. But
first we need a suitable formula for the commutator of two
pseudo -- differential operators.
\begin{lemma}
	\label{l3.3.1}
	Let $G(x,D)$,  $H(x,D)$ be two pseudo -- differential operators with symbol
	$g(x,\xi)$, $h(x,\xi)$ respectively and let
	\begin{equation}
		g_{(\beta)}^{(\alpha)} (x,\xi) =
		\partial_\xi^\alpha D_x^\beta g(x,\xi),
		\label{3.3.1}
	\end{equation}
	\begin{equation}
		\sigma (G_{(\beta)}^{(\alpha)}) (x,\xi) =
		g_{(\beta)}^{(\alpha)} (x,\xi),
		\label{3.3.2}
	\end{equation}
	for any multi-indices $\alpha,\beta$. Also define, for a positive
	integer $M$,
       \begin{equation}
       	\sigma(\{G \circ H\}_M) = \sum_{0\leq \delta \leq M}
       	\frac{1}{\delta!} g^{(\delta)} \xx h_{(\delta)} \xx,
       	\label{3.3.3}
       \end{equation}
       i.e.
       \begin{equation}
       	\{G \circ H \}_M = G \circ H -
       	\left(
       	G \circ H - \op (\sum_{0\leq \delta \leq M}
       	\frac{1}{\delta!} g^{(\delta)} \xx h_{(\delta)} \xx)
       	\right),
       	\label{3.3.4}
       \end{equation}
     where, as usual, $\op (q \xx)$ denotes the pseudo-differential
     operator with symbol $q \xx$. Then for any $M$ non -- negative integer,
     \begin{equation}
     	\{G \circ H\}_M - \{H \circ M\}_M =
     	\sum_{1\leq |\alpha + \beta | \leq M}
        \frac{(-1)^{|\beta|}}{\alpha! \beta!}
        \ \{H_{(\beta)}^{(\alpha)} \circ
        G_{(\beta)}^{(\alpha)}\}_{M-|\alpha+\beta|} .
     	\label{3.3.5}
     \end{equation}
 \end{lemma}
 \begin{corollary}
 	\label{c3.3.2}
 	For any $M$ we have
 	\begin{equation}
 		[G,H] = \sum_{1\leq |\alpha + \beta | \leq M}
        \frac{(-1)^{|\beta|}}{\alpha! \beta!}
        \ H_{(\beta)}^{(\alpha)} \circ G_{(\beta)}^{(\alpha)} +
        R_{[G,H]_{,}M}
 		\label{3.3.6}
 	\end{equation}
 	where
 	\begin{eqnarray}
 		\lefteqn{ R_{[G,H]_{,}M} = [G,H] -(\{G \circ H\}_M - \{H \circ M\}_M) }
 		\label{3.3.7} \\
 		& = & \sum_{1\leq |\alpha + \beta | \leq M}
 		\frac{(-1)^{|\beta|}}{\alpha! \beta!}
 		\left( \{H_{(\beta)}^{(\alpha)} \circ
        G_{(\beta)}^{(\alpha)}\}_{M-|\alpha+\beta|} -
         H_{(\beta)}^{(\alpha)} \circ G_{(\beta)}^{(\alpha)}
 		\right)
 		\nonumber
 	 \end{eqnarray}.
 \end{corollary}
 Corollary \ref{c3.3.2} is just a restatement of Lemma \ref{l3.3.1}.\\
{\bf Proof of Lemma \ref{l3.3.1}.}
By the general calculus we have
\begin{displaymath}
	\sigma ([G,H]) \sim \sum_{|\gamma| \geq 1} \frac{1}{\delta!}
	(g^{(\delta)} h_{(\delta)} - h^{(\delta)} g_{(\delta)})
\end{displaymath}
Using the identity
\begin{displaymath}
	- 1 = \sum_{{|\delta| \geq 1} \atop {\gamma \leq |\delta|}}
	(-1)^{|\gamma|} {\delta \choose \gamma}
\end{displaymath}
we obtain thet
\begin{displaymath}
	- \frac{1}{\delta!}  h^{(\delta)} g_{(\delta)} =
	\sum_{1\leq |\gamma | \leq |\delta|}
	\frac{(-1)^{|\gamma|}}{\gamma! (\delta - \gamma)!}
	h^{(\gamma)(\delta-\gamma)} g_{(\gamma)(\delta-\gamma)},
\end{displaymath}
so that
\begin{displaymath}
	- \sum_{1\leq |\delta | \leq M} \frac{1}{\delta!}  h^{(\delta)} g_{(\delta)}
	= \sum_{1\leq |\gamma | \leq M}
	\sigma( \{ H^{(\gamma)} \circ G_{(\gamma)} \}_{M-|\gamma|}).
\end{displaymath}
Since
\begin{displaymath}
	h_{(\delta)} g^{(\delta)} =
	\sum_{0\leq |\beta | \leq M-|\delta|} \frac{1}{\beta!}
	 h^{(\delta)} g_{(\delta)} -
	 \sum_{1\leq |\beta | \leq M-|\delta|} \frac{1}{\delta!}
	  h^{(\delta)} g_{(\delta)},
\end{displaymath}
using the same identity on the second term we obtain
\begin{displaymath}
	 \sum_{1 \leq |\delta | \leq M} \frac{1}{\delta!}
	 g^{(\delta)} h_{(\delta)} = \sum_{1 \leq |\delta | \leq M}
\frac{1}{\delta !} \
	 \sigma ( \{ H_{(\delta)} \circ G^{(\delta)} \} _{M-|\delta|})
 \end{displaymath}
\begin{displaymath}
	  + \sum_{1 \leq |\delta | \leq M}\  \sum _{1 \leq |\gamma | \leq
	 M-|\delta|}
	 \frac{(-1)^{|\gamma|}}{\delta! \gamma!} \
	 \sigma ( \{ H_{(\delta)}^{(\gamma)} \circ G_{(\gamma)}^
	 {(\delta)} \} _{M-|\delta|} ),
 \end{displaymath}
which proves the Lemma.
\hfil$\blacksquare$
\vskip 1.cm

Our purpose is now to give an estimate of the last term in
(\ref{3.2.10}). We use the convention that $\tilde{H}(x,D)$ denotes the
operator
\begin{displaymath}
	\tilde{H}(x,D) = \Phi_i(x) H(x,D) \Psi_i'(D),
\end{displaymath}
where $\f_i \equiv 1 $ near $\bar{\Omega}$ and $|D^\alpha \Phi_i| \leq
(C N_i)^{|\alpha|}$, for $|\alpha| \leq  N_i$. It will also be useful to
make a small change in the previous notation; namely we set
\begin{equation}
	(T^p)_{\f \psi} = \sum_{|\alpha+\beta|\leq p}
	\frac{(-1)^{|\alpha|}}{\alpha! \beta!} A_{\alpha,\beta}^{(p)} \cdot
	B_{\alpha,\beta}^{(p)},
	\label{3.3.8}
\end{equation}
where
\begin{equation}
	A_{\alpha,\beta}^{(p)} = X''^\alpha X'^\beta T^{p-|\alpha+\beta|}
	\label{3.3.9}
\end{equation}
\begin{equation}
	B_{\alpha,\beta}^{(p)} = ad_{X'}^\alpha ad_{X''}^\beta (\f \psi) =
	B_{\alpha,\beta}(\f \psi) =	B_{\alpha,\beta}.
	\label{3.3.10}
\end{equation}
Using the general formula
\begin{displaymath}
	[AB,\tilde{G}] = A [B,\tilde{G}] + [A,\tilde{G}] B,
\end{displaymath}
and Corollary \ref{c3.3.2}, we may write:
\begin{eqnarray}
	 \lefteqn{ [(T^p)_{(\f_i \psi_i)^{(r)}}, \tilde{G}] =
	 R_{[(T^p)_{(\f_i \psi_i)^{(r)}}, \tilde{G}],N_i}
	 }
	\label{3.3.11} \\
	 &  & + \sum_{   {{{\scriptstyle |\gamma+\delta|\leq N_i} \atop
	 {\scriptstyle \alpha'\leq \alpha,\beta' \leq \beta}} \atop
	 {\scriptstyle \tau+|\alpha+\beta| \leq p}} \atop
	 {\scriptstyle 1\leq |\gamma+\delta|+|\alpha'+\beta'|+\tau}   }
	 (-1)^{|\alpha|} {A \choose \beta'} {\beta \choose \beta'}
	 {p-|\alpha+\beta| \choose \tau}
	 \tilde{G}_{(\alpha',\beta',\tau,\delta)}^{(\gamma)}
	\nonumber \\
	 &  & A_{\alpha-\alpha',\beta-\beta'}^{(p-|\alpha'+\beta'|-\tau)}
	 B_{\alpha,\beta} ((\f_i \psi_i)^{(r)})_{(\gamma)}^{(\delta)}
	 \frac{1}{\alpha!\beta!\gamma!\delta!}
	\nonumber \\
	 & = & R_{[(T^p)_{(\f_i \psi_i)^{(r)}}, \tilde{G}],N_i} +
\Sigma_{r,i,\tilde{G}}^{1} +
	 \Sigma_{r,i,\tilde{G}}^{2},
	\nonumber
\end{eqnarray}
where
\begin{eqnarray}
\lefteqn{ R_{[(T^p)_{ (\f_i \psi_i)^{(r)}}, \tilde{G} ],N_i } =
	\sum_{|\alpha+\beta| \leq p} \frac{(-1)^{|\alpha|}}{\alpha! \beta!}
  }
\label{3.3.12} \\
 & & \cdot \left [ A_{\alpha,\beta}^{(p)}
( B_{\alpha,\beta}^{(p)} \tilde{G}
	 - \{ B_{\alpha,\beta}^{(p)} \circ \tilde{G} \} _{N_i})
		 - A_{\alpha,\beta}^{(p)} ( \tilde{G} B_{\alpha,\beta}^{(p)}
	 - \{ \tilde{G} \circ B_{ \alpha,\beta }^{(p)} \} _{N_i} )
\phantom{\frac{(- 1)^{\lambda}}{\lambda !}}
\right .
 \nonumber \\
	& & \left . +  A_{\alpha,\beta}^{(p)}  \sum_{1 \leq |\lambda+\mu | \leq N_i}
\
	 \frac{(-1)^{|\lambda|}}{\lambda!\mu!} \
	 (   \{ \tilde{G} _{(\mu)}^{(\lambda)}
	 \circ B_{\alpha,\beta (\lambda)}^{(p)(\mu)} \} _{N_i-|\lambda+\mu |}
	 - \tilde{G} _{(\mu)}^{(\lambda)} \circ
	  B_{\alpha,\beta (\lambda)}^{(p)(\mu)}  ) \right ] ,
	  \nonumber
\end{eqnarray}
\begin{equation}
	\tilde{G} _{(\alpha',\beta',\tau,\delta)}^{(\gamma)} (x,D) =
	ad_{X''}^{\alpha'} ad_{X'}^{\beta'} ad_{T}^{\tau }
	(  \tilde{G} _{(\delta)}^{(\gamma)} (x,D))
	\label{3.3.13}
\end{equation}
and where $\Sigma_{r,i,\tilde{G}}^{1}$ and $\Sigma_{r,i,\tilde{G}}^{2}$
denote the sums over $|\gamma+\delta| \leq p - |\alpha' + \beta'| -\tau$
and $|\gamma+\delta| > p - |\alpha'+\beta'| - \tau$ respectively. Here
$\tilde{G} _{(\alpha',\beta',\tau,\delta)}^{(\gamma)}$ has order
$-|\gamma|$ and $B_{\alpha,\beta} ((\f_i \psi_i)^{(r)})_{(\gamma)}^{(\delta)}$
has order $-|\delta|$.\\
In $\Sigma_{r,i,\tilde{G}}^{1}$ we want to decrease the number of
derivatives in
$A_{\alpha-\alpha',\beta-\beta'}^{(p-|\alpha'+\beta'|-\tau)}$ in order to
bring $\tilde{G} _{(\alpha',\beta',\tau,\delta)}^{(\gamma)}$ and
$B_{\alpha,\beta} ((\f_i \psi_i)^{(r)})_{(\gamma)}^{(\delta)}$ up to
order zero: when possible we do this exploiting exclusively powers of $T$
(grouping these  terms in $E_{r,i, \tilde{G}}^1$) and with a mixture of
$T$ and $X$ derivatives otherwise (grouping these terms in
 $E_{r,i, \tilde{G}}^2$):
 \begin{equation}
 	\Sigma _{r,i, \tilde{G}}^1 = E_{r,i, \tilde{G}}^1 +
 	     E_{r,i, \tilde{G}}^2,
 	\label{3.3.14}
 \end{equation}
 where
 \begin{eqnarray}
 \lefteqn{	E_{r,i, \tilde{G}}^1 =
 \sum_{
 	{{ \scriptstyle |\gamma+\delta| \leq p-|\alpha'+\beta'|-\tau} \atop
 	{ \scriptstyle |\gamma+\delta|+|\alpha'+\beta'|+\tau \geq 1}} \atop
 	{ \scriptstyle \alpha'\leq \alpha, \ \beta' \leq \beta}
  }
 	(-1)^{|\alpha|} {\alpha \choose \alpha'} {\beta \choose \beta'}
 	{p-|\alpha+\beta| \choose \tau}
  }
 	\label{3.3.15} \\
 	 &  & \cdot \ \tilde{G} _{(\alpha',\beta',\tau,\delta)}^{(\gamma)}
 	 T^{|\gamma|}
 	 A_{\alpha-\alpha',\beta-\beta'}^{(p-|\alpha'+\beta'|-\tau-|\gamma+\delta|)}
 	 T^{|\delta|}
   \cdot \ B_{\alpha,\beta} ((\f_i \psi_i)^{(r)})_{(\gamma)}^{(\delta)}
   \frac{1}{\alpha! \beta! \delta! \gamma!}
 	\nonumber
 \end{eqnarray}
and
 \begin{eqnarray}
 \lefteqn{	E_{r,i, \tilde{G}}^2 =
 \sum_{
 	{{\scriptstyle |\gamma+\delta| \leq p-|\alpha'+\beta'|-\tau} \atop
 	{\scriptstyle |\gamma+\delta+\alpha'+\beta'|+\tau \geq 1}} \atop
 	{ \scriptstyle \alpha'\leq \alpha \ ,\beta' \leq \beta}       }
 	\frac{(-1)^{|\alpha|}}{\alpha!\beta!\gamma!\delta!}
 {\alpha \choose \alpha'} {\beta \choose \beta'}
 	{p-|\alpha+\beta| \choose \tau}
   }
 	\label{3.3.16} \\
 	 &  & \cdot \ \tilde{G} _{(\alpha',\beta',\tau,\delta)}^{(\gamma)}
 	Z^{|\gamma|}
 	 X^{(p-|\alpha'+\beta'|-\tau-|\gamma+\delta|)}
 	Z^{|\delta|}
 	 \cdot \ B_{\alpha,\beta} ((\f_i \psi_i)^{(r)})_{(\gamma)}^{(\delta)}
   \nonumber
 \end{eqnarray}
where in (\ref{3.3.16}) $Z$ means either an $X$ or a $T$ derivative. That
is, if there are any $T$ vector fields remaining after bringing
 $\tilde{G}_{(\cdot )}^{(\cdot )}$ or $B_{\alpha,\beta (\cdot )}^{\ \ \
(\cdot
 )}$ up to have order zero in $ \Sigma _{r,i, \tilde{G}}^1$, they appear
 explicitely in expressions containing $A$ or $B$.\\
 As for $\Sigma _{r,i, \tilde{G}}^2$ we recall that we actually want to
 estimate the norm of
 \begin{displaymath}
 	V^I \ T^q \ \Sigma _{r,i, \tilde{G}}^2 \ V^K
 \end{displaymath}
 applied to $u$. Our strategy will be to fully exhaust
$ A_{\alpha-\alpha',\beta-\beta'}^{(p-|\alpha'+\beta'|-\tau)}$ and bring
in additional derivatives as well from the vector fields preceding and
following $\Sigma _{r,i, \tilde{G}}^2$. Here, with a change in the
notation used in the preceding sections, we denoted by $V$ either a
derivative along the $X-$ i.e. $X',X''-$direction or along the $Y$
direction.

Thus we may write
\begin{eqnarray}
	\lefteqn{V^I T^q \Sigma_{r,i,\tilde{G}}^{2} V^k}
	\label{3.3.17} \\
	 & = & \underline{C^{N_i}} \tilde{G}_{(\alpha', \beta', \tau, \delta,
	 s_1)}^{(\gamma)} Z^{s_2} Z^{s_3} Z^{s_4} ad_Z^{s_5} \left( B_{\alpha
	 \beta} \left( (\varphi_i \psi_i)^{(r)}\right)_{(\gamma)}^{(\delta)}\right)
	\nonumber
\end{eqnarray}
where $\sum_{i=1}^{5} s_i = |I + K| + p + q - |\alpha' + \beta' | -
\tau$, and, unless $s_3 = 0$, we may take $|s_2 | = |\gamma |$, $|s_4| =
|\delta|$. In any case $|s_2 | \leq |\gamma |$, $|s_4 | \leq |\delta |$
and we may choose to include first any $T$'s present in $Z^{s_2}$,
$Z^{s_4}$, then, when the $T$'s are exhausted we include the possible
$Y$'s present in $Z^{s_2}$, $Z^{s_4}$; only when these are not available
we may include $X$'s in the $Z^{s_2}$, $Z^{s_4}$. At the end of this
process we are left with a term $Z^{s_3}$ of the form $Z^{s_3} = T^{b_1}
V^{b_2}$, where, as above, $V$ denotes either an $X$ or a $Y$ derivative,
$b_1 < q$ provided $q > 0$. We sum up what we did up to now in the
following equality:
\begin{eqnarray}
	\lefteqn{V^I T^q \left[ \tphi p{(\varphi_i \psi_i)^{(r)}} ,
	\tilde{G}\right] V^K}
	\label{3.3.18} \\
	 & = & V^I T^q R_{\left[ \tphi p{(\varphi_i \psi_i)^{(r)}} ,
	\tilde{G}\right] , N_i} V^K
	\nonumber \\
	 &  & + \sum V^I T^q \left( E_{r, i, \tilde{G}}^1 +  E_{r, i,
	 \tilde{G}}^2 \right) V^K
	\nonumber \\
	 &  & + \underline{C^{N_i}} \sum_{{\scriptstyle \tau + |\alpha' + \beta'|
\leq p}
	 \atop {\scriptstyle |\gamma + \delta | \leq N_i}}  \frac{1}{\alpha! \beta !
\gamma ! \delta !} \tilde{G}_{(\alpha', \beta',
	 \tau, \delta, s_1)}^{(\gamma)}   Z^{s_2} Z^{s_3} Z^{s_4}
\nonumber \\
&  & \circ
	 ad_Z^{s_5} \left( B_{\alpha
	 \beta} \left( (\varphi_i \psi_i)^{(r)}\right)_{(\gamma)}^{(\delta)}\right)
	\nonumber
\end{eqnarray}
where $Z^{s_3} = T^{b_1} V^{b_2}$, $b_1 < q$ if $q > 0$ and
$\sum_{i=1}^{5} s_i = |I + K| + p + q - |\alpha' + \beta' | -  \tau$,
$|s_2| \leq |\gamma|$, $|s_4| \leq |\delta |$ (equalities hold unless
$s_3 =0$).

Before proceeding we want to bring all the
$\tilde{G}_{(\cdot)}^{(\cdot)}$ to the left; writing $\tilde{\alpha} =
\alpha - \alpha'$, $\tilde{\beta} = \beta - \beta'$ we get
\begin{eqnarray}
	\lefteqn{V^I T^q \left[ \tphi p{(\varphi_i \psi_i)^{(r)}} ,
	\tilde{G}\right] V^K}
	\label{3.3.19} \\
	 & = & V^I T^q R_{\left[ \tphi p{(\varphi_i \psi_i)^{(r)}} ,
	\tilde{G}\right] , N_i} V^K
	\nonumber \\
	 &  & + \sum_{{ \scriptstyle 1\leq |\alpha'+\beta'+\gamma + \delta | + \tau
\leq p}
	 \atop { \scriptstyle I' \leq I, \ q'\leq q}}(-1)^{|\alpha'|} {\alpha
\choose \alpha'}
	 {\beta \choose \beta'} {|I| \choose |I'|} {q \choose q'}
	 \tilde{G}_{(\alpha', \beta', \tau, \delta, |I'|, q')}^{(\gamma)}
	\nonumber \\
	 &  & \cdot \ T^{|\gamma|} V^{I-I'} T^{q-q'}
	\nonumber \\
  &  & \cdot
  \sum_{|\tilde{\alpha} +
	 \tilde{\beta}| \leq p - |\alpha' + \beta'| - \tau - | \gamma + \delta
	 |} (-1)^{|\tilde{\alpha}|} {{p-|\alpha' + \beta'| - |\tilde{\alpha} +
	 \tilde{\beta} |} \choose \tau}
  \nonumber \\
	 &  &
  \frac{1}{\alpha! \beta! \gamma! \delta!}
  A_{\tilde{\alpha} \tilde{\beta}}^{p - |\alpha' + \beta'| - \tau -
	 | \gamma + \delta |} T^{|\delta|} B_{\tilde{\alpha} + \alpha',
	 \tilde{\beta} + \beta'} \left((\varphi_i
	 \psi_i)^{(r)}\right)_{(\gamma)}^{(\delta)}  V^K
	\nonumber \\
	 &  & + \sum_{{{ \scriptstyle p - \tau - |\gamma + \delta| \leq
|\alpha+\beta| \leq p -
	 \tau } \atop {\scriptstyle \alpha' \leq \alpha, \beta' \leq \beta}} \atop
{ \scriptstyle 1 \leq
	 |\alpha' + \beta' + \gamma + \delta | + \tau \leq p}} \sum_{{I'\leq I}
	 \atop {q'\leq q}} (-1)^{|\alpha|} {\alpha \choose \alpha'} {\beta
	 \choose \beta'} {{p - |\alpha+\beta|} \choose \tau}
	\nonumber \\
  &  & {|I| \choose |I'|}
	 {q \choose q'}
  \frac{1}{\alpha! \beta! \gamma! \delta!}
	 \tilde{G}_{(\alpha', \beta', \tau, \delta, |I'|, q')}^{(\gamma)}
  \nonumber \\
	 &  &
	 \times Z^{|\gamma|} \ T^{q-q'}
V^{I-I'+p-|\alpha'+\beta'+\gamma+\delta|-\tau}
	 Z^{|\delta|} B_{\alpha \beta} \left((\varphi_i
	 \psi_i)^{(r)}\right)_{(\gamma)}^{(\delta)} V^K
	\nonumber \\
	 &  & + \underline{C^{N_i}}
	 \sum_{{{ \scriptstyle |s_1| + |b_1|+|b_2| + |\gamma+\delta| =
	 |I+K|+p+q-|\alpha'+\beta'| -\tau}\atop { \scriptstyle \alpha'\leq \alpha,
\beta'\leq
	 \beta, |\alpha'+\beta'|+\tau\leq p}} \atop { \scriptstyle b_1 < q \ {\rm
if } \ q>0}}
	 \frac{1}{\alpha! \beta! \gamma! \delta!}
	 \tilde{G}_{(\alpha', \beta', \tau, \delta, s_1)}^{(\gamma)}
	\nonumber \\
	 &  & Z^{|\gamma|} T^{b_1} V^{b_2} B_{\alpha \beta} \left((\varphi_i
	 \psi_i)^{(r)}\right)_{(\gamma)}^{(\delta)}
	\nonumber \\
	&  & + \underline{C^{N_i}}
	\sum_{{ \scriptstyle |s_2+s_4|\leq |\gamma+\delta| \leq N_i} \atop
{ \scriptstyle |s_1+s_2+s_4|\leq
	|I+K|+p+q-|\alpha'+\beta'|-\tau}} \frac{1}{\alpha! \beta! \gamma!
	\delta!}  \tilde{G}_{(\alpha', \beta', \tau, \delta, s_1)}^{(\gamma)}
	 \nonumber \\
  &  &
  \cdot Z^{s_2+s_4}  B_{\alpha \beta} \left((\varphi_i
	 \psi_i)^{(r)}\right)_{(\gamma)}^{(\delta)}.
  \nonumber
\end{eqnarray}
The second term on the right needs rewriting if it is to be brought into
a form where the $\Sigma A_{\tilde{\alpha} \tilde{\beta}}
B_{\alpha'+\tilde{\alpha}, \beta'+\tilde{\beta}} ((\f_i \psi_i)^{(r)})
\frac{1}{\tilde{\alpha}! \tilde{\beta}!}
{{p-|\alpha'+\tilde{\alpha}+ \beta' + \tilde{\beta}|} \choose \tau} $
have the right balance, as in (\ref{2.2.1}). First we write
\begin{eqnarray}
	 \lefteqn{  {{p-|\alpha'+\tilde{\alpha}+ \beta' + \tilde{\beta}|}
  \choose \tau} }
  \label{3.3.20} \\
	  & = &\sum_{\ell \leq |\tilde{\alpha}+ \tilde{\beta}|,\tau}
	 (-1)^{\ell} {{p-|\alpha'+ \beta'|- \ell} \choose {\tau - \ell}}
	 {{|\tilde{\alpha}+ \tilde{\beta}|} \choose \ell}
	 \nonumber \\
	 & = & \sum_{\ell \leq |\tilde{\alpha}+ \tilde{\beta}|,\tau}
	 (-1)^{\ell} {{p-|\alpha'+ \beta'|- \ell} \choose {\tau - \ell}}
	 \sum_{ {\scriptstyle |\alpha'' + \beta''| \leq \ell} \atop
{\scriptstyle \alpha'' \leq
	 \tilde{\alpha},\beta'' \leq \tilde {\beta}}
	 } {\tilde{\alpha} \choose \alpha''} {\tilde{\beta} \choose \beta''}.
	\nonumber
\end{eqnarray}
\begin{proposition}
	\label{p3.3.3}
	Any expression of the form
	\begin{displaymath}
		V^\delta \ ad_V^{\tau_1} \ (ad_V^{\tau_2} \ B_{\alpha!\beta!}
		(ad_V^{\tau_3} (\f \psi)^{(r)}))_{(\gamma)}^{(\delta)} \ V^K
	\end{displaymath}
	may be written as a sum of $C^{|\rho_1|}$ terms,
	$|\rho_1| = |\delta + \Sigma \tau_1 + \gamma| +|K|$, each of the form
	\begin{displaymath}
		x^{\tilde{\rho} _1} \ B_{\tilde{a},\tilde{b}}
		 ((\f \psi)^{(|\rho_1| + r)}),
	\end{displaymath}
	where $|\tilde{\rho} _1| \leq |\rho_1| = |\delta + \Sigma \tau_1 + \gamma |
	+ |K|$.
\end{proposition}
\begin{remark}
	This is where the precise form of the vector fields introduced in
	Definition \ref{d2.3.6} is needed. In particular any $x-$ or $\xi-$
	derivative of a symbol can be thourh of as the action of a $Z_\sigma$
	field on the symbol itself. We point out explicitely that tge $W_\sigma$
	and $\Xi_\sigma$ derivatives may be commuted with $B_{\alpha,\beta}$ and
	go directly onto their argument, see e.g. (\ref{2.3.15}).
\end{remark}
The Proposition \ref{p3.3.3} is very easy to prove.
\par\noindent
{\bf Proof.}
Is a repeated use of Lemma \ref{A.3}. An expression of the
form $V^\delta \circ f(x,D)$ may be written as $ C^{|\delta|}$ terms of
the type $ x^{\delta'} W^\delta \circ f(x,D)$, $\delta' \leq \delta$.
Using the definition of the vector fields in Definition \ref{d2.3.6} and
Lemma \ref{lA.3} we prove the result.
\hfill$\blacksquare$
\vskip 1.cm

Using the argument of Proposition \ref{p3.3.3} we may write
\begin{eqnarray}
	 \lefteqn{ T^{|\delta|} \
\left( B_{\tilde{\alpha}+\alpha',\tilde{\beta}+\beta'}
	 ((\f_i \psi_i)^{(r)})_{(\gamma)}^{(\delta)} \right)
	 }
	\label{3.3.21} \\
	 & = & C^{|\gamma+\delta+\alpha'+\delta'|} x^{\rho'}
	 B_{\tilde{\alpha}\tilde{\beta}} \left((\varphi_i
	 \psi_i)^{(r+|\gamma+\delta+\alpha'+\beta'|)}\right)_{(\gamma)}^{(\delta)}
	\nonumber
\end{eqnarray}
where $|\rho'| \leq |\gamma+\delta+\alpha'+\beta'| $.

Let us focus our attention on the second term of (\ref{3.3.19}) and keep
into account (\ref{3.3.20}): the coefficient ${\tilde{\alpha} \choose
\alpha''} {\tilde{\beta} \choose \beta''} \frac{1}{\tilde{\alpha}!
\tilde{\beta}!} = \frac{1}{\alpha''! \underline{\alpha}! \beta''!
\underline{\beta}!}$, where $\underline{\alpha} = \tilde{\alpha} -
\alpha''$, $\underline{\beta} = \tilde{\beta} - \beta''$, makes it appear
that $\underline{\alpha}$, $\underline{\beta}$ should be the running
indices, not $\tilde{\alpha}$, $\tilde{\beta}$. Well,
\begin{equation}
	A_{\tilde{\alpha} \tilde{\beta}}^{(p_1)} = X''^{\alpha''}
	X''^{\underline{\alpha}} X'^{\underline{\beta}} X'^{\beta''} T^{p_1 -
	|\alpha'' + \beta''| - |\underline{\alpha} - \underline{\beta}|},
	\label{3.3.22}
\end{equation}
but now the vector fields $X'^{\beta''}$ are in the wrong position. To
manage this type of terms, where, if we maintain the balance needed for
iteration, then the vector fields are in unadmissible locations, we have
found it helpful to write
\begin{eqnarray}
	\lefteqn{(T^{p_1})_{X'^\sigma, x'^\rho, ad_{X'}^\tau,
	\left((\tilde{\varphi} \tilde{\psi})^{(r)}\right)}}
	\label{3.3.23} \\
	 & = & \sum_{|\alpha_1 + \beta_1| \leq p_1}
	 \frac{(-1)^{|\alpha_1|}}{\alpha_1! \beta_1!}A_{\alpha_1
	 \beta_1}^{(p_1)} X'^\sigma x'^\rho ad_{X'}^\tau B_{\alpha_1
	 \beta_1}\left((\tilde{\varphi} \tilde{\psi})^{(r)}\right).
	\nonumber
\end{eqnarray}
thus the second term on the right of (\ref{3.3.19}) becomes
\begin{eqnarray}
	\lefteqn{\sum_{{\scriptstyle 1\leq|\alpha'+\beta'+\gamma+\delta|+\tau\leq
p}\atop
	{ \scriptstyle I'\leq I, q'\leq q}} (-1)^{|\alpha'|} {\alpha\choose \alpha'}
{\beta
	\choose \beta'} {|I|\choose |I'|} {q\choose q'} \tilde{G}_{(\alpha',
	\beta', \tau, \delta, I', q')}^{(\gamma)}}
	\label{3.3.24} \\
	 &  & \times T^{|\gamma|} V^{I-I'} T^{q-q'}
\sum_{{\scriptstyle |\alpha''+\beta''|=\ell\leq
	 |\tilde{\alpha}+\tilde{\beta}|, \tau}\atop { \scriptstyle \alpha''\leq
	 \tilde{\alpha}, \beta''\leq \tilde{\beta}}} (-1)^\ell
	 {{p-|\alpha'+\beta'|-\ell}\choose {\tau - \ell}}
	\nonumber \\
	 &  & \times\frac{(-1)^{|\alpha''|}}{\alpha''! \beta''!}
	 \underline{C^{|\gamma+\delta+\alpha'+\beta'|}} X''^{\alpha''}
\nonumber \\
&  & \times
	 (T^{p-|\alpha'+\beta'|-|\alpha''+\beta''|-\tau-|\gamma+\delta|})
_{X'^{\beta''},
	 x'^{\rho'}, ad_{X'}^{\alpha''}, \left((\varphi_i
	 \psi_i)^{(r+|\gamma+\delta+\alpha'+\beta'+\beta''|)}\right)}V^K ,
	\nonumber
\end{eqnarray}
where $|\rho'| \leq |\gamma+\delta+\alpha'+\beta'|$.

First of all let us now deal with the $X'^{\beta''}$:
\begin{eqnarray}
	\lefteqn{(T^{p_1})_{X'^{\beta''},
	 x'^{\rho'}, ad_{X'}^{\alpha''}, \left( (\tilde{\varphi}
	 \tilde{\psi})^{(r)}\right)}}
	\label{3.3.25} \\
	 & = & \underline{C^{|\beta''|}}
	 (T^{p_1})_{X'^{\beta''_1}  (x'^{\rho'}),
	 ad_{X'}^{\alpha''+\beta''_2}, \left( (\tilde{\varphi}
	 \tilde{\psi})^{(r)}\right)} X'^{\beta''_3}
	\nonumber
\end{eqnarray}
where $\sum \beta''_i = \beta''$. Now
\begin{displaymath}
	X'^{\beta''_1}  (x'^{\rho'}) = \underline{|\rho'|^{\beta''_1}
	C^{|\beta''_1|}} x'^{\max \{\rho' - \beta''_1, 0\}},
\end{displaymath}
and, as above
\begin{displaymath}
	ad_{X'}^{\alpha''+\beta''_2} \left(B_{\underline{\alpha}
	\underline{\beta}} \left( (\tilde{\varphi}
	 \tilde{\psi})^{(r)}\right) \right) = C^{|\alpha''+\beta''_2|}
	 x'^{\rho_2} B_{\underline{\alpha}
	\underline{\beta}} \left( (\tilde{\varphi}
	 \tilde{\psi})^{(r+|\vartheta|)}\right) ,
\end{displaymath}
where $|\rho_2| \leq |\alpha''+\beta''_2|$, $|\vartheta|\leq |\alpha'' +
\beta''_2|$. Hence
\begin{eqnarray}
	\lefteqn{(T^{p_1})_{X'^{\beta''_1},
	 x'^{\rho'}, ad_{X'}^{\alpha''}, \left( (\varphi_i
	 \psi_i)^{(r+|\gamma+\delta+\alpha'+\beta'+\beta''|)}\right)}}
	\label{3.3.26} \\
	 & = & \underline{C^{|\alpha''+\beta''_2|} |\rho'|^{|\beta''_1|}}
	 (T^{p_1})_{
	 x'^{\rho''+\rho_2}, \left( (\varphi_i
	 \psi_i)^{(r+|\gamma+\delta+\alpha'+\beta'+\beta''|+
	 |\alpha''+\beta''_2| )}\right)} X'^{\beta''_3},
	\nonumber
\end{eqnarray}
where $\sum \beta''_i = \beta''$.

>From the definition we have
\begin{eqnarray}
	\tphi {p_1}{X''_j, (\varphi \psi)^{(r)}} & = & x''_j \tphi
	{p_1}{(\varphi \psi)^{(r)}} + \tphi {p_1-1}{ad_{X'}, (\varphi \psi)^{(r)}}
	\label{3.3.27} \\
	 & = & x''_j \tphi {p_1}{(\varphi \psi)^{(r)}} + \underline{2}
	 \tphi {p_1-1}{X, (\varphi \psi)^{(r+1)}}
	\nonumber \\
	 & = & \cdots
	\nonumber \\
	 & = & \sum_{k=0}^{p_1} \ \underline{C^k} x''_j \tphi {p_1-k}{(\varphi
	 \psi)^{(r+k)}}
	\nonumber
\end{eqnarray}
and
\begin{equation}
	\tphi {p_1}{X'_j, (\varphi \psi)^{(r)}}  =  x'_j \tphi
	{p_1}{(\varphi \psi)^{(r)}} + \tphi {p_1-1}{(\varphi \psi)^{(r+1)}}
	\label{3.3.28}
\end{equation}
so that, for several $X$'s, with $\sum k_i = k$, we obtain
\begin{eqnarray}
	\lefteqn{\tphi {p_1}{X^\varepsilon , (\varphi \psi)^{(r)}}}
	\label{3.3.29} \\
	 & = & \sum_{k_1=0}^{p_1} \sum_{k_2=0}^{p_1-k_1} \cdots
	 \sum_{k_{|\varepsilon| = 0}}^{p_1 - k_1 - \cdots -k_{|\varepsilon|-1}}
	 \underline{C^{k_1+k_2+\cdots + k_{|\varepsilon|}}} x^{\varepsilon'}
	 \tphi {p_1-k}{(\varphi \psi)^{(r+k)}}
	\nonumber \\
	 & = & \sum_{k=0}^{p_1} \underline{C^k {{|\varepsilon| + k - 1} \choose
	 k}} x^{\varepsilon'} \tphi {p_1-k}{(\varphi \psi)^{(r+k)}},
	\nonumber
\end{eqnarray}
where $|\varepsilon'| \leq |\varepsilon|$. Furthermore, since
\begin{eqnarray*}
	\lefteqn{{{p-|\alpha'+\beta'|-|\alpha''+\beta''|} \choose {\tau -
	|\alpha''+ \beta''|}} \frac{1}{\alpha''! \beta''!}}\\
	 & \leq & \frac{C}{\tau!} {\tau \choose {|\alpha''+\beta''|}}
	 \frac{(p-|\alpha'+\beta'| - |\alpha''+\beta''|)!}{(p-|\alpha'+\beta'| -
	 \tau)!}  \\
	 & \leq &  C^\tau \frac{N_i^{\tau-|\alpha''+\beta''|}}{\tau!} ,
\end{eqnarray*}
we may estimate (\ref{3.3.24}) applied to $u$ in $L^2$-norms by:
\begin{eqnarray}
	\lefteqn{ \sup_{ {{{{\scriptstyle 1\leq
|\alpha'+\beta'+\gamma+\delta|+\tau\leq p} \atop
	{ \scriptstyle I'\leq I, q' \leq q, |\alpha''+\beta''|\leq \tau}} \atop
{ \scriptstyle |\rho'|\leq
	|\alpha'+\beta'+\gamma+\delta|}} \atop { \scriptstyle 1\leq
	p-|\alpha'+\beta'+\alpha''+\beta''+\gamma+\delta| - \tau}} \atop
	{\alpha''\leq \alpha', \sum \beta''_i = \beta''\leq \beta'} }
\underline{C^{|\gamma+\delta+\alpha'+\beta'|+\tau+k+q'+|I'|}} }
	\label{3.3.30} \\
&  & \times
	\frac{N_i^{d_3}}{N_i^{d_2} (q'+\tau+|\alpha'+\beta'+\gamma+\delta|+|I|)!}
\nonumber \\
&  & \times
\left\| \tilde{G}_{(\alpha', \beta', \tau, |I'|,
	 q')}^{(\gamma)} T^{|\gamma|} V^{I-I'+\alpha''} T^{q-q'}
	 x^{\rho'+\alpha''+\beta''} \right .
\nonumber \\
	 & \cdot &  \left .  \tphi
	 {p-|\alpha'+\beta'+\alpha''+\beta''+\gamma+\delta| - \tau
	 -k}{(\varphi_i \psi_i)^{(d_2+r)}}  V^{\beta''_3+K} u \right\| ,
	\nonumber
\end{eqnarray}
where $d_3 = \tau+k+ |\gamma+\delta+\alpha'+\beta'+\alpha''+\beta''| +
|I'| +q' - |\alpha'' + \beta''_3|$ and $d_2 =
|\gamma+\delta+\alpha'+\beta'+\alpha''+\beta''+\beta''_2|$. We have
written out $d_3$, $d_2$ since $d_2$ is the number of new derivatives on
$(\varphi_i \psi_i)$ and $d_3$ is the net loss of free derivatives; we
point out that $d_3 \geq 1$. Bringing the powers of $x$ to the left of
$V^{I-I'+\alpha''} T^{q-q'}$ will not alter the form of (\ref{3.3.30})
substantially: it merely adds to the supremum the condition $|\rho''|
\leq |I-I'+\alpha''|$, $|\rho'+\alpha''+\beta''|$, replaces
$V^{I-I'+\alpha''}$ by $V^{I-I'+\alpha''-\rho''}$, moves every $x$'s to
the left of $V^{I-I'+\alpha''}$ and adds a factor
$C^{|\rho'+\alpha''+\beta''|} N_i^{\rho''}$ in front. Thus (\ref{3.3.30})
is bounded by
\begin{eqnarray}
	\lefteqn{ \sup_{ {{\scriptstyle \tilde{I}, \tilde{K}, \tilde{q} \leq q,
\tilde{p} \leq
	p} \atop { \scriptstyle \rho/2, s/2, |d_1+\gamma|\leq \Delta I+\Delta q +
\Delta K +
	\Delta p}} \atop {\scriptstyle \Delta I+\Delta q + \Delta K + \Delta p \geq
1, \Delta
	q \geq 0} } \frac{(CN_i)^{\Delta I +\Delta q +\Delta K +\Delta p}}{N_i^s
	|\gamma+ d_1|!}}
	\label{3.3.31} \\
	 &  & \times \left\| \tilde{G}_{(d_1)}^{(\gamma)} T^{|\gamma|} x^\rho
	 V^{\tilde{I}} T^{\tilde{q}} \tphi {\tilde{p}}{(\varphi_i
	 \psi_i)^{(r+s)}} X^{\tilde{K}} u \right\| ,
	\nonumber
\end{eqnarray}
where $\Delta I = |I| - |\tilde{I}|$, $\Delta p = p - \tilde{p}$, $\Delta
q = q - \tilde{q}$, $\Delta K = |K| - |\tilde{K}|$.
\begin{remark}
	\label{r3.3.4}
	From now on we shall occasionally include a numerical multiple by $N_i$
	as a derivative of a localizing function when this will be enough for
	our purposes. Thus $(\varphi \psi)^{(a+b)}$ could refer to $N^a (\varphi
	\psi)^{(b)}$.
	
	Moreover we shall restrict ourselves to the analytic type
	of estimates, since this will be enough for the propagation of
	regularity theorem.
\end{remark}
Using the fact that
\begin{displaymath}
	{p-|\alpha+\beta| \choose \tau} \leq \frac{N_i^\tau}{\tau!},
\end{displaymath}
\begin{displaymath}
	(|\alpha'+\beta'+\gamma+\delta|+|I'|+q')!\leq
	C^{|\alpha'+\beta'+\gamma+\delta|+|I'|+q'} \alpha'! \beta'! \gamma!
	\delta! |I'|! q'!
\end{displaymath}
and applying Proposition \ref{3.3.3}, the third term in (\ref{3.3.19}) is
bounded by
\begin{eqnarray}
	\lefteqn{ \underline{C^{N_i}} \sup_{ {\scriptstyle |\tilde{I}| + \tilde{q} <
|I| + q
	+ p} \atop {\scriptstyle s < N_i, |d_i+\gamma|\leq N_i, \Delta q \geq 0}}
	\frac{N_i^{\Delta I + \Delta q + p}}{N_i^s |d_1+\gamma|!}}
	\label{3.3.32} \\
	 & \times &  \left\| \tilde{G}_{(d_1)}^{(\gamma)} Z^{|\gamma|}
	 T^{\tilde{q}} V^{\tilde{I}} (\varphi_i \psi_i)^{s+r}\right\| ,
	\nonumber
\end{eqnarray}
again with $\Delta I = |I| - |\tilde{I}|$, $\Delta p = p - \tilde{p}$,
$\Delta q = q - \tilde{q}$, all non negative in this case.

The last two terms in (\ref{3.3.19}) are similarly treated: for the next
to last we get the bound
\begin{eqnarray}
	\lefteqn{ \underline{C^{N_i}} \sup_{ {\scriptstyle |\tilde{I}| \leq |I|,
	\tilde{q}\leq q+ p} \atop {\scriptstyle s < N_i, |d_i+\gamma|\leq N_i}}
	\frac{N_i^{\Delta I +|K| + \Delta q + p}}{N_i^s |d_1+\gamma|!}}
	\label{3.3.33} \\
	 & \times &  \left\| \tilde{G}_{(d_1)}^{(\gamma)} Z^{|\gamma|}
	 T^{\tilde{q}} V^{\tilde{I}} (\varphi_i \psi_i)^{s+r}\right\|.
	\nonumber
\end{eqnarray}
As for the last term, using also Remark \ref{r3.3.4}, we have the bound
\begin{eqnarray}
	\lefteqn{ \underline{C^{N_i}} \sup_{ {{\scriptstyle |s_2|\leq |\gamma|,
|s_4|\leq
	|\delta|} \atop { \scriptstyle |d_1+\gamma|, |\gamma+\delta| \leq N_i}}
\atop { \scriptstyle |s|\leq
	N_i} }  \frac{N_i^{|\gamma|-|s_2|} N_i^K}{N_i^s
	N_i^{|\gamma|-|s_2|+|\delta|-|s_4|}} \frac{1}{\alpha!\beta!\gamma!\delta!}}
	\label{3.3.34} \\
	 & \cdot & \left\| \tilde{G}_{(d_1)}^{(\gamma)} Z^{s_2} (\varphi_i
	 \psi_i)^{(r+s)} V^K u \right \|
	\nonumber \\
	 & \leq & \underline{C^{N_i}}
	 \sup_{{\scriptstyle |s_2|\leq |\gamma|, s\leq N_i} \atop {
\scriptstyle |d_1+\gamma| \leq N_i}}
	 \frac{N_i^{|I|+|K|+p+q}}{N_i^s |d_1+\gamma|!} N_i^{|\gamma| - |s_2|}
	\nonumber \\
	 & \cdot & \left \| \tilde{G}_{(d_1)}^{(\gamma)} Z^{s_2} (\varphi_i
	 \psi_i)^{(r+s)} u \right\|,
	\nonumber
\end{eqnarray}
where the factor $N_i^{|\gamma| - |s_2|}$ will later go with
$\tilde{G}_{(d_1)}^{(\gamma)} Z^{s_2}$.

To bring $\tilde{G}_{(d_1)}^{(\gamma)} Z^{s'}$, $|s'| \leq |\gamma|$, out
of the $L^2$-norm we must recall that the symbol of $\tilde{G}$ has been
cut-off to be zero for $|\xi|\leq N_i$:
\begin{eqnarray*}
	\sigma\left(\tilde{G}_{(d_1)}^{(\gamma)} \right)\xx & = &
	\partial_x^{d_1}\partial_\xi^\gamma \left(\Phi(x) g\xx \Psi'_i(\xi)\right)  \\
	 & = & g_{(d_1)}^{(\gamma)}\xx \Psi'_i(\xi) + \sum_{1\leq |\gamma_1|\leq
	 |\gamma|}  g_{(d_1)}^{(\gamma - \gamma_1)}\xx \Psi'^{(\gamma_1)}_i(\xi) ,
\end{eqnarray*}
when $x\in\Omega$, since $\Phi(x) \equiv 1$ in a neighborhood of
$\bar{\Omega}$.

Since $Z^{s'}$ is at most $C^{|\gamma|} {|s'|\choose |\gamma'|}
|\gamma'|! \leq \tilde{C}^{|\gamma|} |\gamma'|!$ terms of the form
$D^{s'-\gamma'} x^{s'-\gamma'}$, for some $\gamma' \leq s'$, we may write
\begin{eqnarray*}
	\lefteqn{\left\|\op \left( g_{(d_1)}^{(\gamma)}\xx \Psi'_i(\xi) \right)
	D^{s'-\gamma'} w\right\|} \\
	 & \leq & \sup_{{\scriptstyle |\rho|\leq n+1}\atop \scriptstyle |\xi|}
\left|
	 g_{(d_1+\rho)}^{(\gamma)}\xx \xi^{s'-\gamma'}\Psi'_i(\xi) \right| \|w\|  \\
	 & \leq & \sup_{{\scriptstyle |\rho|\leq n+1}\atop \scriptstyle |\xi|}
\left[
	 C^{|d_1+\rho+\gamma|} (d_1+\rho)! |\gamma|! \frac{|\xi|^{|s'-\gamma|}
	 \Psi'_i(\xi)}{(1+|\xi|)^{|\gamma|}} \right] \|w\|  \\
	 & \leq & C C^{|d_1+\gamma|} \frac{d_1! (\gamma-s')!
	 s'!}{N_i^{|\gamma-s'|+|\gamma'|}} ,
\end{eqnarray*}
so that
\begin{equation}
	\left\| \op \left( g_{(d_1)}^{(\gamma)}\xx \Psi'_i(\xi) \right)
	\right\|_{L^2\rightarrow L^2} \leq C C^{|d_1+\gamma|} (d_1+s')!
	\label{3.3.35}
\end{equation}
if $w$ has support in a compact set (so $|x|^{|s'-\gamma'|} \leq
C^{|s'-\gamma'|} \leq C^{|\gamma|}$) and $\xi$ may be taken to lie in
$\tilde{\Gamma}_i$, which we will be able to do since $w$ will contain
$\varphi_i \psi_i$. Again, for sake of simplicity we are assuming that
$g$ (and thus $P$) has ``analytic coefficients".

As for the case $|\gamma_1| \geq 1$ (i.e. when some derivatives land on
$\Psi'_i(\xi)$) the expression $g_{(d_1)}^{(\gamma - \gamma_1)}\xx
\Psi'^{(\gamma_1)}_i(\xi) \xi^{s'-\gamma'} $ has $\xi$-support in $\{
N_i \leq |\xi| \leq 2 N_i\}$ and
\begin{eqnarray*}
	\lefteqn{\sup_{{\scriptstyle |\rho|\leq n+1}\atop {\scriptstyle N_i \leq
|\xi|\leq 2 N_i}}
	\left| \partial_\xi^\rho \left( g_{(d_1)}^{(\gamma - \gamma_1)}\xx
\Psi'^{(\gamma_1)}_i(\xi) \xi^{s'-\gamma'}\right)\right|} \\
	 & \leq &   \sup_{{\scriptstyle |\rho|\leq n+1}\atop \scriptstyle |\xi|}
C^{|d_1+\gamma|+n+1} d_1!
	 (|\gamma-\gamma_1| + |\rho|)!
	 \frac{|\xi|^{|s'-\gamma'|}}{(1+|\xi|)^{|\gamma-\gamma_1|}} \chi_{\{N_i
	 \leq |\xi|\leq 2 N_i\}}   \\
	 & \leq  & C C^{|d_1+\gamma|} \sup_{{\scriptstyle |\rho|\leq n+1,
	 \gamma'\leq s'}\atop { \scriptstyle N_i \leq |\xi|\leq 2 N_i}} \frac{d_1!
	 (|\gamma-\gamma_1| + |\rho|)!
	 }{(1+|\xi|)^{|\gamma-s'+\gamma'-\gamma_1|}}  \\
	 & \leq & C C^{|d_1+\gamma|} \frac{d_1!}{N_i^{|\gamma'-s'|}},
\end{eqnarray*}
so that for $|\gamma| \geq 1$, due to the fact that
$\frac{|\gamma|^{|\gamma'|}}{N_i^{|\gamma'|}} \leq 1$,
\begin{eqnarray}
	\lefteqn{ \left\|\op \left( g_{(d_1)}^{(\gamma-\gamma_1)}\xx
	\Psi'^{(\gamma_1)}_i(\xi)\right) Z^{s_1} w \right\|}
	\label{3.3.36} \\
	 & \leq  & C C^{|d_1+\gamma|} d_1! N_i^{|s'|} \sup_{\gamma'\leq s'}
	 \left\| \op \left( \chi_{\{N_i
	 \leq |\xi|\leq 2 N_i\}} \right)x^{s'-\gamma'} w \right\|.
	\nonumber
\end{eqnarray}
Since in our case $w = T^{\tilde{q}} V^{\tilde{I}}
(\varphi_i\psi_i)^{(r+s)} V^{\tilde{K}} u $, so that the conic support of
$w$ lies in $\tilde{\Gamma}_i$, we write $x^{s'-\gamma'} T^{\tilde{q}}
V^{\tilde{I}}  (\varphi_i\psi_i)^{(r+s)} V^{\tilde{K}} $ with all the
$V$'s at the left:
\begin{displaymath}
	x^{s'-\gamma'} T^{\tilde{q}}  V^{\tilde{I}}  (\varphi_i\psi_i)^{(r+s)}
	V^{\tilde{K}} = \underline{C^{\tilde{K}}} \  x^{s'-\gamma'} T^{\tilde{q}}
V^{\tilde{I} + K' }  (\varphi_i\psi_i)^{(r+s+|\tilde{K}-K'|)},
\end{displaymath}
for some $K'\leq \tilde{K}$, and each of these terms is bounded by
\begin{displaymath}
	N_i^{\rho'} D^{\tilde{q}+\tilde{I}+K'-\rho'}x^{s'-\gamma'+|\tilde{I}+K'|
	-\rho'} (\varphi_i \psi_i)^{(r+s+|\tilde{K}-K'|)},
\end{displaymath}
with $\rho' \leq \min \{ \tilde{q}+\tilde{I}+K',
s'-\gamma'+\tilde{I}+K'\}$. Similarly, using the definition of $\tphi
p{(\varphi_i\psi_i)^{(r+s)}}$, $x^{\rho'} T^{\tilde{q}} V^{\tilde{I}}
\tphi p{(\varphi_i\psi_i)^{(r+s)}} V^{\tilde{K}}$ is at most $C^{N_i}$
terms, each less then $N_i^{\rho'}$ terms of the form
\begin{displaymath}
	\frac{1}{\alpha!\beta!} D^{\tilde{q}+\tilde{I}+\tilde{p}+K'-\rho'}
	x^{\rho+\tilde{p}+\tilde{K}+\tilde{I}-\rho'}
	(\varphi_i\psi_i)^{(|\alpha+\beta|+r+s+|\tilde{K}-K'|)},
\end{displaymath}
for some $\alpha$, $\beta$, $\rho'$, $K'$ with $|\alpha+\beta|\leq
\tilde{p}$, $K'\leq \tilde{K}$ and
\begin{displaymath}
	|\rho|\leq \min \{ \tilde{q}+|\tilde{I}|+\tilde{p}+|K'|,
	|\rho|+|\tilde{I}|+\tilde{p}+|\tilde{K}|\}.
\end{displaymath}
Thus
\begin{eqnarray*}
	\lefteqn{\left\| \chi_{\{N_i
	 \leq |\xi|\leq 2 N_i\}} x^{s'-\gamma'} T^{\tilde{q}} V^{\tilde{I}}
	 (\varphi_i \psi_i)^{(r+s)} V^{\tilde{K}} u \right\|} \\
	 & \leq &  \sup_{K'\leq \tilde{K}} N_i^{\tilde{q} +|\tilde{I} + K'|}
	 \left\| (\varphi_i \psi_i)^{(r+s+|\tilde{K} - K'|)} u \right\|
\end{eqnarray*}
and
\begin{eqnarray}
	\lefteqn{\left\| \chi_{\{N_i
	 \leq |\xi|\leq 2 N_i\}} x^{\rho} T^{\tilde{q}} V^{\tilde{I}}
	 \tphi p{(\varphi_i \psi_i)^{(r+s)}} V^{\tilde{K}} u \right\|  }
	\label{3.3.37} \\
	 & \leq & C^{N_i} \sup_{{\alpha \beta}\atop {K'\leq \tilde{K}}}
	 N_i^{\tilde{q}+\tilde{p}+|\tilde{I}+K'|-|\alpha+\beta|}
	 \left\| (\varphi_i \psi_i)^{(r+s+|\tilde{K} - K'|+
	 |\alpha+\beta|)} u \right\| .
	\nonumber
\end{eqnarray}
Let us now use (\ref{3.3.35}), (\ref{3.3.37}) in (\ref{3.3.31}) --
(\ref{3.3.32}). Thus (\ref{3.3.31}) can be bounded by
\begin{eqnarray}
	\lefteqn{\frac{C (C N_i)^{|I-\tilde{I}|+|K-\tilde{K}|+p-\tilde{p}+
	q-\tilde{q}+\sigma |d_1+\gamma|}}{N_i^s N_i^{|d_1+\gamma|}} }
	\label{3.3.38} \\
	 & \cdot & \left\{  \left\| T^{\tilde{q}} V^{\tilde{I}}
	 \tphi {\tilde{p}}{(\varphi_i \psi_i)^{(r+s)}} V^{\tilde{K}} u
	 \right\|\right.
	\nonumber \\
	 &  & + \left . C^{N_i} \sup_{|\alpha+\beta|\leq \tilde{p}}
	 N_i^{\tilde{q}+\tilde{p} +|\tilde{I}+K'|-|\alpha+\beta|}
	 \left\| (\varphi_i \psi_i)^{(r+s+|\tilde{K} - K'|+
	 |\alpha+\beta|)} u \right\| \right\}
	\nonumber \\
	 & \leq & C \frac{(C N_i)^{|I-\tilde{I}|+|K-\tilde{K}|+p-\tilde{p}+
	q-\tilde{q}}}{N_i^s} \left\|  T^{\tilde{q}} V^{\tilde{I}}
	 \tphi {\tilde{p}}{(\varphi_i \psi_i)^{(r+s)}} V^{\tilde{K}} u \right\|
	\nonumber \\
	 &  & + \frac{C^{N_i} N_i^{|I|+|K| +
	 +p+q+r}}{N_i^{|s'|}}  \sup_{|s'|\leq 2 N_i} \left\| (\varphi_i
	 \psi_i)^{(r+s')} u \right\| .
	\nonumber
\end{eqnarray}
Next, (\ref{3.3.32}) or the third term in (\ref{3.3.19}) is similarly
bounded by
\begin{eqnarray}
	\lefteqn{ C^{N_i}
	\frac{N_i^{|I-\tilde{I}|+|K-K'|+p+q-\tilde{q}+|\gamma|-s'}}{N_i^s
	N_i^{|d_1+\gamma|}} N_i^{|d_1+\gamma|} }
	\label{3.3.39} \\
	 & \cdot & \left\{ \left\| T^{\tilde{q}} V^{\tilde{I}} (\varphi_i
	 \psi_i)^{(r+s)} V^{\tilde{K}} u \right \| + C^{N_i}
	 N_i^{q+|\tilde{I}+K'|} \left\| (\varphi_i
	 \psi_i)^{(r+s+|\tilde{K}-K'|)} u \right\| \right \}
	\nonumber \\
	 & \leq & C^{N_i} \frac{N_i^{|I-I'|+|K-K'|+p+q-\tilde{q}}}{N_i^s}
	 \left\| T^{\tilde{q}} V^{\tilde{I}} (\varphi_i \psi_i)^{(r+s)}
	 V^{\tilde{K}} u \right\|
	\nonumber \\
	 &  & + \frac{N_i^{|I|+|K|+p+q+r}}{N_i^{|s'|}} \sup_{|s'|\leq 2 N_i}
	 \left\| (\varphi_i \psi_i)^{(r+s')} u \right\|.
	\nonumber
\end{eqnarray}
Thus (\ref{3.3.19}) may be rewritten, using (\ref{3.3.11}) and taking
into account the bounds under the suprema in (\ref{3.3.31}) --
(\ref{3.3.34}), as
\begin{proposition}
	\label{p3.3.5}
	Let $|I|+p+q \leq N_i$, $|K|\leq 2$; then
	\begin{eqnarray}
		\lefteqn{\left\| V^I T^q \left[ \tphi p{(\varphi_i\psi_i)^{(r)}} , \
		\tilde{G} \right] V^K u \right\|}
		\label{3.3.40} \\
		 & \leq & \left\| V^I T^q R_{\left[ \tphi p{(\varphi_i\psi_i)^{(r)}} , \
		\tilde{G} \right] , N_i} V^K u \right\|
		\nonumber \\
		 &  & + C \sup_{{{{ \scriptstyle \Delta q + \Delta p + \Delta I + \Delta K
\geq 1}
		 \atop {\scriptstyle \Delta p > 0}} \atop {\scriptstyle s\leq 2 (\Delta I +
\Delta + \Delta p +
		 \Delta q)}} \atop { \scriptstyle -\Delta L = |\tilde{L}|-|K| \leq \Delta
p} }
		 C^{|\Delta|} N_i^{|\Delta|} \left\| V^{\tilde{I}} T^{\tilde{q}}
		 \tphi {\tilde{p}}{(\varphi_i\psi_i)^{(r+s)}} V^{\tilde{L}} u \right\|
		\nonumber \\
		 &  & + C^{N_i} \sup_{{{\scriptstyle \tilde{p} = 0, \Delta q + p \geq 0}
\atop
		 {\scriptstyle \Delta I + \Delta L + \Delta q + \Delta p \geq 1}} \atop
{ \scriptstyle s\leq N_i}
		 } N_i^{|\Delta|} \left\| V^{\tilde{I}} T^{\tilde{q}} (\varphi_i
		 \psi_i)^{(r+s)} V^{\tilde{L}} u \right\|,
		\nonumber
	\end{eqnarray}
	where $\Delta I = |I|-|\tilde{I}|$, $\Delta L = |K| - |\tilde{L}|$,
	$\Delta p = p - \tilde{p}$, $\Delta q = q - \tilde{q}$, $\Delta = \Delta
	I + \Delta L + \Delta p + \Delta q$, $|\Delta| = |\Delta I|+|\Delta L|+
	|\Delta p |+ |\Delta q|$.
\end{proposition}
We remark that $\tilde{L}$ may be large and hence $\Delta L$ may be
negative.

We have
\begin{proposition}
	\label{p3.3.6}
	We may take $\tilde{L} = 0$ on the right hand side of (\ref{3.3.40}).
\end{proposition}
{\bf Proof.}
The idea is to commute the $V$ fields behind $\tphi p{(\varphi_i
\psi_i)^{(r+s)}}$; we need only to be concerned with the $X$ vector
fields; in fact the result of commuting back a $Y$ vector field will give
a term analogous to the last in (\ref{3.2.8}). These terms contain
(positive order) derivetives of $\varphi_i^{\#}$ and will be estimated
easily exploiting the assumption on the WF of u in the leaves of the
characteristic manifold.

As for the $X$ vector fields, it suffices to make the following two
remarks: a) for a bounded number of $X$'s we may commute them back to the
left of $\tphi p{(\varphi_i\psi_i)^{(r)}}$ introducing a new constant and
replacing $p$ by $p - \ell$ and $s$ by $s+\ell$, $\ell$ less or equal
than the number of commuted fields. Doing so we do not alter anything
else except the $C$ in $C^{N_i}$, since the second term in (\ref{2.2.8})
(Proposition \ref{2.2.5}) is readily absorbed in the next to last term of
(\ref{3.3.40}).

b) For more $X$'s (there are at least $|\beta''_3| \leq |\beta''| \leq p$
of them) we use Proposition \ref{2.2.5} ensuring that only terms free of
$\tphi p{(\varphi_i \psi_i)^{(r)}}$, $p > 0$, occur when we commute
$X'^k$ with $\tphi p{(\varphi_i \psi_i)^{(r)}}$ and these terms are
readily absorbed in the next to last term of (\ref{3.3.40}).
\hfil$\blacksquare$
\vskip 1.cm
Let us now go back to (\ref{3.3.19}); we claim, for the remainder term,
that
\begin{equation}
	\left\| V^I T^q R_{\left[ \tphi p{(\varphi_i\psi_i)^{(r)}} , \
		\tilde{G} \right] , N_i} V^K u \right\|
		\leq \left(C K_i\right)^{N_i} N_i^{|I|+p+q+r},
	\label{3.3.41}
\end{equation}
provided that $|I|+p+q \leq N_i$, $r \leq 2 N_i$.

Recall that $\sigma(\tilde{G})\xx = \Phi_i(x) g\xx \Psi_i(\xi)$, where
$\Phi_i (x) \equiv 1 $ in a neighborhood of $\bar{\Omega}$ and satisfies
the inequalities $|D^\alpha \Phi_i(x)| \leq (C N_I)^{|\alpha|}$,
$|\alpha| \leq 3 N_i$. We have the
\begin{lemma}
	\label{l3.3.7}
	Let $a_{\alpha \beta} \left((\varphi_i \psi_i)^{(s)}\right) =
	\sigma\left( B_{\alpha \beta} \left( (\varphi_i
	\psi_i)^{(s)}\right)\right)$. Then
	\begin{eqnarray}
		\lefteqn{\left | \left(\frac{\partial}{\partial x}\right)^\rho  \left(
		a_{\alpha \beta}
		\left((\varphi_i\,\psi_i)^{(s)}\right)\right)_{(\gamma)}^{(\delta)} \xx
		\right | }
		\label{3.3.42} \\
		 & \leq & (C \tilde{K}_i)^{s+|\alpha+\beta+\rho+\gamma+\delta|}
		 N_i^{s+|\alpha+\beta+\rho+\gamma|}
		 \left(\frac{N_i}{|\xi|}\right)^{|\delta|}.
		\nonumber
	\end{eqnarray}
\end{lemma}
{\bf Proof.}
We may write
\begin{eqnarray*}
	 \lefteqn {\left(\frac{\partial}{\partial x}\right)^\rho  \left(
		a_{\alpha \beta}
		\left((\varphi_i\,\psi_i)^{(s)}\right)\right)_{(\gamma)}^{(\delta)} }\\
	 & = & \sum_{r\leq |\alpha+\beta|, |\rho+\gamma|}
	 \left(\frac{\partial}{\partial x}\right)^{r+\rho}
	 \left(\frac{\partial}{\partial \xi}\right)^\delta ad_{X'}^\alpha
	 ad_{X''}^\beta \left((\varphi_i \psi_i)^{(s)}\right)  \\
	 & = &  \sum_{r\leq |\alpha+\beta|, |\rho+\gamma|}
	 \underline{C^{|\alpha+\beta|} N_i^r} \
	 \left(\frac{\partial}{\partial x}\right)^{r+\rho}
	 \left(\frac{\partial}{\partial \xi}\right)^\delta
	 x^{|\alpha+\beta|-r} \left((\varphi_i
	 \psi_i)^{(s+|\alpha+\beta|-r)}\right) \\
	 & = & \sum_{r\leq |\alpha+\beta|, |\rho+\gamma|}
	 \underline{C^{|\alpha+\beta+\rho+\gamma|} N_i^r} \
	 \left(\frac{\partial}{\partial \xi}\right)^\delta
	 x^{|\alpha+\beta|-r} \left((\varphi_i
	 \psi_i)^{(s+|\alpha+\beta+\rho+\gamma|-r)}\right),
\end{eqnarray*}
and the result follows from Proposition \ref{2.3.8}.
\hfil$\blacksquare$
\vskip 1.cm

To estimate the remainder term, we write
\begin{eqnarray*}
	\lefteqn{V^I T^q R_{\left[ \tphi p{(\varphi_i\psi_i)^{(s)}} , \
		\tilde{G} \right] , N_i} V^J u} \\
	 & = &  \sum_{{r\leq |I|+p} \atop {r'\leq| J|}} \frac{C^{N_i}
	 N_i^{r+r'}}{\alpha!\beta!\gamma!\delta!} x^{\leq |I|+P+q} D^{|I|+P+q-r}  \\
	 &  & \left(F_{\alpha \beta, \gamma}^\delta H_{\alpha \beta ,
	 \delta}^\gamma -  \left\{ F_{\alpha \beta, \gamma}^\delta \circ
	 H_{\alpha \beta, \delta}^\gamma \right\}_{N_i - |\delta+\gamma|}\right)  \\
	 &  & D^{|J|-r'} x^{|J|-r'} u ,
\end{eqnarray*}
where either
\par\noindent
(i) \qquad $F_{\alpha\beta, \gamma}^\delta =
\tilde{G}_{(\delta)}^{(\gamma)}$, $H_{\alpha\beta, \delta}^\gamma =
\left( B_{\alpha\beta} (\varphi_i
\psi_i)^{(s)}\right)_{(\delta)}^{(\gamma)}$,
\par\noindent
or
\par\noindent
(ii)  \qquad $F_{\alpha\beta, \gamma}^\delta =
\left( B_{\alpha\beta} (\varphi_i
\psi_i)^{(s)}\right)_{(\delta)}^{(\gamma)}$, $H_{\alpha\beta, \delta}^\gamma =
\tilde{G}_{(\delta)}^{(\gamma)}$,
\par\noindent
with $\alpha=\beta=0$ in case (ii). In all cases $p+q+|I| \leq N_i$,
$s\leq N_i$ and $|\gamma+\delta|\leq N_i$, $|\alpha+\beta|\leq p$.

The claim is proved by an application of Lemma \ref{lA.1}.
\hfil$\blacksquare$
\vskip 1.cm
Summing up we have achieved the proof of the result of this section:
\begin{proposition}
	\label{p3.3.8}
	Assume that $u$ verifies the hypotheses of Theorem \ref{th0.1} (with
	$s=1$ for sake of simplicity). Let $\tilde{G}\xd$ be an analytic
	pseudo -- differential operator of degree 0, $|I'|+p+q\leq N_i$, $r\leq
	N_i$, $|J| \leq 2$. Then there exists a positive constant $C =
	C(\tilde{G})$, independent of $N$, such that
	\begin{eqnarray}
		\lefteqn{\left\|  V^{I'} T^q \left[ \tphi p{(\varphi_i\psi_i)^{(r)}} , \
		\tilde{G} \right] V^J u \right\|}
		\label{3.3.43} \\
		 & \leq & C \sup_{{\scriptstyle \Delta\geq 0, \Delta p \geq 0, \Delta q
\geq 0}
		 \atop {\scriptstyle s \leq 2 \Delta, |\tilde{I}|-|I'| \leq \Delta p}}
		 \frac{C^{|\Delta |} N_i^{|\Delta |}}{N_i^s} \left\| V^{\tilde{I}}
		 T^{\tilde{q}} \tphi{\tilde{p}}{(\varphi_i \psi_i)^{(r+s)}} u \right\|
		\nonumber \\
		 &  & + C^{N_I} \sup_{{\scriptstyle \Delta q \geq 0, \Delta \geq 1} \atop
{ \scriptstyle s \leq
		 N_i}} N_i^{\Delta - s} \left\| V^{\tilde{I}}
		 T^{\tilde{q}} (\varphi_i \psi_i)^{(r+s)} u \right\|
		\nonumber \\
		 &  & + C C^{N_i} N_i^{|I'|+p+q+r} K_i^{N_i+r} \|u\| ,
		\nonumber
	\end{eqnarray}
	where $\Delta = \Delta I + \Delta p + \Delta q$, $\Delta I =
	|I'|+|J|-|\tilde{I}|$, $\Delta p = p - \tilde{p}$, $\Delta q = q -
	\tilde{q}$, $|\Delta | = |\Delta I |+ |\Delta q|+ \Delta p$ and $V$
	denotes, as usual in this section, a derivative either in the $X$ or in
	the $Y$ directions.
\end{proposition}

\subsection{Reducing the order by half and the end of the proof of Theorem
2.1}
\label{s3.4}
\setcounter{equation}{0}
\setcounter{theorem}{0}
\setcounter{proposition}{0}
\setcounter{lemma}{0}
\setcounter{corollary}{0}
\setcounter{definition}{0}
\setcounter{remark}{0}

Combining (\ref{3.2.10}), Lemma \ref{l3.2.2} and Proposition \ref{p3.3.8}
we obtain the estimate
\begin{eqnarray}
	\lefteqn{
    \sup_{{\scriptstyle I=I'+J} \atop {\scriptstyle |J| \leq 2}}
    \Vert Z^I T^q	(T^p)_{(\f_i \psi_i)^{(r)}} u \Vert _{L^2(\Omega)}
	}
	\label{3.4.1}\\
	& \leq & C \left\{
  \Vert Z^{I'} T^q (T^p)_{(\f_i \psi_i)^{(r)}} \Phi_i
  Pu \Vert _{L^2(\Omega)}
\phantom{\sup_{{|I'|}\atop {q'' = q+1}} }
	 \right .
	\nonumber \\
	 &  & +  \sup_{ {{\scriptstyle \Delta q \geq 0,\Delta p \geq 0}\atop
{\scriptstyle \Delta \geq 1}}
	 \atop {\scriptstyle s \leq 2 \Delta } }
	 C^{|\Delta|} N_i^{\Delta - s} \Vert Z^{\tilde {I}} T^{\tilde {q}}
	 	(T^{\tilde{p}})_{(\f_i \psi_i)^{(r+s)}} u \Vert _{L^2(\Omega)}
	\nonumber \\
	 &  & + C^{N_i} \sup_ { {\scriptstyle \Delta q \geq 0,\Delta p = p} \atop
{\scriptstyle s \leq N_i} }
	 N_i^{\Delta - s} \Vert Z^{\tilde{I}} T^{\tilde{q}} (\f_i \psi_i)^{(r+s)}
	 u \Vert _{L^2(\Omega)}
	\nonumber \\
	 &  & + C^{N_i}N_i^{|I'|+p+q+r} K_i^{N_i+r} \Vert u \Vert _{L^2(\Omega)}
	\nonumber \\
	 &  & + \left . \underline{|I'|} \sup_{{ \scriptstyle |I'|}\atop {
\scriptstyle q'' = q+1}}
	 \Vert Z^{I'} T^{q''} (T^p)_{(\f_i \psi_i)^{(r)}} u \Vert _{L^2(\Omega)}
	 \right\},
	\nonumber
	\end{eqnarray}
where $Z$ denotes as usual either a $X-$ or a $Y-$ derivative.\\
We point out explicitly that in deriving (\ref{3.4.1}) we used the
assumption that $\tilde{\f}_i \f_i^{\sharp (\alpha)} u$ is analytic, if
 $|\alpha| > 0$, near our base point $\rho$.

 Now starting with $i=q=r=I'=0$, $p=p_0$ we use estimate (\ref{3.4.1})
 with $i=0$ repeatedly. Each time we resubject the 2nd, 3rd and 5th terms
 in (\ref{3.4.1}) to (\ref{3.4.1}) to reduce $|I|+p+q$ (in the 3rd term
 $p=0$ and it cannot be subjected to (\ref{3.4.1}) again).\\
 The aim is to obtain only the first four terms in the right hand side of
  (\ref{3.4.1}), which will happen eventually. That is, we claim that the
  last term in (\ref{3.4.1}) will eventually disappear. To see this,
  observe what happens after each iteration: in the last term $q$ has
  risen by $1$, but $|I|$ has dropped by $2$. After at most $\frac{N}{2}$
  iteration, every term will either contain $Pu$ or, we claim, have $p=0$
  and have at most $1$ free $Z$, with $q$ at most equal to
  $\frac{(p_0+2)}{2}$. To see that this is the case, note that in
  Proposition \ref{p3.3.8}, $p$ decreases in each term on the right,
  while in (\ref{3.2.8}), $p$ may keep it value. Now these terms arise
  only if $|I'|+q > 0$. Thus a given application of Proposition \ref{p3.3.8}
  may reduce $|I'|$ or $p$ or $q$, but once $|I'|+q$ reaches $0$, $p$ must
  decrease.\\
  The value of $q$ need not decrease, and may rise (via $[X,X]=T$) to one
  half the original value of $|I'|+2+p$ or $\frac{p_0+2}{2}$ (new $Z'$s
  arise in Proposition \ref{p3.3.8} - first term on the right hand side -
  only by a corresponding decrease of $p$).\\
  Thus we have proved:
  \begin{proposition}
  	\label{p3.4.1}
  	For $p \leq N_i$,
  	\begin{eqnarray}
  		\lefteqn{
     \sup_{|J|\leq2} \Vert Z^J (T^p)_{(\f_i \psi_i)} u
    		\Vert _{L^2(\Omega)}
      }
  			\label{3.4.2} \\
  		 & \leq &  C^{N_i} \left\{
  		 \sup_{s\leq N_i+2\Delta}  N_i^{\Delta - s}\Vert Z^{\tilde {I}}
  		 T^{\tilde {q}} (T^{\tilde{p}})_{(\f_i \psi_i)^{(s)}} \Phi_i Pu
  		 \Vert _{L^2(\Omega)}
  		 \right .
  		\nonumber \\
  		 &  & + \sup_{{\scriptstyle s\leq2N_i,2\Delta}\atop {
\scriptstyle |J|\leq 2}}
  		 N_i^{\Delta - s} \Vert Z^J (\f_i \psi_i)^{(s)} T^{\tilde {q}}
  		 u \Vert _{L^2(\Omega)}
  		\nonumber \\
  		 &  &    +  \left . \phantom{\sup_{N=0} N^N}  N_i^{p-s} K_i^{N_i} \Vert
u
\Vert
    _{L^2(\Omega)}
  		 \right\},
  		\nonumber
  	\end{eqnarray}
where $K_i$ satisfies the bounds (\ref{1.1.1}) - (\ref{1.1.8}) and
  	$\Delta = p-\tilde{p} + q + \tilde{q}- |\tilde{I}|$.
  \end{proposition}

\begin{lemma}
	\label{l3.4.2}
	For $|I|+p+q\leq N_i, s\leq 2N_i$
	\begin{displaymath}
		\frac{1}{N_i^s} \Vert X^I T^q (T^p)_{(\f_i \psi_i)^{(s)}}
		\Phi_i Pu \Vert _{L^2(\Omega)} \leq C_{{\mit f}}^{N_i} K_i^{N_i}
		N_i^{|I|+p+q}
	\end{displaymath}
\end{lemma}
{\bf Proof.}
Clear due to the assumptions of Theorem \ref{th0.1} and by Proposition
\ref{p2.3.8}.
\hfill$\blacksquare$
\vskip 1.cm

For the second term on the right in Proposition \ref{p3.4.1} we pass to a
new pair of localizing functions, $\f_{i'}$,$\psi_{i'}$, where $i'$ is the
largest integer such that $q+2\leq N_{i'}$ (so that $\tilde{q}!^{-1}) \leq
C^{N_i} N_i^{-\tilde{q}}$.
\begin{proposition}
	\label{p3.4.3}
	Let $|J|\leq 2$ and $s\leq 2N_i$, $\tilde{q}\leq N_i$; then
	\begin{equation}
\label{3.4.3}
	 \frac {1}{N_i^s} \Vert Z^J (\f_i \psi_i)^{(s)} (T^{\tilde{q}} -
		(T^{\tilde{q}})_{\f_i' \psi_i'}) u \Vert  _{L^2(\Omega)}
	\leq  C^{N_i+1} N_i^{\tilde{q}} .
	\end{equation}
\end{proposition}
{\bf Proof.}
Arguing exactly as in the proof of Proposition \ref{2.3.1}.
\hfill$\blacksquare$

\begin{proposition}
	\label{p3.4.4}
	For $S\leq 2N_i$,
\begin{eqnarray*}
	\lefteqn{ \sup_{|J|\leq 2} \Vert Z^J (\f_i \psi_i)^{(s)}
	(T^{\tilde{q}})_{\f_i' \psi_i'}) u \Vert  _{L^2(\Omega)} / N_i^s
	}  \\
	 & \leq &  C^{N_i} K_i^{N_i}
	 sup_{|J|\leq 2} \Vert Z^J
	 (T^{\tilde{q}})_{\f_i' \psi_i'}) u \Vert  _{L^2(\Omega)} .
\end{eqnarray*}
\end{proposition}
{\bf Proof.}
Is a consequence of the bounds (\ref{1.1.1}) - (\ref{1.1.8}) and the
bounds on $(\f_i\psi_i)^{(s)}$ (Proposition \ref{p2.3.8}).
\hfill$\blacksquare$
\vskip 1.cm
\noindent
Thus we have shown:
\begin{proposition}
	\label{p3.4.5}
\begin{eqnarray}
	 \lefteqn { \sup_{{\scriptstyle |J|\leq 2}\atop { \scriptstyle p\leq N_i}}
	 \Vert Z^I (T^p)_{\f_i \psi_i} u \Vert  _{L^2(\Omega)}  }
	\label{3.4.4} \\
	 & \leq & C^{N_i} K_i^{N_i} \left\{
	 N_i^p \Vert u \Vert  _{L^2(\Omega)} +
	  \sup_{{ \scriptstyle p_1\leq (N_i+2)/2}\atop { \scriptstyle |J|\leq 2}}
	  N_i^{p-p_1} \Vert Z^J (T^{p_1})_{\f_i' \psi_i'} u \Vert  _{L^2(\Omega)}
	 \right\}
	\nonumber \\
	 &  & + C_{{\mit f}}^{N_i} N_i^p K_i^{N_i}.
	\nonumber
\end{eqnarray}
\end{proposition}
Iterating $\log_2 N$ times (\ref{3.4.4}), starting with $i=0$, $p\leq N$, we
obtain
\begin{eqnarray}
	 \lefteqn  {  \sup_{{\scriptstyle p\leq N}\atop {\scriptstyle |J|\leq 2}}
\Vert Z^J (T^p)_{\f_0
\psi_0}
	 u \Vert  _{L^2(\Omega)} }
	\label{3.4.5} \\
	 & \leq &  C^N ( \Pi K_i^{N_i}) N^N \{ \Vert u \Vert  _{L^2(\Omega)} +
	  C_f^N \}
	\nonumber \\
	 & \leq & C^N \ N! ,
	\nonumber
\end{eqnarray}
thus proving Theorem \ref{th0.1} .

\section{Proof of Theorem 2.2}
\label{s4}
\setcounter{equation}{0}
\setcounter{theorem}{0}
\setcounter{proposition}{0}
\setcounter{lemma}{0}
\setcounter{corollary}{0}
\setcounter{definition}{0}
\setcounter{remark}{0}
\renewcommand{\thetheorem}{\thesection.\arabic{theorem}}
\renewcommand{\theproposition}{\thesection.\arabic{proposition}}
\renewcommand{\thelemma}{\thesection.\arabic{lemma}}
\renewcommand{\thedefinition}{\thesection.\arabic{definition}}
\renewcommand{\thecorollary}{\thesection.\arabic{corollary}}
\renewcommand{\theequation}{\thesection.\arabic{equation}}
\renewcommand{\theremark}{\thesection.\arabic{remark}}

To prove Theorem \ref{th0.2} we prove a slightly more general result.

We denote by $Q$ a pseudodifferential operator
of order zero whose symbol has small conic support (to be used
soon to microlocalize), and by
$D$ any of the partial derivatives
${\partial \over {\partial x_j}}, j\leq n.$
And our starting point in terms of estimates is
\be
\label{est:apriori}
\sum_{j=1}^m \normL2{X_jv}^2 +
\normdelta{v}
\leq C\{\Re (Pv,v) +\normL2{v}^2\},
\ee
where $\delta$ is a real number, $0 < \delta \leq 1$ and the operator $P$ has
the form
$$P = \sum_{|I|\leq 2} a_IX_I + ibX_0.$$
\begin{theorem}
If $(x_0,\xi_0)
\notin WF_{s} (f=Pu)$ for some distribution $u,$ and some
$s \geq 1/\delta,$ then
$(x_0,\xi_0) \notin WF_{s} (u).$
\end{theorem}
\vskip0.15 in
\noindent
{\bf Proof.}
Assuming the solution
$u$ belongs to $C^\infty,$ we need to obtain reasonable
estimates for $D^\alpha Qu$ of the form
$$|D^\alpha Qu| \leq C^{|\alpha|+1}\alpha !^s$$ or,
what amounts to the same thing,
$$\normL2{D^\alpha Qu} \leq
C^{|\alpha|+1}N^{s|\alpha|}
\qquad |\alpha| \leq N.$$  Now in view of the
maximality and subellipticity of $P,$ we shall
profit from all aspects of the{\em a priori}
estimate and apply it to $v=D^{\alpha} Qu:$
$$
\sum_1^m \normL2{X_jD^{\alpha} Qu}^2 +
\normdelta{D^{\alpha} Qu} \leq $$
\begin{equation}%
\leq C\{\Re (PD^{\alpha} Qu,D^{\alpha} Qu)_{L^2}
+\normL2{D^{\alpha} Qu}^2\}
\label{eqn:start}
\end{equation}
$$ \leq C\{ \|D^{\alpha} QPu\|_{L^2}^2
+ \|D^{\alpha} Qu\|_{L^2}^2
+ \Re ([P,D^{\alpha} Q]u,D^{\alpha} Qu)_{L^2}\}
$$
The essential work is to commute $P$ with $D^{\alpha} Q,$
and we write this as follows. Schematically writing
$P=(aX_j)(aX_k),$ (the first order term with $iX_0$ will not pose a
problem) we write
\begin{equation}%
\Re ([P,D^{\alpha}Q]u,D^{\alpha} Qu)_{L^2} =
(aX ad_{aX}(D^{\alpha}Q)u,
D^{\alpha} Qu)_{L^2} + \end{equation}%
$$ +
(ad_{aX}^2 (D^{\alpha}Q)u,
D^{\alpha} Qu)_{L^2}$$
$$= E_1 + E_2.$$
If we write $$S^{(\gamma)} =
ad_{aX}^\gamma (S)$$ for any
pseudo-differential operator $S,$ and
$E_1$ then consists of terms with an X free, which may be
integrated by parts (modulo zero order terms). A weighted
Schwarz inequality is used on such terms, with the right hand
member being absorbed on the left hand side of the inequality.
Thus:
\begin{equation}%
|E_1| \leq {1 \over 2} \sum_{j=1}^m
\normL2{X_jD^{\alpha} Qu}^2 +
C\normdelta{\Lambda^{-|\delta |}
(D^{\alpha}Q)^{(1)}u}.
\end{equation}%
In much the same way, $E_2$ may be
estimated
$$ |E_2| \leq {1 \over 2} \normdelta{D^{\alpha}
Qu} +
C\normdelta{\Lambda^{-2|\delta |}
(D^{\alpha}Q)^{(2)}u}.$$
Thus using (\ref{eqn:start}), using the above
estimates for the errors $E_1$ and $E_2,$
and iterating $|\alpha|/|\delta|$ times we have:
$$
\sum_{j=1}^m \normL2{X_jD^{\alpha} Qu}^2 +
\normdelta{D^{\alpha} Qu}\leq $$
$$\leq
C^{|\alpha|} \left \{
\sup_{{\sum j_k \leq {|\alpha|/
|\delta|}}\atop{1\leq j_k\leq2}}
\normL2{\Pi_k(\Lambda^{-|\delta|j_k}
(D^{\alpha}Q)^{(j_k)})Pu}^2
\right .
$$
$$
\left .
+ \sup_{{{|\alpha|\over
|\delta|}\leq\sum j_k \leq {{|\alpha|\over
|\delta|}+1}}\atop{1\leq j_k\leq2}}
\normdelta{\Pi_k(\Lambda^{-|\delta|j_k}
(D^{\alpha}Q)^{(j_k)})u} \right \}$$
What has happened, as can already be seen in the
first step, is that for every derivative that lands on
$D^\alpha Q$
as a bracket, there is a 'gain' of
$|\delta|$ derivatives. Thus either multiple iterations
with the same Q are required, (or one could perhaps
replace the Q with another Q after each such gain,
although this has not seemed to be any simpler). In
any case, estimating the behavior, in $L^2$ or
in $H^{|\delta|}$ of these
very high order commutators of pseudo --
differential operators is very far from simple.
However we have already carried out an extremely precise
 analysis of each such bracket in (\ref{3.3.11}) and what
follows (\ref{3.3.11}). While that analysis is for the
carefully balanced operators required to localize $T^p,$
the present situation requires no such balance and is
thus essentially simpler. We do not carry out the
full details since we would be repeating much of the
analysis following (\ref{3.3.11}).

\section{Appendix}
\label{app}
\renewcommand{\thetheorem}{A.\arabic{theorem}}
\renewcommand{\theproposition}{A.\arabic{proposition}}
\renewcommand{\thelemma}{A.\arabic{lemma}}
\renewcommand{\thedefinition}{A.\arabic{definition}}
\renewcommand{\thecorollary}{A.\arabic{corollary}}
\renewcommand{\theremark}{A.\arabic{remark}}
\renewcommand{\theequation}{A.\arabic{equation}}
\renewcommand{\thesubsection}{A.\arabic{subsection}}
\setcounter{equation}{0}
\setcounter{theorem}{0}
\setcounter{proposition}{0}
\setcounter{lemma}{0}
\setcounter{corollary}{0}
\setcounter{definition}{0}
\setcounter{remark}{0}

In this Appendix we gather some general--purpose results, mostly well
known, which are used and referenced throughout the paper. Sometimes we
shall prove statements a bit more general than needed at a specific point.
\begin{lemma}
	\label{lA.1}
	Let $F\xd \in L^m(\R_x^n)$, $H\xd \in L^{m'}(\R_x^n)$ be pseudo
	differential operators with (full) symbol $f\xx$, $h\xx$ respectively.
	Denote by $\hat{h}(\eta, \xi) = \int e^{-i x \cdot \eta} h\xx dx$ the
	Fourier transform of $h$ with respect to $x$. Then for any
	multi--indices $a$, $b \in \Z_+^n$, $\forall M \in \Z_+$, $\forall w \in
	\czi (\R^n)$ we have
	\begin{eqnarray}
		\lefteqn{D_x^a \left( F\xd \circ H\xd - \left\{ F \circ H \right\}_M
		\xd \right) D_x^b w}
		\label{A.1} \\
		 & = & \underline{C^{|a + b| + n}} \int e^{i x \cdot \xi} \left(
		 \int r_{c_1}(x, \eta, \xi) \dt \eta\right) \hat{w}(\xi) \dt \xi,
		\nonumber
	\end{eqnarray}
where
	\begin{eqnarray*}
		\lefteqn{r_{c_1}(x, \eta, \xi) }\\
		 & = & \sum_{|\varepsilon | = M} \frac{(1 + |\eta |)^{-n-1}}{\varepsilon
		 !} \left(h^{(\varepsilon + b_2 +
c_1)}\right)\hat{\phantom{h}} \ e^{i x
\cdot \eta} \\
   &   & \cdot \left( \int_{0}^{1} f_{(\varepsilon)}^{(a_2)}
   (x, \xi + \rho \eta)
		 (\xi + \rho\eta)^{a_1 + b_1 - c_1} (1 - \rho)^{|c_1 + \varepsilon|}
		 d\rho \right);
	\end{eqnarray*}
	and where
	\begin{equation}
		\left\{F\circ H\right\}_M \xx = \sum_{|\alpha| < M} \frac{1}{\alpha !}
		f_{(\alpha)}\xx h^{(\alpha)}\xx ,
		\label{A.2}
	\end{equation}
	$f_{(\alpha)}^{(\beta)} \xx = \partial_\xi^\alpha D_x^\beta f \xx$ and
	$\varepsilon$, $a_1$, $a_2$, $b_1$, $b_2$, $c_1 \in \Z_+^n$ are
	multi--indices such that
	\begin{displaymath}
		|b_1 + b_2 | = |b| + n + 1 \quad \mbox{{\rm and}} \quad c_1 \leq a_1 +
		b_1.
	\end{displaymath}
\end{lemma}
{\bf Proof.}
By a direct calculation:
\begin{eqnarray*}
\lefteqn{D_x^a \left( F\xd \circ H\xd - \left\{ F \circ H \right\}_M
		\xd \right) D_x^b w}  \\
	 & = & D_x^a \int\int e^{i x\cdot (\xi+\eta)} \left[ f(x, \xi+\eta) -
	 \sum_{|\varepsilon|<M} f_{(\varepsilon)}\xx
	 \frac{\eta^\varepsilon}{\varepsilon !} \right] \hat{h}(\eta, \xi) \xi^b
	 \hat{w}(\xi) \dt \eta \dt \xi  \\
	 & = & \underline{C^{|a+b| +n}} \int\int e^{i x\cdot (\xi+\eta)}
	 (\xi+\eta)^{a_1 + b_1} \sum_{|\varepsilon| = M}\frac{M}{\varepsilon !}
	 (1 + |\eta|)^{-n-1} \\
   &  & \cdot \int_{0}^{1} f_{(\varepsilon)}^{(a_2)}(x,
	 \xi+\rho\eta) (1 - \rho)^{|\varepsilon|} d\rho \left(h^{(\varepsilon +
	 b_2)}\right)\hat{\phantom{h}}(\eta,\xi)\hat{w}(\xi) \dt \eta \dt \xi,
\end{eqnarray*}
where $a_1 + a_2 = a$, $|b_1 + b_2 | \leq b + n + 1$. The latter can be
written as a sum (over $c_1 \leq a_1 + b_1$) of expressions of the form
\begin{displaymath}
	\underline{C^{|a+b| +n}} \int e^{i x \cdot \xi} \left(
		 \int r_{c_1}(x, \eta, \xi) \dt \eta\right) \hat{w}(\xi) \dt \xi,
\end{displaymath}
where
\begin{eqnarray*}
	\lefteqn{r_{c_1}(x, \eta, \xi)}\\
&  = & \sum_{|\varepsilon| = M}
\frac{1}{\varepsilon
	!} (1 + |\eta|)^{-n-1} \left(h^{(\varepsilon +
	 b_2 + c_1 )}\right)\hat{\phantom{h}}(\eta,\xi) e^{i x \cdot \eta} \\
&  &\cdot
	 \left( \int_{0}^{1} f_{(\varepsilon)}^{(a_2)} (x, \xi+\rho\eta) (\xi +
	 \rho\eta)^{a_1 + b_1 - c_1} (1 - \rho)^{|c_1 + \varepsilon|}d\rho
	 \right).
\end{eqnarray*}
Here we wrote $\xi + \eta = \xi + \rho\eta + (1 - \rho)\eta$. And this
completes the proof of the Lemma.
\begin{lemma}
	\label{lA.2}
	Let $\alpha$, $\beta \in \Z_+^n$ be multi--indices. Then
	\begin{equation}
		\partial_x^\alpha x^\beta = \sum_{\delta \leq \alpha , \beta} {\alpha
		\choose \delta} {\beta\choose \delta} \delta ! x^{\beta - \delta}
		\partial_x^{\alpha - \delta};
		\label{A.3}
	\end{equation}
	\begin{equation}
		x^\alpha \partial_x^\beta = \sum_{\delta \leq \alpha, \beta}
		(-1)^{|\delta|} {\alpha  \choose \delta} {\beta\choose \delta} \delta !
		\partial_x^{\beta - \delta} x^{\alpha - \delta}.
		\label{A.4}
	\end{equation}
\end{lemma}
{\bf Proof.}
(\ref{A.3}) is proved by Leibniz' formula:
\begin{displaymath}
	\partial_x^\alpha x^\beta u = \sum_{\delta \leq \alpha, \beta} {\alpha
	\choose \delta} \beta (\beta - 1) \cdots (\beta - \delta + 1) x^{\beta -
	\delta} \partial_x^{\alpha - \delta} u,
\end{displaymath}
whereas (\ref{A.4}) can be proved taking the Fourier transform of
(\ref{A.3}).

\begin{lemma}
	\label{lA.3}
	Let $\alpha_j$, $\beta_j \in \Z_+^n$ be multi--indices, $j = 1, \ldots ,
	r$. Denoting by $\alpha = \sum_{i=1}^{r}\alpha_i$, $\beta =
	\sum_{j=1}^{r}\beta_j$, we have that
	\begin{equation}
		x^{\alpha_1}\partial_x^{\beta_1} x^{\alpha_2}\partial_x^{\beta_2}
		\cdots x^{\alpha_r}\partial_x^{\beta_r}  = \sum_{\delta \leq \alpha,
		\beta}\underline{C^{|\beta|}} \underline{|\alpha|^{|\delta|}} x^{\alpha
		- \delta} \partial_x^{\beta - \delta}.
		\label{A.5}
	\end{equation}
\end{lemma}
{\bf Proof.}
Using (\ref{A.3}) we want to move the last group of $x$'s in (\ref{A.5}),
i.e. $x^{\alpha_r}$, to the left. This yields:
\begin{eqnarray*}
	\lefteqn{x^{\alpha_1}\partial_x^{\beta_1} x^{\alpha_2}\partial_x^{\beta_2}
		\cdots x^{\alpha_r}\partial_x^{\beta_r}}\\
	 & = & \sum_{{\delta_r \leq \alpha_r}\atop {\delta_r \leq
	 \beta_{r-1}}}{\beta_{r-1}\choose \delta_r}
	 \underline{|\alpha_r|^{|\delta_r|}} x^{\alpha_1} \partial_x^{\beta_1}
	 \cdots \partial_x^{\beta_{r - 2}} x^{\alpha_{r-1}+\alpha_r - \delta_r}
	 \partial_x^{\beta_{r-1} + \beta_r - \delta_r}.
\end{eqnarray*}
Iterating this procedure, i.e. moving $x^{\alpha_{r-1} + \alpha_r -
\delta_r}$ to the left and using again (\ref{A.3}) we obtain
\begin{eqnarray*}
	\lefteqn{x^{\alpha_1}\partial_x^{\beta_1} x^{\alpha_2}\partial_x^{\beta_2}
		\cdots x^{\alpha_r}\partial_x^{\beta_r}}\\
	 & = & \sum_{{\scriptstyle \delta_r \leq \alpha_r , \beta_{r-1}}\atop
{\scriptstyle \delta_{r-1} \leq
	 \alpha_r + \alpha_{r-1} - \delta_r , \beta_{r-2}}} {\beta_{r-1}\choose
	 \delta_r}{\beta_{r-2} \choose \delta_{r-1}}
	 \underline{|\alpha_r|^{|\delta_r|}} \ \ \underline{|\alpha_{r-1} +
\alpha_r  - \delta_r |^{|\delta_{r-1}|}} \\
&  & \cdot
x^{\alpha_1}\partial_x^{\beta_1}
\cdots
	 \partial_x^{\beta_{r-3}} x^{\alpha_r + \alpha_{r-1} + \alpha_{r-2} -
	 \delta_r - \delta_{r-1}} \partial_x^{\beta_{r-2} + \beta_{r-1} +
	 \beta_r - \delta_r - \delta_{r-1}}\\
	 & = & \cdots \\
	 & = & \sum_{j=0}^{r} \sum_{{\scriptstyle \delta_{r-j} \leq \beta_{r-j-1}}
\atop
	 {\scriptstyle \delta_{r-j} \leq \sum_{\ell=0}^{j}\left (\alpha_{r-\ell} -
	 \delta_{r-\ell-1}\right ) }} \prod_{s=1}^{r} {\beta_{r-s}\choose
\delta_{r-s+1}}
	 \underline{|\alpha|^{|\delta|}} \  x^{\alpha - \sum_{i=1}^{r}\delta_i}
	 \partial_x^{\beta - \sum_{i=0}^{r}\delta_i},
\end{eqnarray*}
and this proves the assertion since $\prod_{s=1}^{r}{\beta_{r-s}\choose
\delta_{r-s+1}} = \underline{C^{|\beta|}}$.
\begin{proposition}
	\label{pA.4}
	If $X^J$ means $X_1^{j_1} \cdots X_{2k}^{j_{2k}}$, where the $X$'s are
	the same vector fields of Section \ref{s2.1}, then
	\begin{equation}
		X^J = \underline{C^{|J|}} \sum_{J' \leq J} |J|^{|J'|} x^{J - J'}
		\partial_x^{J-J'}.
		\label{A.6}
	\end{equation}
\end{proposition}
{\bf Proof.}
This is an application of Lemma \ref{lA.3}.

\end{document}